\documentclass[hidelinks, 10pt, fleqn]{article}

\usepackage[UKenglish]{babel}
\usepackage[utf8]{inputenc}
\usepackage{amsmath}
\usepackage{amsthm}
\usepackage{amsfonts}
\usepackage{mathtools}
\usepackage{dsfont} 
\usepackage{amssymb}

\usepackage[dvipsnames]{xcolor}
\usepackage{hyperref}
\usepackage{fullpage}
\usepackage{microtype}
\usepackage{newpxmath,newpxtext}
\usepackage{booktabs}
\usepackage{enumitem} 

\usepackage{tikz}
\usetikzlibrary{positioning}
\usetikzlibrary{arrows,backgrounds}
\usetikzlibrary{decorations.pathmorphing}
\usepackage{tikz-cd}
\usepackage{adjustbox}

\usepackage[autostyle,italian=guillemets]{csquotes}
\usepackage[style=numeric,backend=biber]{biblatex}
\tikzset{node distance=2cm, auto}
\tikzcdset{row sep/normal=2.7em,column sep/normal=2.7em}

\allowdisplaybreaks

\setlist[description]{style=multiline,leftmargin=3cm}

\renewcommand\labelenumi{(\roman{enumi})}
\renewcommand\theenumi\labelenumi

\hypersetup{
    colorlinks=true,
    linkcolor=gray,
    filecolor=gray,      
    urlcolor=gray,
    citecolor=gray
    }
\numberwithin{equation}{section}

\newtheorem{definition}{Definition}[section]
\newtheorem{proposition}[definition]{Proposition}
\newtheorem{lemma}[definition]{Lemma}
\newtheorem{theorem}[definition]{Theorem}
\newtheorem{corollary}[definition]{Corollary}

\theoremstyle{definition}
\newtheorem{example}[definition]{Example}
\newtheorem{remark}[definition]{Remark}

\newcommand{\CategoryFont}[1]{\mathsf{#1}}
\newcommand{\FunctorFont}[1]{\mathsf{#1}}
\newcommand{\OperadFont}[1]{\mathscr{#1}}

\renewcommand{\=}{\colon\kern-1ex=}
\renewcommand{\,}{,\dots,}
\renewcommand{\epsilon}{\varepsilon}
\renewcommand{\o}{\circ}

\newcommand{\<}{\langle}
\renewcommand{\>}{\rangle}

\newcommand{\x}{\otimes}

\newcommand{\ntimes}{\times\dots\times}
\renewcommand{\.}{\cdot}

\newcommand{\op}{\mathsf{op}}
\newcommand{\co}{\mathsf{co}}

\newcommand{\lnr}{\mathsf{lnr}}
\renewcommand{\bar}[1]{\overline{{#1}}}
\renewcommand{\^}[1]{^{(#1)}}

\renewcommand{\*}{\ast}
\newcommand{\nto}{\nrightarrow}
\newcommand{\Nat}{\mathsf{Nat}}

\newcommand{\N}{\mathbb{N}}

\renewcommand{\ker}{{\mathsf{ker}\:}}
\newcommand{\id}{\mathsf{id}}
\renewcommand{\S}{\mathbb{S}}
\newcommand{\Lieg}{\mathfrak{g}}

\newcommand{\X}{\mathbb{X}}

\newcommand{\Operad}{\CategoryFont{Operad}}
\newcommand{\Mod}{\CategoryFont{Mod}}
\newcommand{\Alg}{\CategoryFont{Alg}}
\newcommand{\Geom}{\CategoryFont{Geom}}

\newcommand{\OprPair}{\CategoryFont{OprPair}}
\newcommand{\TngPair}{\CategoryFont{TngPair}}
\newcommand{\Cat}{\CategoryFont{Cat}}
\newcommand{\TngCat}{\CategoryFont{TngCat}}

\newcommand{\TTngCat}{\mathbb{T}\CategoryFont{ngCat}}
\newcommand{\trmTngCat}{\CategoryFont{trmTngCat}}

\newcommand{\adjTngCat}{\CategoryFont{adjTngCat}}

\newcommand{\DObj}{\CategoryFont{DObj}}
\newcommand{\DBnd}{\CategoryFont{DBnd}}

\newcommand{\T}{\mathrm{T}}
\newcommand{\TT}{\mathbb{T}}
\newcommand{\FootT}{\mathchoice
{\rotatebox[origin=c]{180}{$\T$}} 
{\rotatebox[origin=c]{180}{$\T$}} 
{\rotatebox[origin=c]{180}{$\T$}} 
{\rotatebox[origin=c]{180}{\scalebox{.7}{$\T$}}} 
}
\newcommand{\FootTT}{\mathchoice
{\rotatebox[origin=c]{180}{$\TT$}} 
{\rotatebox[origin=c]{180}{$\TT$}} 
{\rotatebox[origin=c]{180}{$\TT$}} 
{\rotatebox[origin=c]{180}{\scalebox{.7}{$\TT$}}} 
}
\newcommand{\smallFootT}{{\scriptscriptstyle\FootT}}

\renewcommand{\d}{\mathsf{d}}
\newcommand{\ITC}{\mathfrak{I}}

\newcommand{\Sym}{\FunctorFont{Sym}}

\newcommand{\Kahler}{\Omega}
\newcommand{\Shur}{\FunctorFont{S}}
\newcommand{\Free}{\FunctorFont{Free}}
\newcommand{\Env}{\FunctorFont{Env}}
\newcommand{\Term}{\FunctorFont{Term}}
\newcommand{\Slice}{\FunctorFont{Slice}}

\newcommand{\Init}{\FunctorFont{Init}}

\renewcommand{\P}{\OperadFont{P}}
\newcommand{\OprQ}{\OperadFont{Q}}

\newcommand{\Com}{{\OperadFont{C\!o\!m}}}
\newcommand{\Ass}{{\OperadFont{A\!s\!s}}}
\newcommand{\Lie}{{\OperadFont{L\!i\!e}}}
\newcommand{\Pois}{{\OperadFont{P\!o\!i\!s}}}

\newcommand{\EnvAlg}{\Env}

\title{The differential bundles\\ of the geometric tangent category of an operad}
\author{Marcello Lanfranchi}
\date{}
\makeatother
\makeatletter

\begin{document}

\maketitle

\begin{abstract}
Affine schemes can be understood as objects of the opposite of the category of commutative and unital algebras. Similarly, $\P$-affine schemes can be defined as objects of the opposite of the category of algebras over an operad $\P$. An example is the opposite of the category of associative algebras. The category of operadic schemes of an operad carries a canonical tangent structure. This paper aims to initiate the study of the geometry of operadic affine schemes via this tangent category. For example, we expect the tangent structure over the opposite of the category of associative algebras to describe algebraic non-commutative geometry. In order to initiate such a program, the first step is to classify differential bundles, which are the analogs of vector bundles for differential geometry. In this paper, we prove that the tangent category of affine schemes of the enveloping operad $\P\^A$ over a $\P$-affine scheme $A$ is precisely the slice tangent category over $A$ of $\P$-affine schemes. We are going to employ this result to show that differential bundles over a $\P$-affine scheme $A$ are precisely $A$-modules in the operadic sense.
\end{abstract}

\noindent
\textbf{Acknowledgements.}\quad
{\small
We want to thank Sacha Ikonicoff and Jean-Simon Lemay for the work done together in~\cite{ikonicoff:operadic-algebras-tagent-cats} which led to this paper, and for the informal discussions we had around this topic. We are also thankful to Dorette Pronk and Geoffrey Cruttwell (PhD supervisors) for the discussions, advice, support and precious help during the realization of this article.
}

\tableofcontents

\section{Introduction}
\label{section:introduction}
Cruttwell and Lemay showed that some key geometrical features of affine schemes, in the sense of algebraic geometry, can be captured by defining a suitable tangent structure $\TT$ (cf.~\cite{cruttwell:algebraic-geometry}). A tangent structure $\TT$ over a category $\X$ provides a categorical axiomatization for the tangent bundle functor of differential geometry.
\par Concretely, a tangent structure (cf.~\cite[Definition~2.3]{cockett:tangent-cats}) consists of an endofunctor $\T$ of $\X$ together with a projection $p\colon\T\Rightarrow\id_\X$, a zero section $z\colon\id_\X\Rightarrow\T$ of the projection, a sum morphism $s\colon\T_2\Rightarrow\T$, whose domain $\T_2$ is the pullback of the projection along itself, so that for every object $A\in\X$, $p\colon\T A\to A$ becomes an additive bundle (cf.~\cite[Definition~2.1]{cockett:tangent-cats}), that is a commutative monoid in the slice category over $A$. Moreover, a tangent structure carries two other structures: a vertical lift $l\colon\T\Rightarrow\T^2$, where $\T^2$ denotes the composition of $\T$ with itself, and a canonical flip $c\colon\T^2\Rightarrow\T^2$.
\par The vertical lift defines an abstract version of the Euler vector field and, by satisfying a key universal property (cf.~\cite[Section 2.5]{cockett:tangent-cats}), introduces a notion of linearity for morphisms of differential bundles (cf.~\cite{cockett:differential-bundles}). Moreover, when the tangent category has negatives (cf.~\cite[Section~3.3]{cockett:tangent-cats}) this universal property is also used to equip the set of sections of the projection, i.e. the vector fields, with Lie brackets. Finally, the canonical flip encodes the symmetry of the Hessian matrix.
\par Tangent categories (with negatives) were first introduced by Rosick\'y (\cite{rosicky:tangent-cats}). Recently, the ideas of Rosick\'y were revisited and generalized by Cockett and Cruttwell (\cite{cockett:tangent-cats}) and expanded into a flourishing research program. In the tangent category of affine schemes described by Cruttwell and Lemay, the tangent bundle functor is the functor that maps a commutative algebra $A$ into the symmetric algebra of the $A$-module of K\"{a}hler differentials $\Omega A$ of $A$, i.e. $\T A\=\Sym_A\Omega A$ (cf.~\cite{cruttwell:algebraic-geometry}). One striking result of their paper is the complete classification of differential bundles in this tangent category. Differential bundles, first introduced by Cockett and Cruttwell (cf.~\cite{cockett:differential-bundles}), play the same role as vector bundles in the category of smooth finite-dimensional manifolds for an abstract tangent category (cf.~\cite{macadam:vector-bundles}). Interestingly, Cruttwell and Lemay show that the category of differential bundles and linear morphisms over an affine scheme $A$ is equivalent to the opposite of the category of modules over $A$.
\par The author of this paper together with Sacha Ikonicoff and Jean-Simon Lemay extended the idea of studying the algebraic geometry of affine schemes with tangent categories to a new plethora of contexts. In~\cite{ikonicoff:operadic-algebras-tagent-cats}, they showed that the category of algebras $\Alg_\P$ of a (symmetric) operad $\P$ over the category of $R$-modules (for a commutative and unital ring $R$) comes equipped with a tangent structure. In the following, we refer to this as the \textbf{algebraic tangent structure} of the operad $\P$ which will be denoted by $\FootTT\^\P$, or simply by $\FootTT$ when the operad $\P$ is clear from the context. Moreover, the corresponding tangent category will be denoted as $\Alg(\P)\=(\Alg_\P,\FootTT\^\P)$. In the aforementioned paper, it was proven that every operad comes with a coCartesian differential monad (cf.~\cite[Theorem 4.1.1]{ikonicoff:operadic-algebras-tagent-cats}) and that this tangent category is precisely the tangent category of algebras of this monad.
\par Crucially, $\FootTT\^\P$ admits an adjoint tangent structure (cf.~\cite[Theorem~4.4.4]{ikonicoff:operadic-algebras-tagent-cats}) which makes the opposite of the category of operadic algebras into a tangent category. In the following, we refer to this tangent structure as the \textbf{geometric tangent structure} of the operad $\P$ which will be denoted by $\TT\^\P$, or simply by $\TT$ when the operad $\P$ is clear from the context. This tangent category can be interpreted as the tangent category of affine schemes over the operad $\P$, and will be denoted by $\Geom(\P)\=(\Alg_\P^\op,\TT\^\P)$.
\par To properly appreciate the relevance of this result, notice that before the article~\cite{ikonicoff:operadic-algebras-tagent-cats}, the most revelant available examples of tangent categories were differential geometry, synthetic differential geometry, algebraic geometry, commutative rings etc. In particular, there was no example of non-commutative geometry completely described by tangent category theory. The existence of the geometric tangent category $\Geom(\Ass)$ of the associative operad $\Ass$, whose algebras are associative algebras, proves that tangent categories are suitable to describe a wider variety of geometries, including non-commutative geometry. In Example~\ref{example:com-ass-Geom} we discuss in detail this particular case, with a comparison with the commutative one.\newline
\par In the same paper, differential objects (cf.~\cite[Definition~4.8]{cockett:tangent-cats}) of the geometric tangent category of an operad $\P$ were classified and proved to be in bijective correspondence with left $\P(1)$-modules, where $\P(1)$ denotes the unital and associative ring defined over the first entry of the operad $\P$ and whose unit and multiplication are defined by the unit and the multiplication of the operad.
\par In the same way as the tangent category described by Cruttwell and Lemay captures some key geometrical features of (commutative and unital) affine schemes, we expect the geometric tangent category of an operad $\P$ to capture similar geometrical properties of the affine schemes over $\P$. The goal of this paper is to investigate this assumption by covering the intimate relationship between operads and their corresponding geometric tangent categories. One of the main results of the paper will be the complete classification of differential bundles over operadic affine schemes. We will reinterpret Cruttwell and Lemay's result as a special case of a larger phenomenon: the category of differential bundles and linear morphisms over an operadic affine scheme is equivalent to the opposite of the category of modules of the affine scheme.
\par To prove this, we will first show another key result: the geometric tangent category of the enveloping operad over a $\P$-algebra $A$ is equivalent to the slice tangent category over $A$ of the geometric tangent category of $\P$. The classification of differential bundles will follow directly from this insight: differential bundles are precisely differential objects in the slice tangent category.

\subsection{Outline}
\label{subsection:outline}
The paper is organized as follows. In Section~\ref{section:operadic-affine-schemes}, we first recall the main result of~\cite{ikonicoff:operadic-algebras-tagent-cats} which establishes that every operad $\P$ produces two tangent categories: the algebraic and the geometric tangent categories of $\P$. Once this is established, we show that the operation which takes an operad to its associated tangent categories is functorial (Section~\ref{subsection:functoriality}). In particular, we provide four distinct functors from the category of operads to the category of tangent categories. In Section~\ref{section:slice-tangent-category} we recall the notion of the slice tangent category of a tangent category over an object and we give a new characterization of this construction. In particular, in Section~\ref{subsection:universality-slicing} we show that the operation which takes a tangent pair to its associated slice tangent category extends to a right adjoint of the functor $\Term$, which sends a tangent category with terminal object to the tangent pair formed by the tangent category and its terminal object. The main result of the paper is proved in Section~\ref{section:slicing-operadic-affine-schemes}. We first recall the definition of the enveloping operad of an operadic pair and then prove that the geometric tangent category of the enveloping operad $\P\^A$ of the operadic pair $(\P;A)$ is equivalent to the slice tangent category over the geometric tangent category of the operad $\P$ over $A$. In Section~\ref{subsection:differential-bundles} we employ this result to classify the differential bundles over an operadic affine scheme as modules over the affine scheme. Finally, we dedicate Section~\ref{section:conclusion} to exploring some ideas for future work.

\subsection{Background}
\label{subsection:background}
We assume the reader is comfortable with the theory of symmetric operads over a symmetric monoidal category (see~\cite{loday:operads} for reference), and with fundamental notions of category theory like functors, adjunctions, limits, colimits, pullbacks, pushouts etc. We also assume the reader is knowledgeable about basic notions of tangent category theory (see~\cite{cockett:tangent-cats} for reference). Even if we summarize in the first section the main results of the previous paper, we also recommend reading~\cite{ikonicoff:operadic-algebras-tagent-cats} to fully appreciate the whole story.

\subsection{Notation and naming conventions}
\label{subsection:notation}
We denote by $R$ a fixed commutative and unital ring and by $\Mod_R$ the associated category of left $R$-modules. For an operad $\P$ we refer to a symmetric operad over the symmetric monoidal category $\Mod_R$, where the symmetric monoidal structure is defined by the usual tensor product over $R$, simply denoted by $\x$. The symmetric group that acts over $n$ distinct elements is denoted by $\S_n$. The generators of the free $\P$-algebra over an $R$-module $M$ are denoted by $(\mu;v_1\,v_m)$, where $\mu\in\P(m)$, $v_1\,v_m\in M$. Given $\mu\in\P(m)$ and $\mu_1\in\P(k_1)\,\mu_m\in\P(k_m)$, for positive integers $m,k_1\,k_m$, the operadic composition of $\mu$ with $\mu_1\,\mu_m$ is denoted by $\mu(\mu_1\,\mu_m)$. The unit of the operad $\P$ is denoted by $1_\P$; the monad associated with $\P$ is denoted by $\Shur_\P$, with $\gamma_\P$ for the composition. We denote by $\Operad$ the category of symmetric operads over $\Mod_R$ and their morphisms.
\par The category of $\P$-algebras is denoted by $\Alg_\P$. Given a $\P$-algebra $A$, the action of the abstract $m$-ary operation $\mu\in\P(m)$ over $m$ elements $a_1\,a_m$ of $A$ induced by the structure map of $A$ is denoted by $\mu_A(a_1\,a_m)$ and when $A$ is clear from the context simply by $\mu(a_1\,a_m)$.
\par The category of modules (in the operadic sense) over a $\P$-algebra $A$ is denoted by $\Mod_A\^\P$, or simply by $\Mod_A$ when $\P$ is clear from the context. We will write expressions like $\sum_{k=1}^m\mu(a_1\,x_k\,a_m)$ to denote the sum over the index $k$ of $\mu\.\sigma_k(a_1\,a_{k-1},a_{k+1}\,a_m,x_k)$ where $\sigma_k$ denotes the cylic permutation $(k\quad k+1\quad\dots\quad m)$, where $x_k\in M$, $a_1\,a_{k-1},a_{k+1}\,a_m\in A$ and $M$ is an $A$-module.\newline

\par Given a tangent category $(\X,\TT)$, we denote the tangent bundle functor $\T$ by using the same letter as used for the tangent structure. For the projection, the zero morphism, the sum morphism, the lift, the canonical flip, and the negation (in case of a tangent category with negatives) we will use the letters $p\^{\T},z\^{\T},s\^{\T},l\^{\T},c\^{\T}$ and $n\^{\T}$, respectively. When the tangent structure is clear from the context, we will simplify the notation by omitting the superscript $\^{\T}$.
\par Morphisms of tangent categories come in different flavours. We need to distinguish among them therefore we introduce the following convention. Given two tangent categories $(\X,\TT)$ and $({\X'},\TT')$, we refer to a \textbf{lax tangent morphism} $(F,\alpha)\colon(\X,\TT)\to({\X'},\TT')$ as a functor $F\colon\X\to{\X'}$ together with a natural transformation $\alpha\colon F\o\T\Rightarrow\T'\o F$ compatible with the tangent structures (cf.~\cite[Definition~2.7]{cockett:tangent-cats}). We refer to $\alpha$ as the \textbf{lax distributive law} of the morphism.
\par By a \textbf{colax tangent morphism} $(G,\beta)\colon(\X,\TT)\nto({\X'},\TT')$ we mean a functor $G\colon\X\to{\X'}$ together with a natural transformation $\beta\colon\T'\o G\Rightarrow G\o\T$ compatible with the tangent structures (the compatibilities are similar to the ones of a lax tangent morphism, where the distributive law goes in the opposite direction). We refer to $\beta$ as the \textbf{colax distributive law} of the morphism. We also adopt the notation $\nto$ to denote colax tangent morphisms.
\par By a \textbf{strong tangent morphism} we mean a lax tangent morphism where the distributive law is an isomorphism. Notice that the underlying functor of a strong tangent morphism together with the inverse of the lax distributive law defines a colax tangent morphism. Finally, by a \textbf{strict tangent morphism} we refer to a strong tangent morphism whose distributive law is the identity. Since, in this case, the distributive law is trivial, we will omit it completely in the notation and simply refer to the functor as the strict tangent morphism.
\par We denote by $\TngCat$ the category of tangent categories and lax tangent morphisms. When required, we abuse notation and denote by $\TngCat$ the $2$-category with the same objects and $1$-morphisms and whose $2$-morphisms are natural transformations compatible with the lax distributive laws. Similarly, we denote by $\TngCat_\cong$ the category of tangent categories and strong tangent morphisms, and finally, by $\TngCat_=$ the category of tangent categories and strict tangent morphisms.

\par Adopting the same naming convention used in~\cite{ikonicoff:operadic-algebras-tagent-cats}, a category $\X$ is called \textbf{semi-additive} if $\X$ has finite biproducts, which means that it admits finite products, finite coproducts and the canonical morphism between $n$ products and $n$ coproducts is an isomorphism. We denote by $\oplus$ the biproducts of $\X$. In particular, in $\Mod_R$, given two $R$-modules $X$ and $Y$, we denote the elements of $X\oplus Y$ as pairs $(x,y)$ for each $x\in X$ and $y\in Y$. In such a category, the empty biproduct is denoted by $0$ and is the zero object, which is an object that is both initial and terminal. Note that for a category to be semi-additive is equivalent to being enriched over the category of commutative monoids.

\par A semi-additive category $\X$ comes equipped with a canonical tangent structure $\FootTT$ whose tangent bundle functor $\FootT$ is the diagonal functor $\FootT X=X\oplus X$, the projection is the projection on the first coordinate, i.e. $p=\pi_1\colon X\oplus X\to X$, the zero morphism is the injection in the first coordinate, i.e. $z=\iota_1\colon X\to X\oplus X$, the $n$-fold pullback of the projection along itself is (isomorphic to) the $n+1$ tuple $\FootT_n X=X\oplus X\oplus\dots\oplus X$, the sum morphism is the identity in the first coordinate and the sum on the second and the third, i.e. $s\colon X\oplus X\oplus X\xrightarrow{\id_X\oplus+}X\oplus X$; the vertical lift maps the first coordinate to the first one and the second coordinate to the fourth one, i.e. $l\colon X\oplus X\xrightarrow{\<\pi_1\o\iota_1,\pi_4\o\iota_2\>}X\oplus X\oplus X\oplus X$; the canonical flip flips the internal coordinates, i.e. $c\colon X\oplus X\oplus X\oplus X\xrightarrow{\id_X\oplus\tau\oplus\id_X}X\oplus X\oplus X\oplus X$, where $\tau=\<\pi_2\o\iota_1,\pi_1\o\iota_2\>\colon X\oplus Y\to Y\oplus X$; finally, if $\X$ is \textbf{additive}, i.e. $\CategoryFont{Ab}$-enriched, then the negation morphism is the identity on the first coordinate and the negation on the second, i.e. $n\colon X\oplus X\xrightarrow{\id_X\oplus-}X\oplus X$. In this paper, we refer to the tangent structure $\FootTT$ induced by additivity over $\Mod_R$ as the \textbf{canonical tangent structure} and to $(\Mod_R,\FootTT)$ as the \textbf{canonical tangent category}.\newline

\par For two composable morphisms $f\colon A\to B$ and $g\colon B\to C$ of a category $\X$, we denote by $g\o f$ their composition. We will also often use the diagrammatic notation, i.e. $fg\=g\o f$. For functors, we adopt a similar notation with a single variation: when an object $X\in\X$ is specified, we denote by $GFX$ the object $(G\o F)(X)$ and similarly for morphisms. An adjunction between two functors $F\colon\X\to{\X'}$ and $G\colon{\X'}\to\X$ with unit $\eta$ and counit $\epsilon$ is denoted by $(\eta,\epsilon)\colon F\dashv G$. A similar notation will be adopted for conjunctions in the context of double categories.

\section{The geometry of affine schemes over an operad}
\label{section:operadic-affine-schemes}
In~\cite{ikonicoff:operadic-algebras-tagent-cats}, the author of this paper, Sacha Ikonicoff, and Jean-Simon Lemay showed that every operad provides a tangent structure over the category of operadic algebras (\cite[Theorem~4.3.3]{ikonicoff:operadic-algebras-tagent-cats}) as well as a tangent structure over the opposite of the same category (\cite[Theorem~4.4.4]{ikonicoff:operadic-algebras-tagent-cats}). Since this is the starting point for this paper, we dedicate this section to recall this construction.
\par Concretely, the tangent structure $\FootTT\^\P$ over the category of $\P$-algebras is defined as follows:
\begin{description}
\item[tangent bundle functor] The tangent bundle functor $\FootT\colon\Alg_\P\to\Alg_\P$ maps every $\P$-algebra $A$ to the semi-direct product $A\ltimes A$, which is the $\P$-algebra over the $R$-module $A\times A$ and with structure map defined as follows:
\begin{align*}
&\mu\left((a_1,b_1)\,(a_m,b_m)\right)\=\left(\mu(a_1\,a_m),\sum_{k=1}^m\mu(a_1\,b_k\,a_m)\right)
\end{align*}
\item[projection] The projection $p\^{\smallFootT}\colon\FootT A\to A$ projects along the first component, that is:
\begin{align*}
&p\^{\smallFootT}(a,b)\=a
\end{align*}
\item[$n$-fold pullbacks] The $n$-fold pullback along the projection of the tangent bundle functor $\FootT_n\colon\Alg_\P\to\Alg_\P$ maps every $\P$-algebra into the semi-direct product $A\ltimes(A\ntimes A)$. Moreover, the $k$-th projection $\pi_k\colon\FootT_n A\to\FootT A$ is defined as follows:
\begin{align*}
&\pi_k(a;b_1\,b_n)\=(a,b_k)
\end{align*}
\item[zero morphism] The zero morphism $z\^{\smallFootT}\colon A\to\FootT A$ injects into the first component, that is:
\begin{align*}
&z\^{\smallFootT}(a)\=(a,0)
\end{align*}
\item[sum morphism] The sum morphism $s\^{\smallFootT}\colon\FootT_2A\to\FootT A$ is defined by:
\begin{align*}
&s\^{\smallFootT}(a;b_1,b_2)\=(a;b_1+b_2)
\end{align*}
\item[vertical lift] The vertical lift $l\^{\smallFootT}\colon\FootT A\to\FootT^2A$ is defined by:
\begin{align*}
&l\^{\smallFootT}(a,b)\=(a,0,0,b)
\end{align*}
\item[canonical flip] The canonical flip $c\^{\smallFootT}\colon\FootT^2A\to\FootT^2A$ is defined by:
\begin{align*}
&c\^{\smallFootT}(a_1,b_1,a_2,b_2)\=(a_1,a_2,b_1,b_2)
\end{align*}
\item[negation] The negation morphism $n\^{\smallFootT}\colon\FootT A\to\FootT A$ is defined by:
\begin{align*}
&n\^{\smallFootT}(a,b)\=(a,-b)
\end{align*}
\end{description}
On the other hand, the tangent structure $\TT\^\P$ over the opposite of the category of the category of $\P$-algebras is defined as follows:
\begin{description}
\item[tangent bundle functor] The tangent bundle functor $\T\colon\Alg_\P^\op\to\Alg_\P^\op$ maps a $\P$-algebra $A$ to the $\P$-algebra $\Free_A\Kahler A$, where $\Free_A\colon\Mod_A\to\Alg_\P$ is the functor that maps a $A$-module $M$ to the free $\P$-algebra under $A$ (cf.~\cite[Proposition~4.4.1]{ikonicoff:operadic-algebras-tagent-cats}) and $\Kahler A$ is the module of K\"ahler differentials of $A$. Concretely, $\T A$ is the $P$-algebra generated by all elements $a$ of $A$ and symbols $\d\^\P a$, for each $a\in A$ such that the following relations are fulfilled:
\begin{align*}
&\mu_{\T A}(a_1\,a_m)=\mu_A(a_1\,a_m)\\
&\d\^\P(ra+sb)=r\d\^\P a+s\d\^\P b\\
&\d\^\P\left(\mu(a_1\,a_m)\right)=\sum_{k=1}^m\mu(a_1\,\d\^\P a_k\,a_m)
\end{align*}
for every $r,s\in R$ and $a,b,a_1\,a_m\in A$. In the following, we will omit the superscript $\^\P$ in $\d\^\P$ whenever the operad $\P$ is clear from the context. 
\item[projection] The projection, regarded as an $\Alg_\P$-morphism, $p\^\T\colon A\to\T A$ injects $a\in A$ into $a\in\T A$.
\item[$n$-fold pullbacks] The $n$-fold pushout (in $\Alg_\P$) along the projection of the tangent bundle functor $\T_n\colon\Alg_\P^\op\to\Alg_\P^\op$ is the $\P$-algebra generated by all the elements $a$ of $A$ and by symbols $\d_1a,\d_2a\,\d_na$, for each $a\in A$, such that the following relations are fulfilled:
\begin{align*}
&\mu_{\T_nA}(a_1\,a_m)=\mu_A(a_1\,a_m)\\
&\d_i(ra+sb)=r\d_ia+s\d_ib\\
&\d_i\left(\mu(a_1\,a_m)\right)=\sum_{k=1}^m\mu(a_1\,\d_ia_k\,a_m)
\end{align*}
for every $r,s\in R$, $a,b,a_1\,a_m\in A$, and for every $i=1\,n$. Moreover, the injections $\iota_k\colon\T A\to\T_nA$ map each $a$ to $a$ and $\d a$ to $\d_k a$, for every $k=1\,n$.
\item[zero morphism] The zero morphism, regarded as a $\Alg_\P$-morphism, $z\^\T\colon\T A\to A$ projects each $a$ to itself $a$ and each $\d a$ to $0$.
\item[sum morphism] The sum morphism, regarded as a $\Alg_\P$-morphism, $s\^\T\colon\T A\to\T_2A$ maps each $a$ to $a$ and each $\d a$ into $\d_1a+\d_2a$.
\item[vertical lift] The vertical lift, regarded as a $\Alg_\P$-morphism, $l\^\T\colon\T^2A\to\T A$ maps each $a\in A$ to $a$, $\d a$ and $\d'a$ to $0$ and $\d'\d a$ to $\d a$.
\item[canonical flip] The canonical flip, regarded as a $\Alg_\P$-morphism, $c\^\T\colon\T^2A\to\T^2A$ maps each $a$ to $a$, $\d a$ to $\d'a$, $\d'a$ to $\d a$ and $\d'\d a$ to $\d'\d a$.
\item[negation] The negation morphism, regarded as a $\Alg_\P$-morphism, $n\^\T\colon\T A\to\T A$ maps each $a$ to $a$ and each $\d a$ to $-\d a$.
\end{description}
In the following, given an operad $\P$, we refer to $\Alg(\P)\=(\Alg_\P,\FootTT\^\P)$ as the \textbf{algebraic tangent category} of $\P$ and to $\Geom(\P)\=(\Alg_\P^\op,\TT\^\P)$ as the \textbf{geometric tangent category} of $\P$.

\subsection{The functoriality of the algebraic and the geometric tangent categories}
\label{subsection:functoriality}
So far we recapped the main result of \cite{ikonicoff:operadic-algebras-tagent-cats}: every operad produces two distinct tangent categories, $\Alg(\P)$ and $\Geom(\P)$. In this section, we explore the relationship between morphisms of operads and the corresponding morphisms of tangent categories; we will also show that this operation is functorial.
\par First, we briefly recall that a morphism of operads $\varphi\colon\P\to\OprQ$ is a sequence of $R$-linear morphisms $\{\varphi_n\colon\P(n)\to\OprQ(n)\}_{n\in\N}$, compatible with the operadic structures, that is, given $\mu\in\P(m),\mu_1\in\P(k_1)\,\mu_m\in\P(k_m)$:
\begin{align*}
&\varphi_1(1_\P)=1_{\OprQ}\\
&\varphi_M(\mu(\mu_1\,\mu_m))=\varphi_m(\mu)(\varphi_{k_1}(\mu_1)\,\varphi_{k_m}(\mu_m))
\end{align*}
where $M\=k_1+\dots+k_m$. For the sake of simplicity, in the following, we will omit the index and simply denote by $\varphi$ any of the morphisms in the sequence. A morphism of operads induces a forgetful functor $\varphi^\*\colon\Alg_\OprQ\to\Alg_\P$, which sends a $\OprQ$-algebra $B$ into the $\P$-algebra $\varphi^\*B$ over the $R$-module underlying $B$ and with structure map defined by:
\begin{align*}
&\mu_{\varphi^\*B}(b_1\,b_m)\=(\varphi(\mu))_B(b_1\,b_m)
\end{align*}
The functor $\varphi^\*$ admits a left adjoint $\varphi_!\colon\Alg_\P\to\Alg_\OprQ$, which sends each $\P$-algebra $A$ to the $\OprQ$-algebra $\varphi_!A$ obtained by identifying the two structure maps induced by the operadic composition and by the structure map of $A$ over the free $\OprQ$-algebra over the underlying $R$-module of $A$. Concretely, $\varphi_!A$ can be understood as the coequalizer:
\begin{equation*}
\begin{tikzcd}
{\Shur_\OprQ\Shur_\P A} & {\Shur_\OprQ\Shur_\OprQ A} & {\Shur_\OprQ A} & {\varphi_!A}
\arrow[dashed, from=1-3, to=1-4]
\arrow["{\Shur_\OprQ\Shur_\varphi A}", from=1-1, to=1-2]
\arrow["{\gamma_\OprQ}", from=1-2, to=1-3]
\arrow["{\Shur_\OprQ\theta}"', bend right, from=1-1, to=1-3]
\end{tikzcd}
\end{equation*}
where $\theta$ is the structure map of $A$.
\par As already mentioned in the introduction, the existence of the algebraic tangent category $\Alg(\P)$ of an operad $\P$ is a consequence of the fact that the monad $\Shur_\P$ associated to $\P$ carries a differential combinator $\partial_\P$, so that $\Shur_\P$ becomes a coCartesian differential monad over $\Mod_R$ (see~\cite[Section~4.1]{ikonicoff:operadic-algebras-tagent-cats} for details). As shown by~\cite[Corollary~3.2.6]{ikonicoff:operadic-algebras-tagent-cats}, over a semi-additive category $\X$ there is a bijective correspondence between coCartesian differential monads over $\X$ and tangent monads over the tangent category $(\X,\FootTT)$, where $\FootTT$ is defined by the existence of biproducts in $\X$ (cf.~\cite[Lemma~3.1.1]{ikonicoff:operadic-algebras-tagent-cats}).
\par We recall that a tangent monad, first introduced in~\cite[Definition~19]{cockett:tangent-monads}, is a monad in the $2$-category $\TngCat$ of tangent categories, lax tangent morphisms, and tangent natural transformations, which are natural transformations compatible in an obvious way with the distributive laws. We also recall that the distributive law associated to a tangent monad lifts the tangent structure over the base tangent category to the category of algebras of the monad. Concretely, the tangent bundle functor $\FootT\^S\colon\Alg_S\to\Alg_S$ over the category of algebras of a tangent monad $(S,\alpha)$ over the canonical tangent category $(\Mod_R,\FootTT)$ sends an $S$-algebra $A$ with structure map $\theta\colon SA\to A$ into the $S$-algebra $\FootT A$ with structure map $S\FootT A\xrightarrow{\alpha}\FootT SA\xrightarrow{\smallFootT\theta}\FootT A$, where $\alpha\colon S\o\FootT\Rightarrow\FootT\o S$ is the lax distributive law of $S$.
\par This is precisely the origin of the tangent structure of $\Alg(\P)$, which is lifted from the canonical tangent structure on $\Mod_R$. On the other hand, the tangent structure of $\Geom(\P)$ is the adjoint tangent structure of the algebraic one (see~\cite[Section~4.4]{ikonicoff:operadic-algebras-tagent-cats} for details). Concretely, this means that the tangent bundle functor $\T$, regarded as an endofunctor over $\Alg_\P$, is the left adjoint of $\FootT$ and that the projection, the zero morphism, the sum morphism, the vertical lift, the canonical flip and the negation of $\TT$ are the mates of the corresponding natural transformations of $\FootTT$ along the adjunction $\T\dashv\FootT$.
\par The intimate connection between operads and tangent monads plays a crucial role in understanding the relationship between morphisms of operads and corresponding morphisms of tangent categories. It is not hard to see that a morphism of operads $\varphi\colon\P\to\OprQ$ induces a morphism of the corresponding tangent monads $\varphi\colon(\Shur_\P,\alpha_\P)\to(\Shur_\OprQ,\alpha_\OprQ)$, where we recall that the distributive law $\alpha_\P\colon\Shur_\P\o\FootT\Rightarrow\FootT\o\Shur_\P$ associated to an operad $\P$ is the natural transformation:
\begin{align*}
&\alpha_\P\left(\mu;(x_1,y_1)\,(x_m,y_m)\right)=\left((\mu;x_1\,x_m),\sum_{k=1}^m(\mu;x_1\,y_k\,x_m)\right)
\end{align*}
In this context, a morphism of tangent monads $\varphi\colon(S,\alpha)\to(W,\beta)$ over $(\Mod_R,\FootTT)$ consists of a natural transformation $\varphi\colon S\Rightarrow W$, compatible with the lax distributive laws $\alpha$ and $\beta$, that is:
\begin{equation}
\label{equation:morphism-tangent-monads}
\begin{tikzcd}
S\o\FootT & W\o\FootT \\
{\FootT\o S} & {\FootT\o W}
\arrow["{\varphi_{\smallFootT}}", from=1-1, to=1-2]
\arrow["\FootT\varphi"', from=2-1, to=2-2]
\arrow["\alpha"', from=1-1, to=2-1]
\arrow["\beta", from=1-2, to=2-2]
\end{tikzcd}
\end{equation}
Moreover, since the tangent structure $\FootTT\^S$ over the category of algebras $\Alg_S$ of a tangent monad $(S,\alpha)$ is lifted along the distributive law $\alpha$ from the base tangent category $(\Mod_R,\FootTT)$, a morphism of tangent monads $\varphi\colon(S,\alpha)\to(W,\beta)$ induces a strict tangent morphism $\varphi^\*\colon(\Alg_W,\FootTT\^W)\to(\Alg_S,\FootTT\^S)$, whose underlying functor is the forgetful functor which sends a $W$-algebra $B$ with structure map $\psi\colon WB\to B$ to the $S$-algebra $B$ with structure map $SB\xrightarrow{\varphi}WB\xrightarrow{\psi}B$. To see this, take a $W$-algebra $B$ with structure map $\psi\colon WB\to B$. So, $\varphi^\*\FootT\^WB$ is the $S$-algebra $B$ with structure map:
\begin{align*}
&S\FootT B\xrightarrow{\varphi_{\smallFootT}}W\FootT B\xrightarrow{\beta}\FootT WB\xrightarrow{\smallFootT\psi}\FootT B
\end{align*}
On the other hand, $\FootT\^S\varphi^\*B$ is the $S$-algebra $B$ with structure map:
\begin{align*}
&S\FootT B\xrightarrow{\alpha}\FootT SB\xrightarrow{\smallFootT\varphi}\FootT WB\xrightarrow{\smallFootT\beta}\FootT B
\end{align*}
Thanks to Equation~\eqref{equation:morphism-tangent-monads} and to the naturality of $\varphi$, $\varphi^\*\FootT\^WB$ is precisely $\FootT\^S\varphi^\*B$.
\par By putting together that morphisms of operads induce morphisms of tangent monads and that morphisms of tangent monads induce strict tangent morphisms of the corresponding tangent categories, we find that:

\begin{proposition}
\label{proposition:functoriality-Alg-star}
The operation which takes an operad to its algebraic tangent category extends to a functor $\Alg^\*\colon\Operad^\op\to\TngCat_=$ which sends a morphism of operads $\varphi\colon\P\to\OprQ$ to the strict tangent morphism $\varphi^\*\colon\Alg(\OprQ)\to\Alg(\P)$.
\end{proposition}

As previously recalled, a morphism of operads $\varphi\colon\P\to\OprQ$ induces a left adjoint $\varphi_!\colon\Alg_\P\to\Alg_\OprQ$. Given a tangent morphism $(G,\beta)\colon({\X'},\TT')\to(\X,\TT)$ between two tangent categories whose underlying functor $G\colon{\X'}\to\X$ admits a left adjoint $F\colon\X\to{\X'}$, it is natural to ask whether or not the functor $F$ inherits from $(G,\beta)$ a distributive law $\alpha$ which makes $(F,\alpha)$ into a new tangent morphism.
\par It turns out that this works only if $(G,\beta)$ is a colax tangent morphism. In that case, $F$ becomes a lax tangent morphism. This interesting role played by colax tangent morphisms is better contextualized within the settings of double categories. Heuristically, a double category is a collection of objects together with two classes of morphisms, called horizontal and vertical morphisms, denoted by $\to$ and the second ones by $\nto$, respectively, and a collection of double cells, that are squares:
\begin{equation*}
\begin{tikzcd}
\bullet & \bullet \\
\bullet & \bullet
\arrow[""{name=0, anchor=center, inner sep=0}, from=1-1, to=1-2]
\arrow[""{name=1, anchor=center, inner sep=0}, from=2-1, to=2-2]
\arrow["\shortmid"{marking}, from=1-2, to=2-2]
\arrow["\shortmid"{marking}, from=1-1, to=2-1]
\arrow["\theta"{description}, draw=none, from=0, to=1]
\end{tikzcd}
\end{equation*}
which can be composed horizontally and vertically. We invite the interested reader to consult~\cite{ehresmann:double-cats} for more details on double categories. Notice that double categories can also be characterized as internal categories in the $2$-category of categories.

\begin{proposition}
\label{proposition:double-category-tangent-categories}
Tangent categories can be organized into a double category $\TTngCat$ whose horizontal morphisms are lax tangent morphisms, vertical morphisms are colax tangent morphisms and double cells:
\begin{equation*}
\begin{tikzcd}
{(\X_1,\TT_1)} & {(\X'_1,\TT'_1)} \\
{(\X_2,\TT_2)} & {(\X'_2,\TT'_2)}
\arrow["{(G,\beta)}"', "\shortmid"{marking}, from=1-1, to=2-1]
\arrow["{(G',\beta')}", "\shortmid"{marking}, from=1-2, to=2-2]
\arrow["{(F_1,\alpha_1)}", from=1-1, to=1-2]
\arrow["{(F_2,\alpha_2)}"', from=2-1, to=2-2]
\arrow["\varphi"{description}, draw=none, from=2-1, to=1-2]
\end{tikzcd}
\end{equation*}
are \textbf{tangent double cells}, which are natural transformations $\varphi\colon F_2\o G\Rightarrow G'\o F_1$, fulfilling the commutativity of the following diagram:
\begin{equation*}
\begin{tikzcd}
{F_2\o\T_2\o G} & {F_2\o G\o\T_1} & {G'\o F_1\o\T_1} \\
{\T'_2\o F_2\o G} & {\T'_2\o G'\o F_1} & {G'\o\T'_1\o F_1}
\arrow["{(\alpha_2)_{G}}"', from=1-1, to=2-1]
\arrow["{F_2\beta}", from=1-1, to=1-2]
\arrow["{\varphi_\T}", from=1-2, to=1-3]
\arrow["{G'\alpha_1}", from=1-3, to=2-3]
\arrow["{\T_2'\varphi}"', from=2-1, to=2-2]
\arrow["{\beta'_{F_1}}"', from=2-2, to=2-3]
\end{tikzcd}
\end{equation*}
\begin{proof}
The proof that $\TTngCat$ is a double category is straightforward but tedious, thus is left to the reader.
\end{proof}
\end{proposition}

Proposition~\ref{proposition:double-category-tangent-categories} shows that tangent categories can be organized into a double category. Conjunctions in this double category play a fundamental role in our story. Intuitively speaking, a conjunction in an arbitrary double category is the analog of an adjunction of $1$-morphisms in a $2$-category. Concretely, a conjunction consists of a vertical morphism $G\colon{\X'}\to\X$ together with a horizontal morphism $F\colon\X\to{\X'}$ and two double cells $\eta$ and $\epsilon$
\begin{equation*}
\begin{tikzcd}
\X & {\X'} \\
\X & \X
\arrow["G", "\shortmid"{marking}, from=1-2, to=2-2]
\arrow["F", from=1-1, to=1-2]
\arrow["\shortmid"{marking}, Rightarrow, no head, from=1-1, to=2-1]
\arrow[Rightarrow, no head, from=2-1, to=2-2]
\arrow["\eta"{description}, Rightarrow, from=2-1, to=1-2]
\end{tikzcd}\hfill
\begin{tikzcd}
{\X'} & {\X'} \\
\X & {\X'}
\arrow[Rightarrow, no head, from=1-1, to=1-2]
\arrow["\shortmid"{marking}, Rightarrow, no head, from=1-2, to=2-2]
\arrow["G"', "\shortmid"{marking}, from=1-1, to=2-1]
\arrow["F"', from=2-1, to=2-2]
\arrow["\epsilon"{description}, Rightarrow, from=2-1, to=1-2]
\end{tikzcd}
\end{equation*}
fulfilling the triangle identities.

\begin{proposition}
\label{proposition:conjuctions-tangent-morphisms}
If the underlying functor $G$ of a colax tangent morphism $(G,\beta)\colon({\X'},\TT')\nto(\X,\TT)$ is the right adjoint in a functorial adjunction $(\eta,\epsilon)\colon F\dashv G$, then the left adjoint $F$ becomes a lax tangent morphism with the lax distributive law defined as the mate of $\beta$ along the adjunction, that is:
\begin{align}
\label{equation:distributive-law-conjoints}
&\alpha\colon F\o\T\xrightarrow{F\T\eta}F\o\T\o G\o F\xrightarrow{F\beta_F}F\o G\o\T'\o F\xrightarrow{\epsilon_{\T' F}}\T'\o F
\end{align}
In particular, $(\eta,\epsilon)\colon(F,\alpha)\dashv(G,\beta)$ forms a conjunction in the double category $\TTngCat$. Finally, also the opposite holds: any conjunction in $\TTngCat$ is of the form $(\eta,\epsilon)\colon(F,\alpha)\dashv(G,\beta)$ where $\alpha$ is defined as in Equation~\eqref{equation:distributive-law-conjoints} and $(\eta,\epsilon)\colon F\dashv G$ is a functorial adjunction.
\begin{proof}
Let's start by proving that $(F,\alpha)$ is a lax tangent morphism. The first step is to show that $\alpha$ is compatible with the projections, i.e. $\alpha p\^{\T'}_F=Fp\^\T$, where $p\^{\T'}$ denotes the projection of the tangent structure $\TT'$ and $p\^\T$ the projection of $\TT$. We will adopt a similar notation for the other natural transformations of the tangent structures. This amounts to showing the commutativity of the following diagram:
\begin{equation*}
\begin{tikzcd}
F\o\T & {F\o\T\o G\o F} & {F\o G\o\T'\o F} & {\T'\o F} \\
& F\o G\o F & F\o G\o F \\
F &&& F
\arrow["F\eta"', from=3-1, to=2-2]
\arrow[""{name=0, anchor=center, inner sep=0}, Rightarrow, no head, from=3-1, to=3-4]
\arrow["Fp\^\T"', from=1-1, to=3-1]
\arrow["{p\^{\T'}_F}", from=1-4, to=3-4]
\arrow["(F\o\T)\eta", from=1-1, to=1-2]
\arrow[""{name=1, anchor=center, inner sep=0}, "{F\beta_F}", from=1-2, to=1-3]
\arrow["{\epsilon_{\T'\o F}}", from=1-3, to=1-4]
\arrow["{Fp\^\T_{G\o F}}"', from=1-2, to=2-2]
\arrow["{(F\o G)p\^{\T'}_F}", from=1-3, to=2-3]
\arrow["{\epsilon_F}"', from=2-3, to=3-4]
\arrow[""{name=2, anchor=center, inner sep=0}, Rightarrow, no head, from=2-2, to=2-3]
\arrow["\Nat"', draw=none, from=1-1, to=2-2]
\arrow["\Nat", draw=none, from=1-4, to=2-3]
\arrow["\Delta"{description}, draw=none, from=2, to=0]
\arrow["{(\beta;p\^\T,p\^{\T'})}"{description}, draw=none, from=1, to=2]
\end{tikzcd}
\end{equation*}
To express the commutativity of the diagrams that compose the whole diagram we adopted the following convention: with $\Nat$ we denoted commutativity by naturality, by $(\beta;p\^\T,p\^{\T'})$ we denoted the compatibility between $\beta$ and the projections, and $\Delta$ indicates the triangle identities between the unit and the counit of the adjunction. In the following, we adopt a similar notation.
\par The second step is to prove the compatibility with the zero morphisms. This amounts to showing that $Fz\^\T\alpha=z\^{\T'}_F$, i.e.:
\begin{equation*}
\begin{tikzcd}
F\o\T & {F\o\T\o G\o F} & {F\o G\o\T'\o F} & {\T'\o F} \\
& F\o G\o F & F\o G\o F \\
F &&& F
\arrow["F\eta"', from=3-1, to=2-2]
\arrow[""{name=0, anchor=center, inner sep=0}, Rightarrow, no head, from=3-1, to=3-4]
\arrow["Fz\^\T", from=3-1, to=1-1]
\arrow["{z\^{\T'}_F}"', from=3-4, to=1-4]
\arrow["(F\o\T)\eta", from=1-1, to=1-2]
\arrow[""{name=1, anchor=center, inner sep=0}, "{F\beta_F}", from=1-2, to=1-3]
\arrow["{\epsilon_{\T'\o F}}", from=1-3, to=1-4]
\arrow["{Fz\^\T_{G\o F}}", from=2-2, to=1-2]
\arrow["{(F\o G)z\^{\T'}_F}"', from=2-3, to=1-3]
\arrow["{\epsilon_F}"', from=2-3, to=3-4]
\arrow[""{name=2, anchor=center, inner sep=0}, Rightarrow, no head, from=2-2, to=2-3]
\arrow["\Nat"', draw=none, from=1-1, to=2-2]
\arrow["\Nat", draw=none, from=1-4, to=2-3]
\arrow["\Delta"{description}, draw=none, from=2, to=0]
\arrow["{(\beta;z\^\T,z\^{\T'})}"{description}, draw=none, from=1, to=2]
\end{tikzcd}
\end{equation*}
Let's show the compatibility with the sum morphism, which is $(\alpha)_2s\^{\T'}_F=Fs\^\T\alpha$:
\begin{equation*}
\begin{tikzcd}
{F\o\T_2} & {F\o\T_2\o G_2\o F_2} & {F\o G_2\o\T'_2\o F_2} & {\T'_2\o F_2} \\
F\o\T & {F\o\T\o G\o F} & {F\o G\o\T'\o F} & F
\arrow["Fs\^\T"', from=1-1, to=2-1]
\arrow["{s\^{\T'}_F}", from=1-4, to=2-4]
\arrow[""{name=0, anchor=center, inner sep=0}, "{(F\o\T_2)\eta_2}", from=1-1, to=1-2]
\arrow[""{name=1, anchor=center, inner sep=0}, "{F(\beta_2)_{F_2}}", from=1-2, to=1-3]
\arrow[""{name=2, anchor=center, inner sep=0}, "{(\epsilon_2)_{\T'_2\o F_2}}", from=1-3, to=1-4]
\arrow[""{name=3, anchor=center, inner sep=0}, "(F\o\T)\eta"', from=2-1, to=2-2]
\arrow[""{name=4, anchor=center, inner sep=0}, "{F\beta_F}"', from=2-2, to=2-3]
\arrow[""{name=5, anchor=center, inner sep=0}, "{\epsilon_{\T'\o F}}"', from=2-3, to=2-4]
\arrow["{Fs\^\T_{G\o F}}"', from=1-2, to=2-2]
\arrow["{(F\o G)s\^{\T'}_F}", from=1-3, to=2-3]
\arrow["\Nat"{description}, draw=none, from=0, to=3]
\arrow["\Nat"{description}, draw=none, from=2, to=5]
\arrow["{(\beta;s\^\T,s\^{\T'})}"{description}, draw=none, from=1, to=4]
\end{tikzcd}
\end{equation*}
Let's now show the compatibility with the vertical lifts, i.e. $\alpha l\^{\T'}_F=Fl\^\T\alpha_\T\T'\alpha$:
\begin{equation*}
\adjustbox{scale=.7,center}{
\begin{tikzcd}
F\o\T & {F\o\T\o G\o F} &&& {F\o G\o\T'\o F} & F \\
& {F\o\T^2\o G\o F} & {F\o\T\o G\o\T'\o F} & {F\o\T\o G\o\T'\o F} & {F\o G\o\T'^2\o F} \\
&& {F\o\T\o G\o F\o G\o\T'\o F} & {F\o\T\o G\o F\o G\o\T'\o F} \\
& {F\o\T\o G\o F\o\T\o G\o F} &&& {F\o G\o\T'\o F\o G\o\T'\o F} \\
&& {F\o G\o\T'\o F\o\T\o G\o F} & {F\o G\o\T'\o F\o\T\o G\o F} \\
{F\o\T^2} & {F\o\T\o G\o F\o\T} & {F\o G\o\T'\o F\o\T} & {\T'\o F\o\T\o G\o F} & {\T'\o F\o\T\o G\o F} & {\T'^2\o F} \\
&& {\T'\o F\o\T} & {\T'\o F\o\T}
\arrow["{(F\o\T)\eta}", from=1-1, to=1-2]
\arrow["{F\beta_F}", from=1-2, to=1-5]
\arrow["{\epsilon_{\T'\o F}}", from=1-5, to=1-6]
\arrow[""{name=0, anchor=center, inner sep=0}, "{Fl\^\T}"', from=1-1, to=6-1]
\arrow[""{name=1, anchor=center, inner sep=0}, "{l\^{\T'}_F}", from=1-6, to=6-6]
\arrow["{(F\o\T)\eta_{G\o\T'\o F}}"', from=2-3, to=3-3]
\arrow[""{name=2, anchor=center, inner sep=0}, Rightarrow, no head, from=2-3, to=2-4]
\arrow["{(F\o\T)\eta_\T}"', from=6-1, to=6-2]
\arrow[""{name=3, anchor=center, inner sep=0}, "{F\beta_{F\o\T}}"', from=6-2, to=6-3]
\arrow["{(\T'\o F\o\T)\eta}"', from=7-4, to=6-4]
\arrow[""{name=4, anchor=center, inner sep=0}, "{(\T'\o F)\beta_F}"', from=6-4, to=6-5]
\arrow["{\T'\epsilon_{\T'\o F}}"', from=6-5, to=6-6]
\arrow["{\epsilon_{\T'\o F\o\T\o G\o F}}"{pos=0.8}, from=5-4, to=6-4]
\arrow[""{name=5, anchor=center, inner sep=0}, "{(F\o G\o\T'\o F)\beta_F}"{pos=0.3}, from=5-4, to=4-5]
\arrow["{(F\o\T\o G\o F\o\T)\eta}"', from=6-2, to=4-2]
\arrow["{\epsilon_{\T'\o F\o G\o\T'\o F}}"', from=4-5, to=6-5]
\arrow[""{name=6, anchor=center, inner sep=0}, "{(F\o\T\o G\o F)\beta_F}"'{pos=0.8}, from=4-2, to=3-3]
\arrow["{(F\o\T^2)\eta}", from=6-1, to=2-2]
\arrow["{(F\o\T)\eta_{\T\o G\o F}}", from=2-2, to=4-2]
\arrow["{(F\o G\o\T')\epsilon_{\T'\o F}}", from=4-5, to=2-5]
\arrow["{\epsilon_{\T'^2\o F}}", from=2-5, to=6-6]
\arrow["{(F\o\T)\beta_F}", from=2-2, to=2-3]
\arrow["{F\beta_{\T'\o F}}", from=2-4, to=2-5]
\arrow["{Fl\^\T_{G\o F}}"', from=1-2, to=2-2]
\arrow["{(F\o G)l\^{\T'}_F}", from=1-5, to=2-5]
\arrow["\Nat"{description}, draw=none, from=1-1, to=4-2]
\arrow["\Nat"{description}, draw=none, from=1-6, to=4-5]
\arrow["\Nat"{description}, draw=none, from=3-3, to=5-4]
\arrow["{(F\o\T\o G)\epsilon_{\T'\o F}}"', from=3-4, to=2-4]
\arrow[""{name=7, anchor=center, inner sep=0}, "{F\beta_{F\o G\o\T'\o F}}"'{pos=0.3}, from=3-4, to=4-5]
\arrow[""{name=8, anchor=center, inner sep=0}, Rightarrow, no head, from=3-3, to=3-4]
\arrow[""{name=9, anchor=center, inner sep=0}, "{F\beta_{F\o\T\o G\o F}}"{pos=0.8}, from=4-2, to=5-3]
\arrow["{(F\o G\o\T'\o F\o\T)\eta}"{pos=0.2}, from=6-3, to=5-3]
\arrow[""{name=10, anchor=center, inner sep=0}, Rightarrow, no head, from=5-3, to=5-4]
\arrow["{\epsilon_{\T'\o F\o\T}}"', from=6-3, to=7-3]
\arrow[""{name=11, anchor=center, inner sep=0}, Rightarrow, no head, from=7-3, to=7-4]
\arrow["{(\beta;l\^\T,l\^{\T'})}"{description}, draw=none, from=1-2, to=2-5]
\arrow["\Nat"{description, pos=0.7}, draw=none, from=0, to=6-2]
\arrow["\Nat"{description, pos=0.3}, draw=none, from=6-5, to=1]
\arrow["\Nat"{description, pos=0.3}, draw=none, from=2-2, to=6]
\arrow["\Nat"{description}, draw=none, from=5, to=4]
\arrow["\Nat"{description}, draw=none, from=9, to=3]
\arrow["\Nat"{description}, draw=none, from=10, to=11]
\arrow["\Delta"{description}, draw=none, from=2, to=8]
\arrow["\Nat"{description, pos=0.3}, draw=none, from=2-5, to=7]
\end{tikzcd}
}
\end{equation*}
Finally, the compatibility with the canonical flips, i.e. $\alpha_\T\T'\alpha c\^{\T'}_F=Fc\^\T\alpha_\T\T'\alpha$:
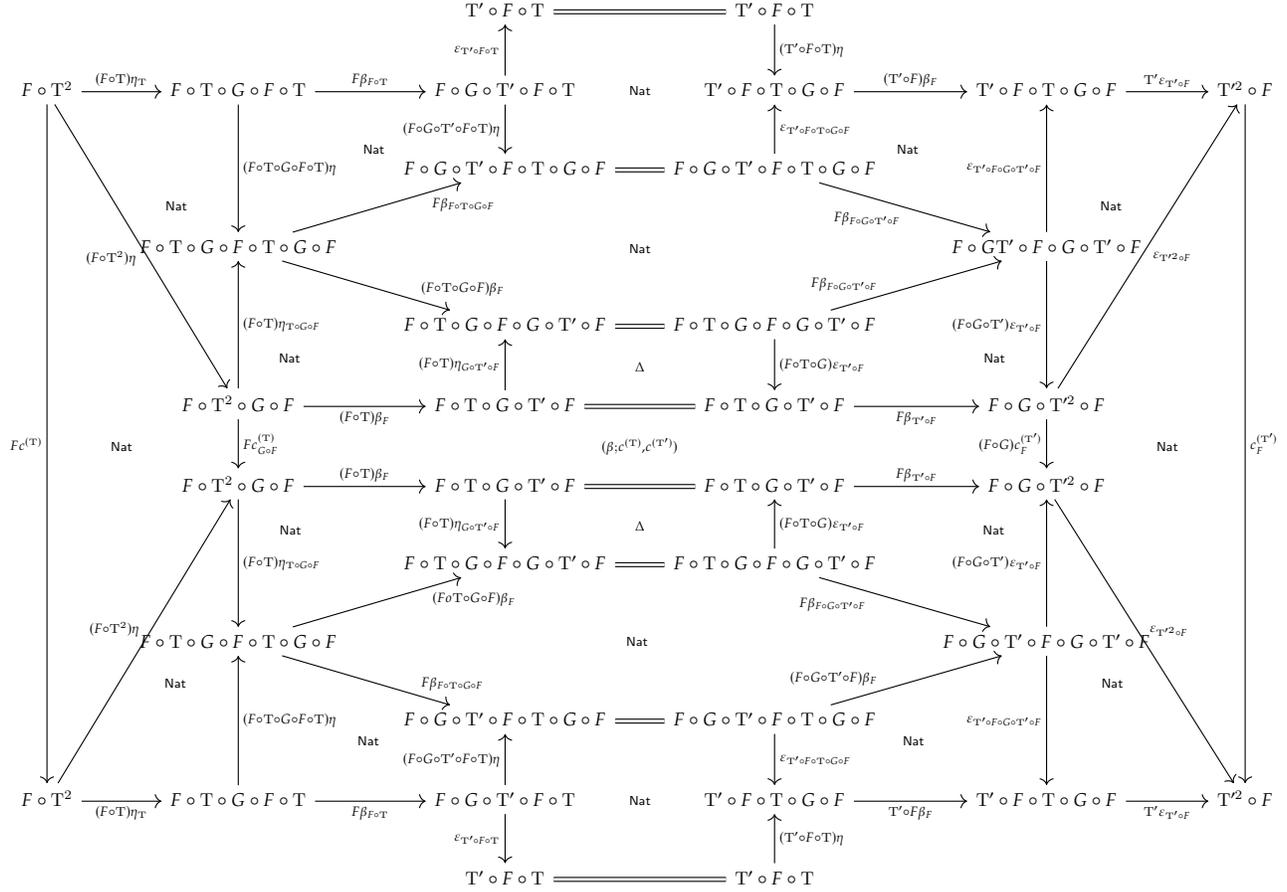
\begin{equation*}
\adjustbox{scale=.7,center}{
\begin{tikzcd}
&& {\T'\o F\o\T} & {\T'\o F\o\T} \\
{F\o\T^2} & {F\o\T\o G\o F\o \T} & {F\o G\o\T'\o F\o\T} & {\T'\o F\o\T\o G\o F} & {\T'\o F\o\T\o G\o F} & {\T'^2\o F} \\
&& {F\o G\o\T'\o F\o\T\o G\o F} & {F\o G\o\T'\o F\o\T\o G\o F} \\
& {F\o\T\o G\o F\o\T\o G\o F} &&& {F\o G\T'\o F\o G\o\T'\o F} \\
&& {F\o\T\o G\o F\o G\o\T'\o F} & {F\o\T\o G\o F\o G\o\T'\o F} \\
& {F\o\T^2\o G\o F} & {F\o\T\o G\o\T'\o F} & {F\o\T\o G\o\T'\o F} & {F\o G\o\T'^2\o F} \\
& {F\o\T^2\o G\o F} & {F\o\T\o G\o\T'\o F} & {F\o\T\o G\o\T'\o F} & {F\o G\o\T'^2\o F} \\
&& {F\o\T\o G\o F\o G\o\T'\o F} & {F\o\T\o G\o F\o G\o\T'\o F} \\
& {F\o\T\o G\o F\o\T\o G\o F} &&& {F\o G\o\T'\o F\o G\o\T'\o F} \\
&& {F\o G\o\T'\o F\o\T\o G\o F} & {F\o G\o\T'\o F\o\T\o G\o F} \\
{F\o\T^2} & {F\o\T\o G\o F\o\T} & {F\o G\o\T'\o F\o\T} & {\T'\o F\o\T\o G\o F} & {\T'\o F\o\T\o G\o F} & {\T'^2\o F} \\
&& {\T'\o F\o\T} & {\T'\o F\o\T}
\arrow[""{name=0, anchor=center, inner sep=0}, "{c\^{\T'}_F}", from=2-6, to=11-6]
\arrow["{(F\o\T)\eta_{G\o\T'\o F}}"', from=7-3, to=8-3]
\arrow[""{name=1, anchor=center, inner sep=0}, Rightarrow, no head, from=7-3, to=7-4]
\arrow[""{name=2, anchor=center, inner sep=0}, "{(F\o\T)\eta_\T}"', from=11-1, to=11-2]
\arrow[""{name=3, anchor=center, inner sep=0}, "{F\beta_{F\o\T}}"', from=11-2, to=11-3]
\arrow["{(\T'\o F\o\T)\eta}"', from=12-4, to=11-4]
\arrow[""{name=4, anchor=center, inner sep=0}, "{\T'\o F\beta_F}"', from=11-4, to=11-5]
\arrow[""{name=5, anchor=center, inner sep=0}, "{\T'\epsilon_{\T'\o F}}"', from=11-5, to=11-6]
\arrow["{\epsilon_{\T'\o F\o\T\o G\o F}}", from=10-4, to=11-4]
\arrow[""{name=6, anchor=center, inner sep=0}, "{(F\o G\o\T'\o F)\beta_F}"{pos=0.3}, from=10-4, to=9-5]
\arrow["{(F\o\T\o G\o F\o\T)\eta}"', from=11-2, to=9-2]
\arrow["{\epsilon_{\T'\o F\o G\o\T'\o F}}"', from=9-5, to=11-5]
\arrow[""{name=7, anchor=center, inner sep=0}, "{(Fo\T\o G\o F)\beta_F}"'{pos=0.8}, from=9-2, to=8-3]
\arrow["{(F\o\T^2)\eta}", from=11-1, to=7-2]
\arrow["{(F\o\T)\eta_{\T\o G\o F}}", from=7-2, to=9-2]
\arrow["{(F\o G\o\T')\epsilon_{\T'\o F}}", from=9-5, to=7-5]
\arrow["{\epsilon_{\T'^2\o F}}", from=7-5, to=11-6]
\arrow["{(F\o\T)\beta_F}", from=7-2, to=7-3]
\arrow["{F\beta_{\T'\o F}}", from=7-4, to=7-5]
\arrow["\Nat"{description}, draw=none, from=8-3, to=10-4]
\arrow["{(F\o\T\o G)\epsilon_{\T'\o F}}"', from=8-4, to=7-4]
\arrow[""{name=8, anchor=center, inner sep=0}, "{F\beta_{F\o G\o\T'\o F}}"'{pos=0.3}, from=8-4, to=9-5]
\arrow[""{name=9, anchor=center, inner sep=0}, Rightarrow, no head, from=8-3, to=8-4]
\arrow[""{name=10, anchor=center, inner sep=0}, "{F\beta_{F\o\T\o G\o F}}"{pos=0.8}, from=9-2, to=10-3]
\arrow["{(F\o G\o\T'\o F\o\T)\eta}", from=11-3, to=10-3]
\arrow[""{name=11, anchor=center, inner sep=0}, Rightarrow, no head, from=10-3, to=10-4]
\arrow["{\epsilon_{\T'\o F\o\T}}"', from=11-3, to=12-3]
\arrow[""{name=12, anchor=center, inner sep=0}, Rightarrow, no head, from=12-3, to=12-4]
\arrow[""{name=13, anchor=center, inner sep=0}, "{Fc\^\T}"', from=2-1, to=11-1]
\arrow["{(F\o\T)\beta_F}"', from=6-2, to=6-3]
\arrow[""{name=14, anchor=center, inner sep=0}, Rightarrow, no head, from=6-3, to=6-4]
\arrow["{F\beta_{\T'\o F}}"', from=6-4, to=6-5]
\arrow["{Fc\^\T_{G\o F}}", from=6-2, to=7-2]
\arrow["{(F\o G)c\^{\T'}_F}"', from=6-5, to=7-5]
\arrow[""{name=15, anchor=center, inner sep=0}, Rightarrow, no head, from=5-3, to=5-4]
\arrow["{(F\o\T)\eta_{G\o\T'\o F}}", from=6-3, to=5-3]
\arrow["{(F\o\T\o G)\epsilon_{\T'\o F}}", from=5-4, to=6-4]
\arrow["{(F\o\T)\eta_{\T\o G\o F}}"', from=6-2, to=4-2]
\arrow["{(F\o G\o\T')\epsilon_{\T'\o F}}"', from=4-5, to=6-5]
\arrow[""{name=16, anchor=center, inner sep=0}, "{F\beta_{F\o G\o\T'\o F}}"{pos=0.3}, from=5-4, to=4-5]
\arrow[""{name=17, anchor=center, inner sep=0}, "{(F\o\T\o G\o F)\beta_F}"{pos=0.8}, from=4-2, to=5-3]
\arrow[""{name=18, anchor=center, inner sep=0}, Rightarrow, no head, from=3-3, to=3-4]
\arrow[""{name=19, anchor=center, inner sep=0}, "{F\beta_{F\o\T\o G\o F}}"'{pos=0.8}, from=4-2, to=3-3]
\arrow[""{name=20, anchor=center, inner sep=0}, "{F\beta_{F\o G\o\T'\o F}}"', from=3-4, to=4-5]
\arrow[""{name=21, anchor=center, inner sep=0}, "{(F\o\T)\eta_\T}", from=2-1, to=2-2]
\arrow[""{name=22, anchor=center, inner sep=0}, "{F\beta_{F\o\T}}", from=2-2, to=2-3]
\arrow[""{name=23, anchor=center, inner sep=0}, Rightarrow, no head, from=1-3, to=1-4]
\arrow["{\epsilon_{\T'\o F\o\T}}", from=2-3, to=1-3]
\arrow["{(\T'\o F\o\T)\eta}", from=1-4, to=2-4]
\arrow[""{name=24, anchor=center, inner sep=0}, "{(\T'\o F)\beta_F}", from=2-4, to=2-5]
\arrow[""{name=25, anchor=center, inner sep=0}, "{\T'\epsilon_{\T'\o F}}", from=2-5, to=2-6]
\arrow["{\epsilon_{\T'\o F\o G\o\T'\o F}}", from=4-5, to=2-5]
\arrow["{(F\o\T\o G\o F\o\T)\eta}", from=2-2, to=4-2]
\arrow["{(F\o G\o\T'\o F\o\T)\eta}"', from=2-3, to=3-3]
\arrow["{\epsilon_{\T'\o F\o\T\o G\o F}}"', from=3-4, to=2-4]
\arrow["{\epsilon_{\T'^2\o F}}"', from=6-5, to=2-6]
\arrow["{(F\o\T^2)\eta}"', from=2-1, to=6-2]
\arrow["\Nat"{description, pos=0.3}, draw=none, from=7-2, to=7]
\arrow["\Delta"{description}, draw=none, from=1, to=9]
\arrow["\Nat"{description}, draw=none, from=6, to=4]
\arrow["\Nat"{description}, draw=none, from=10, to=3]
\arrow["\Nat"{description}, draw=none, from=11, to=12]
\arrow["\Nat"{description, pos=0.7}, draw=none, from=13, to=11-2]
\arrow["\Delta"{description}, draw=none, from=14, to=15]
\arrow["\Nat"{description, pos=0.3}, draw=none, from=6-2, to=17]
\arrow["\Nat"{description, pos=0.3}, draw=none, from=6-5, to=16]
\arrow["{(\beta;c\^\T,c\^{\T'})}"{description}, draw=none, from=14, to=1]
\arrow["\Nat"{description}, draw=none, from=20, to=24]
\arrow["\Nat"{description}, draw=none, from=19, to=22]
\arrow["\Nat"{description, pos=0.3}, draw=none, from=11-5, to=0]
\arrow["\Nat"{description, pos=0.3}, draw=none, from=2-5, to=0]
\arrow["\Nat"{description, pos=0.3}, draw=none, from=2-2, to=13]
\arrow["\Nat"{description}, draw=none, from=23, to=18]
\arrow["\Nat"{description}, Rightarrow, draw=none, from=18, to=15]
\arrow["\Nat"{description}, draw=none, from=21, to=2]
\arrow["\Nat"{description}, draw=none, from=25, to=5]
\arrow["\Nat"{description, pos=0.3}, draw=none, from=7-5, to=8]
\end{tikzcd}
}
\end{equation*}
So far, we proved that $(F,\alpha)$ is a lax tangent morphism. The next step is to prove that:
\begin{equation*}
\begin{tikzcd}
{(\X,\TT)} & {({\X'},\TT')} \\
{(\X,\TT)} & {(\X,\TT)}
\arrow["\shortmid"{marking}, Rightarrow, no head, from=1-1, to=2-1]
\arrow[Rightarrow, no head, from=2-1, to=2-2]
\arrow["{(G,\beta)}", "\shortmid"{marking}, from=1-2, to=2-2]
\arrow["{(F,\alpha)}", from=1-1, to=1-2]
\arrow["\eta"{description}, draw=none, from=2-1, to=1-2]
\end{tikzcd}\hfill
\begin{tikzcd}
{({\X'},\TT')} & {({\X'},\TT')} \\
{(\X,\TT)} & {({\X'},\TT')}
\arrow["{(G,\beta)}"', "\shortmid"{marking}, from=1-1, to=2-1]
\arrow["{(F,\alpha)}"', from=2-1, to=2-2]
\arrow[Rightarrow, no head, from=1-1, to=1-2]
\arrow["\shortmid"{marking}, Rightarrow, no head, from=1-2, to=2-2]
\arrow["\epsilon"{description}, draw=none, from=2-1, to=1-2]
\end{tikzcd}
\end{equation*}
are tangent double cells. This amounts to showing the commutativity of the following diagrams:
\begin{equation*}
\begin{tikzcd}
\T & {G\o F\o\T} \\
{\T\o G\o F} & {G\o F\o\T\o G\o F} \\
{G\o\T'\o F} & {G\o F\o G\o\T'\o F} \\
{G\o\T'\o F} & {G\o\T'\o F}
\arrow["{\eta_\T}", from=1-1, to=1-2]
\arrow["{(G\o F\o\T)\eta}", from=1-2, to=2-2]
\arrow["\T\eta"', from=1-1, to=2-1]
\arrow["{\eta_{\T\o G\o F}}"', from=2-1, to=2-2]
\arrow["\Nat"{description}, draw=none, from=1-2, to=2-1]
\arrow["{(G\o F)\beta_F}", from=2-2, to=3-2]
\arrow["{G\epsilon_{Q\o F}}", from=3-2, to=4-2]
\arrow[Rightarrow, no head, from=4-1, to=4-2]
\arrow["{\beta_F}"', from=2-1, to=3-1]
\arrow[Rightarrow, no head, from=3-1, to=4-1]
\arrow["{\eta_{G\o\T'\o F}}", from=3-1, to=3-2]
\arrow["\Nat"{description}, draw=none, from=2-2, to=3-1]
\arrow["\Delta"{description}, draw=none, from=3-1, to=4-2]
\end{tikzcd}\hfill
\begin{tikzcd}
{F\o\T\o G} & {F\o\T\o G} \\
{F\o\T\o G\o F\o G} & {F\o\T\o G} \\
{F\o G\o\T'\o F\o G} & {F\o G\o\T'} \\
{\T'\o F\o G} & \T'
\arrow[Rightarrow, no head, from=1-1, to=1-2]
\arrow[Rightarrow, no head, from=1-2, to=2-2]
\arrow["{(F\o\T)\eta_G}"', from=1-1, to=2-1]
\arrow["{(F\o\T\o G)\epsilon}"', from=2-1, to=2-2]
\arrow["F\beta", from=2-2, to=3-2]
\arrow["{F\beta_{F\o G}}"', from=2-1, to=3-1]
\arrow["{(F\o G\o\T')\epsilon}"', from=3-1, to=3-2]
\arrow["\Nat"{description}, draw=none, from=2-2, to=3-1]
\arrow["\Delta"{description}, draw=none, from=1-2, to=2-1]
\arrow["{\epsilon_{\T'\o F\o G}}"', from=3-1, to=4-1]
\arrow["\T'\epsilon"', from=4-1, to=4-2]
\arrow["{\epsilon_\T'}", from=3-2, to=4-2]
\arrow["\Nat"{description}, draw=none, from=3-2, to=4-1]
\end{tikzcd}
\end{equation*}
The converse is a straightforward computation we leave for the reader to spell out.
\end{proof}
\end{proposition}

Thanks to Proposition~\ref{proposition:conjuctions-tangent-morphisms} we can extends $\Alg(\P)$ to a covariant pseudofunctor which sends each morphism of operads $\varphi\colon\P\to\OprQ$ to a lax tangent morphism $(\varphi_!,\alpha_!)\colon\Alg(\P)\to\Alg(\OprQ)$.

\begin{proposition}
\label{proposition:functoriality-Alg-shriek}
The operation which takes an operad to its algebraic tangent category extends to a pseudofunctor $\Alg_!\colon\Operad\to\TngCat$ which sends each morphism of operads $\varphi\colon\P\to\OprQ$ to the lax tangent morphism $(\varphi_!,\beta_!)\colon\Alg(\P)\to\Alg(\OprQ)$, whose underlying functor is the left adjoint of $\varphi^\*$ and $\beta_!$ is defined as follows:
\begin{align*}
&\beta_!\colon\varphi_!\o\FootT\^\P\xrightarrow{\varphi_!\smallFootT\^\P\eta}\varphi_!\o\FootT\^\P\o\varphi^\*\o\varphi_!=\varphi_!\o\varphi^\*\o\FootT\^\OprQ\o\varphi_!\xrightarrow{\epsilon_{\smallFootT\^\OprQ\varphi_!}}\FootT\^\OprQ\o\varphi_!
\end{align*}
\end{proposition}

\begin{remark}
\label{remark:pseudofunctoriality-Alg-shriek}
Notice that $\Alg_!$ is only pseudofunctorial. This comes from the fact that the left adjoint of a functor is only unique up to a unique natural isomorphism. Such a natural isomorphism equips $\Alg_!$ with an associator and a left and a right unitor.
\end{remark}

In Proposition~\ref{proposition:functoriality-Alg-shriek}, we used that $\varphi^\*$ is a colax tangent morphism, since $\varphi^\*$ is a strict tangent morphism. To unwrap the definition of $\beta_!$ notice that, given a $\P$-algebra $A$, $\varphi_!(\FootT\^\P A)$ is the $\OprQ$-algebra generated by pairs $(a,b)$ for $a,b\in A$, satisfying some suitable relations defined by the coequalizer that defines $\varphi_!$. Similarly, also $\FootT\^\OprQ(\varphi_!A)$ is generated by pairs $(a,b)$ for $a,b\in A$. So, $\beta_!$ sends each generator $(a,b)$ to the corresponding generator $(a,b)$.

\begin{example}
\label{example:com-ass-Alg}
Consider the operads $\Ass$ and $\Com$, respectively known as the associative and the commutative operads. The corresponding algebras are the associative $R$-algebras and the commutative $R$-algebras, respectively. Concretely, $\Ass$ is generated by a $2$-ary operation $\mu$ which satisfies the following relation:
\begin{align*}
&\mu(1_\Ass,\mu)=\mu(\mu,1_\Ass)
\end{align*}
Similarly, $\Com$ is generated by a $2$-ary operation $\nu$ which satisfies the same associativity condition as $\mu$ and moreover is symmetric, i.e.:
\begin{align*}
&\nu\.\tau=\nu
\end{align*}
where $\tau\in\S_2$ is the permutation $(1\quad 2)$. Since $\nu$ satisfies the same relation as $\mu$, there is a quotient morphism $\varphi\colon\Ass\to\Com$ of operads, which sends $\mu$ to $\nu$, and that induces an adjunction:
\begin{align*}
\varphi_!\colon\Alg_\Ass\leftrightarrows\Alg_\Com\colon\varphi^\*
\end{align*}
$\varphi^\*$ sends a commutative algebra $B$ to the underlying associative algebra $\varphi^\*B$, while $\varphi_!$ sends an associative algebra $A$ to its abelianization $A/[A,A]$, where $[A,A]$ denotes the commutator, i.e. the ideal generated by symbols $ab-ba$, for any $a,b\in A$.
\par The functor $\Alg^\*$ maps the morphism of operads $\varphi$ to the strict tangent morphism over the pullback functor $\varphi^\*$, which makes $\Alg(\Com)$ a tangent subcategory of $\Alg(\Ass)$.
\par The functor $\Alg_!$ maps the morphism of operads $\varphi$ to the lax tangent morphism whose underlying functor is the abelianization functor $\varphi_!$. To understand what is the corresponding distributive law $\varphi_!\o\FootT\^\Ass\to\FootT\^\Com\o\varphi_!$, first notice that, for an associative algebra $A$, $\varphi_!(\FootT\^\Ass(A))$ is the abelianization of $A\ltimes A$. It is not hard to see that this is isomorphic to $\varphi_!(A)\ltimes\varphi_!(A)$ which is precisely $\FootT\^\Com(\varphi_!(A))$. On the other hand, the distributive law sends the generator $(a,b)\in \varphi_!(A\ltimes A)$ to $(a,b)\in\varphi_!(A)\ltimes\varphi_!(A)$. Thus, the distributive law is precisely the isomorphism between the abelianization of $A\ltimes A$ and the semi-direct product of the abelianization of $A$ with itself.
\end{example}

\begin{example}
\label{example:lie-ass-Alg}
Consider the operad $\Lie$, which generates Lie algebras. Concretely, $\Lie$ is the operad generated by a binary operation $\mu$ satisfying the following relations:
\begin{align*}
&\mu\.\tau=-\mu\\
&\mu(\mu,1_\Lie)+\mu(\mu,1_\Lie)\.\sigma+\mu(\mu,1_\Lie)\.\sigma^2=0
\end{align*}
Note that the first relation encodes the antisymmetry of Lie brackets, while the second one corresponds to the Jacobi identity, where $\tau\in\S_2$ is the permutation $(1\quad 2)$ and $\sigma\in\S_3$ is the cyclical permutation $(1\quad 2\quad 3)$. The interested reader can find detailed equivalent constructions of $\Lie$ in \cite[Section~13.2]{loday:operads}. There is a canonical morphism of operads $\varphi\colon\Lie\to\Ass$ (see~\cite[Sections~5.2.12 and~13.2.5]{loday:operads}). Consider the induced adjunction:
\begin{align*}
&\varphi_!\colon\Alg_\Lie\leftrightarrows\Alg_\Ass\colon\varphi^\*
\end{align*}
The pullback functor $\varphi^\*$ sends an associative algebra $A$ to the underlying Lie algebra with Lie brackets defined by the commutator $[a,b]\=ab-ba$. On the other hand, the left adjoint $\varphi_!$ sends a Lie algebra $\Lieg$ to its universal enveloping algebra $\EnvAlg_\Lie\Lieg$.
\par The functor $\Alg^\*$ sends $\varphi$ to the strict tangent morphism whose underlying functor is the pullback functor $\varphi^\*$. The functor $\Alg_!$ sends $\varphi$ to the lax tangent morphism whose underlying functor is the universal enveloping algebra functor $\varphi_!$. To understand the distributive law $\varphi_!\o\FootT\^\Lie\to\FootT\^\Ass\o\varphi_!$, we first take a closer look at $\varphi_!(\FootT\^\Lie(\Lieg))$ and $\FootT\^\Ass(\varphi_!(\Lieg))$, for a Lie algebra $\Lieg$. The former is the universal enveloping algebra of the semi-direct product $\Lieg\ltimes\Lieg$. Concretely, this is the associative algebra generated by pairs $(g,h)$ for each $g,h\in\Lieg$, satisfying the relation:
\begin{align}
\label{equation:relations-lie-ass-1}
&(g,h)(g',h')-(g',h')(g,h)=([g,g'],[g,h']+[h,g'])
\end{align}
The second one is the semi-direct product of the universal enveloping algebra with itself. Concretely, this is the associative algebra of pairs $(g,h)$ for $g,h\in\EnvAlg_\Lie\Lieg$, satisfying the relations:
\begin{align}
\label{equation:relations-lie-ass-2}
\begin{split}    
&(g,h)(g',h')=(gg',gh'+hg')\\
&gh-hg=[g,h]
\end{split}
\end{align}
It is straightforward to see that the latter relations imply the former ones, thus there is a canonical morphism of Lie algebras $\varphi_!(\FootT\^\Lie(\Lieg))\to\FootT\^\Ass(\varphi_!(\Lieg))$, which corresponds to the distributive law.
\end{example}

\begin{remark}
\label{remark:alpha-shriek-not-iso}
Given a morphism of operads $\varphi\colon\P\to\OprQ$, we could be tempted to think that $(\varphi_!,\beta_!)$ is strong, or maybe even strict, since $\varphi^\*$ is a strict tangent morphism. A counterexample is given by Example~\ref{example:lie-ass-Alg}: relations~\eqref{equation:relations-lie-ass-2} imply relations~\eqref{equation:relations-lie-ass-1}, but not vice versa. To see that, notice that, given a Lie algebra $\Lieg$, the associative multiplication of $\varphi_!(\FootT\^\Lie(\Lieg))$ is an operation of pairs, e.g. $(g,h)(g',h')$ and, a priori, there is no well-defined multiplication on single elements of $\Lieg$, while in $\FootT\^\Ass(\varphi_!(\Lieg))$ there is indeed a multiplication on the elements of $\Lieg$ itself that comes from the universal enveloping algebra $\varphi_!(\Lieg)$. To understand the reasons why the distributive law $\beta_!$ of $\varphi_!$ is not an isomorphism notice that even if $\beta_!$ is the mate of an isomorphism $\beta$ and that mating preserves pasting diagrams (see~\cite[Proposition~2.2]{kelly:mates-adjunctions}) this holds as long as the mates of the diagrams are well-defined. This is not the case for $\beta^{-1}$, which does not admit a mate along the adjunction $\varphi^\*\dashv\varphi_!$.
\end{remark}

So far we proved that the operation which takes an operad to its algebraic tangent categories extends to a pair of functors $\Alg^\*$ and $\Alg_!$. Now, we focus our attention on the geometric tangent category of an operad. We are going to employ the fact that the geometric tangent structure is the adjoint tangent structure of the algebraic one.
\par We briefly recall that a tangent structure $\FootTT$ over a category $\X$ is called \textbf{adjunctable} (in~\cite[Proposition~5.17]{cockett:tangent-cats} the authors introduced the ``dual tangent structure'' while in~\cite[Definition~2.2.1]{ikonicoff:operadic-algebras-tagent-cats} the authors use the expression ``having an adjoint tangent structure''. Here we use ``adjunctable tangent structure'') if for any positive integer $n$, the functor $\FootT_n\colon\X\to\X$, which sends each object $A\in\X$ to the $n$-fold pullback $\FootT_nA$ along the projection over $A$, admits a left adjoint $\T_n$. Cockett and Cruttwell proved in~\cite[Proposition~5.17]{cockett:tangent-cats} that if $\FootTT$ is adjunctable, then the opposite category $\X^\op$ of $\X$ admits a tangent structure $\TT$, called the \textbf{adjoint tangent structure} of $\FootTT$, whose tangent bundle functor is the left adjoint $\T$ of $\FootT$ and whose projection, zero morphism, sum morphism, vertical lift and canonical flip are mates of the corresponding natural transformations of $\FootTT$.
\par Thanks to~\cite[Corollary~2.2.4]{ikonicoff:operadic-algebras-tagent-cats}, if $\X$ has enough finite colimits, e.g. $\X$ is cocomplete, then a tangent structure $\FootTT$ over $\X$ is adjunctable if and only if the tangent bundle functor $\FootT$ admits a left adjoint $\T$. In the following we denote by $\adjTngCat$ the $2$-category of adjunctable tangent categories, lax tangent morphisms and tangent natural transformations.

\begin{proposition}
\label{proposition:functoriality-of-adjoint}
The operation which takes an adjunctable tangent category $(\X,\FootTT)$ to its associated adjoint tangent category $(\X^\op,\TT)$ extends to a pseudofunctor $(-)^\op\colon\adjTngCat\to\adjTngCat$, which equips the $2$-category $\adjTngCat$ with a \textbf{pseudoinvolution}, that is an endofunctor together with a natural isomorphism $(-)^\op\o(-)^\op\Rightarrow\id_\adjTngCat$. In particular, given two adjunctable tangent categories $(\X,\FootTT)$ and $({\X'},\FootTT')$ with adjoint tangent categories $(\X^\op,\TT)$ and $({\X'}^\op,\TT')$,respectively, and a lax tangent morphism $(F,\alpha)\colon(\X,\FootTT)\to({\X'},\FootTT')$, $(F,\alpha)^\op\colon(\X^\op,\TT)\to({\X'}^\op,\TT')$ is the lax tangent morphism whose underlying functor is $F^\op$ and whose lax distributive law $\alpha^\op$, is the mate of $\alpha$ along the adjunctions $(\theta,\tau)\colon\T\dashv\FootT$ and $({{\theta'}},{\tau'})\colon\T'\dashv\FootT'$, that is:
\begin{align*}
&\alpha^\op\colon\T'\o F\xrightarrow{\T' F\theta}\T'\o F\o\FootT\o\T\xrightarrow{\T'\alpha_{\T}}\T'\o\FootT'\o F\o\T\xrightarrow{{\tau'}_{F\T}}F\o\T
\end{align*}
regarded as a morphism in ${\X'}$.
\begin{proof}
By definition, the natural transformations (i.e. projection etcetera) of the adjoint tangent structure $\TT$ of a tangent structure $\FootTT$ are mates along the adjunction $(\theta,\tau)\colon\T\dashv\FootT$ between the tangent bundle functors of the corresponding natural transformations of $\FootTT$. Thanks to~\cite[Proposition~2.2]{kelly:mates-adjunctions}, the mate of a pasting diagram is the pasting diagram of the mates, as long as the mate of each morphism of the diagram is well-defined. Therefore, given a lax tangent morphism $(F,\alpha)$ the distributive law $\alpha^\op$ is compatible with the tangent structures and thus $(F^\op,\alpha^\op)$ is a lax tangent morphism between the corresponding adjoint tangent categories. To prove that $(-)^\op$ is a pseudofunctor notice first that, given three adjunctable tangent categories $(\X,\FootTT),({\X'},\FootTT')$ and $({\X''},\FootTT'')$ with adjoint tangent categories $(\X^\op,\TT),({\X'}^\op,\TT')$ and $({\X''}^\op,\TT'')$, respectively, and two lax tangent morphisms $(F,\alpha)\colon(\X,\FootTT)\to({\X'},\FootTT')$ and $(G,\beta)\colon({\X'},\FootTT')\to({\X''},\FootTT'')$, the composition of $(F^\op,\alpha^\op)$ with $(G^\op,\beta^\op)$ is $(G^\op\o F^\op,G\alpha^\op\o\beta^\op_F)$. This must be compared with the opposite of the composition $(G\o F,\beta_F\o G\alpha)$. However, for the pasting diagram property of mates, these are the same lax tangent morphism. Similarly, we can argue that $(\id_\X^\op,\id_{\smallFootT}^\op)$ corresponds precisely to $(\id_{\X^\op},\id_{\T})$. Finally, notice that if $(\X,\FootTT)$ is adjunctable, then so is its adjoint tangent category $(\X^\op,\TT)$ and its adjoint is (isomorphic to) $(\X,\FootTT)$. This makes $(-)^\op$ a pseudoinvolution over $\adjTngCat$.
\end{proof}
\end{proposition}

\begin{remark}
\label{remark:pseudofunctoriality-adjoints}
We point out that $(-)^\op$ defined by Proposition~\ref{proposition:functoriality-of-adjoint} is only a pseudofunctor and not a strict functor because the choice of a left adjoint for the tangent bundle functor $\FootT$ is only unique up to a unique isomorphism. This implies that associativity and unitality are only defined up to a unique isomorphism, which defines the associator and the left and the right unitors of $(-)^\op$. 
\end{remark}

\begin{remark}
\label{remark:colax-tangent-morphisms-dont-adjoint-well}
One could hope that a similar pseudoinvolution $(-)^\op$ could also occur in the $2$-category $\adjTngCat_\co$ of adjunctable tangent categories, colax tangent morphisms, and corresponding tangent natural transformations. However, this is not the case. The reason is that mates of the colax distributive laws along the adjunctions of the tangent bundle functors are simply not well-defined. This breaking of symmetry plays a crucial role in understanding the differences between non-commutative algebraic geometry and the geometry of commutative affine schemes. We will come back to this point later in Example~\ref{example:com-ass-Geom}.
\end{remark}

Before proving the functoriality of the operation which takes an operad to its geometric tangent category, we notice an interesting fact.

\begin{lemma}
\label{lemma:double-mates-preserve-strength}
Consider a strong tangent morphism $(G,\alpha)\colon({\X'},\FootTT')\to(\X,\FootTT)$ between two adjunctable tangent categories. Suppose also that the functor $G$ has a left adjoint $F\dashv G$ and denote by $\beta\=\alpha^{-1}\colon\FootT\o G\Rightarrow G\o\FootT'$ the inverse of $\alpha\colon G\o\FootT'\Rightarrow\FootT\o G$. Then the corresponding tangent morphism $(F^\op,(\beta_!)^\op)\colon(\X^\op,\TT)\to({\X'}^\op,\TT')$ over the left adjoint $F$ and between the adjoint tangent categories is also strong.
\begin{proof}
By Proposition~\ref{proposition:conjuctions-tangent-morphisms}, the mate of $\beta$ along the adjunction $F\dashv G$ defines a lax tangent morphism $(F,\beta_!)\colon(\X,\FootTT)\to({\X'},\FootTT')$, where $\beta_!\colon F\o\FootT\Rightarrow\FootT'\o F$.
\par By Proposition~\ref{proposition:functoriality-of-adjoint}, the mate of the distributive law $\alpha$ along the adjunctions between the tangent bundle functors and their left adjoint defines a lax tangent morphism $(G^\op,\alpha^\op)\colon({\X'}^\op,\TT')\to(\X^\op,\TT)$, so that, as an $\X$-morphism, $\alpha^\op\colon\T\o G\Rightarrow G\o\T'$. Similarly, $\beta_!$ defines, again by mating, a lax tangent morphism $(F^\op,(\beta_!)^\op)\colon(\X^\op,\TT)\to({\X'}^\op,\TT')$ , so that, as a ${\X'}$-morphism, $(\beta_!)^\op\colon\T'\o F\Rightarrow F\o\T$. Interestingly, $\alpha^\op$ admits a second mate along the adjunction $(\eta,\epsilon)\colon F\dashv G$:
\begin{align*}
&(\alpha^\op)_!\colon F\o\T\xrightarrow{F\T\eta}F\o\T\o G\o F\xrightarrow{F\alpha^\op_F}F\o G\o\T'\o F\xrightarrow{\epsilon_{\T' F}}\T'\o F
\end{align*}
regarded as a morphism in ${\X'}$. Thus, we also obtain a colax tangent morphism $(F^\op,(\alpha^\op)_!)\colon(\X^\op,\TT)\nto({\X'}^\op,\TT')$. To prove that $(\alpha^\op)_!$ is the inverse of $(\beta_!)^\op$, consider the following diagram:
\begin{equation*}
\adjustbox{scale=.6,center}{
\begin{tikzcd}
F\o\T & {F\o\T\o G\o F} & {F\o\T\o G\o\FootT'\o\T'\o F} & {F\o\T\o\FootT\o G\o\T'\o F} & {F\o G\o\T'\o F} & {\T'\o F} \\
F\o\T\o\FootT\o\T & {F\o\T\o G\o F\o\FootT\o\T} & {F\o\T\o G\o\FootT'\o\T'\o F\o\FootT\o\T} &&& {\T'\o F\o\FootT\o\T} \\
{F\o\T\o\FootT\o G\o F\o\T} & {F\o\T\o G\o F\o\FootT\o G\o F\o\T} & {F\o\T\o G\o\FootT'\o\T'\o F\o\FootT\o G\o F\o\T} &&& {\T'\o F\o\FootT\o G\o F\o\T} \\
{F\o\T\o G\o\FootT'\o F\o\T} & {F\o\T\o G\o F\o G\o\FootT'\o F\o\T} & {F\o\T\o G\o\FootT'\o\T'\o F\o G\o\FootT'\o F\o\T} &&& {\T'\o F\o G\o\FootT'\o F\o\T} \\
& {F\o\T\o G\o\FootT'\o F\o\T} & {F\o\T\o G\o\FootT'\o\T'\o\FootT'\o F\o\T} &&& {\T'\o\FootT'\o F\o\T} \\
&& {F\o\T\o G\o\FootT'\o F\o\T} & {F\o\T\o\FootT\o G\o F\o\T} & {F\o G\o F\o\T} & {F\o T}
\arrow["F\T\eta", from=1-1, to=1-2]
\arrow["{F\T G{{\theta'}}_F}", from=1-2, to=1-3]
\arrow["{F\T\alpha_{\T' F}}", from=1-3, to=1-4]
\arrow["{F\tau_{G\T' F}}", from=1-4, to=1-5]
\arrow["{\epsilon_{\T' F}}", from=1-5, to=1-6]
\arrow["{\T' F\theta}", from=1-6, to=2-6]
\arrow["{\T' F\FootT\eta_{\T}}", from=2-6, to=3-6]
\arrow["{\T' F\beta_{F\T}}", from=3-6, to=4-6]
\arrow["{\T'\epsilon_{\FootT' F\T}}", from=4-6, to=5-6]
\arrow["{{\tau'}_{F\T}}", from=5-6, to=6-6]
\arrow["{F\T G\FootT'\T' F\theta}", from=1-3, to=2-3]
\arrow["{F\T G\FootT'\T' F\FootT\eta_{\T}}", from=2-3, to=3-3]
\arrow["{F\T G\FootT'\T' F\beta_{F\T}}", from=3-3, to=4-3]
\arrow["{F\T G\FootT'\T'\epsilon_{\FootT' F\T}}", from=4-3, to=5-3]
\arrow["{F\T G\FootT'{\tau'}_{F\T}}", from=5-3, to=6-3]
\arrow["{F\T\alpha_{F\T}}"', from=6-3, to=6-4]
\arrow["{F\tau_{GF\T}}"', from=6-4, to=6-5]
\arrow["{\epsilon_{F\T}}"', from=6-5, to=6-6]
\arrow["\Nat"{description}, draw=none, from=1-6, to=6-3]
\arrow[""{name=0, anchor=center, inner sep=0}, Rightarrow, no head, from=5-2, to=6-3]
\arrow["{F\T G{{\theta'}}_{\FootT' F\T}}", from=5-2, to=5-3]
\arrow["{F\T G F\theta}", from=1-2, to=2-2]
\arrow["{F\T GF\FootT\eta_{\T}}", from=2-2, to=3-2]
\arrow["{F\T G F\beta_{F\T}}", from=3-2, to=4-2]
\arrow["{F\T G\epsilon_{\FootT' F\T}}", from=4-2, to=5-2]
\arrow["\Nat"{description}, draw=none, from=1-2, to=6-3]
\arrow["F\T\theta"', from=1-1, to=2-1]
\arrow["{F\T\FootT\eta_{\T}}"', from=2-1, to=3-1]
\arrow["{F\T\beta_{F\T}}"', from=3-1, to=4-1]
\arrow[""{name=1, anchor=center, inner sep=0}, Rightarrow, no head, from=4-1, to=5-2]
\arrow["{F\T\eta_{G\FootT' F\T}}", from=4-1, to=4-2]
\arrow["\Nat"{description}, draw=none, from=1-1, to=4-2]
\arrow["\Delta"{description}, draw=none, from=0, to=5-3]
\arrow["\Delta"{description}, draw=none, from=1, to=4-2]
\end{tikzcd}
}
\end{equation*}
where $(\eta,\epsilon)\colon F\dashv G$, $(\theta,\tau)\colon\T\dashv\FootT$ and $({{\theta'}},{\tau'})\colon\T'\dashv\FootT'$. This shows that the following diagram commutes:
\begin{equation*}
\begin{tikzcd}
F\o\T &&& \T'\o F \\
F\o\T\o\FootT\o\T \\
{F\o\T\o\FootT\o G\o F\o\T} \\
{F\o\T\o G\o\FootT'\o F\o\T} & {F\o\T\o\FootT\o G\o F\o\T} & {F\o G\o F\o\T} & F\o\T
\arrow["{(\alpha^\op)_!}", from=1-1, to=1-4]
\arrow["{(\beta_!)^\op}", from=1-4, to=4-4]
\arrow["F\T\theta"', from=1-1, to=2-1]
\arrow["{F\T\FootT\eta_{\T}}"', from=2-1, to=3-1]
\arrow["{F\T\beta_{F\T}}"', from=3-1, to=4-1]
\arrow["{F\T\alpha_{F\T}}"', from=4-1, to=4-2]
\arrow["{F\tau_{GF\T}}"', from=4-2, to=4-3]
\arrow["{\epsilon_{F\T}}"', from=4-3, to=4-4]
\end{tikzcd}
\end{equation*}
However:
\begin{equation*}
\begin{tikzcd}
F\o\T &&& F\o\T \\
F\o\T\o\FootT\o\T && F\o\T \\
{F\o\T\o\FootT\o G\o F\o\T} & {F\o\T\o\FootT\o G\o F\o\T} \\
{F\o\T\o G\o\FootT'\o F\o\T} & {F\o\T\o\FootT\o G\o F\o\T} & {F\o G\o F\o\T} & F\o\T
\arrow["F\T\theta"', from=1-1, to=2-1]
\arrow["{F\T\FootT\eta_{\T}}"', from=2-1, to=3-1]
\arrow["{F\T\beta_{F\T}}"', from=3-1, to=4-1]
\arrow["{F\T\alpha_{F\T}}"', from=4-1, to=4-2]
\arrow["{F\tau_{GF\T}}"', from=4-2, to=4-3]
\arrow["{\epsilon_{F\T}}"', from=4-3, to=4-4]
\arrow[""{name=0, anchor=center, inner sep=0}, "{F\tau_{\T}}"{description}, from=2-1, to=2-3]
\arrow[""{name=1, anchor=center, inner sep=0}, "{F\eta_{\T}}"{description}, from=2-3, to=4-3]
\arrow[""{name=2, anchor=center, inner sep=0}, Rightarrow, no head, from=1-1, to=1-4]
\arrow[""{name=3, anchor=center, inner sep=0}, Rightarrow, no head, from=1-4, to=4-4]
\arrow[Rightarrow, no head, from=1-4, to=2-3]
\arrow[Rightarrow, no head, from=3-1, to=3-2]
\arrow["\Nat"{description}, draw=none, from=2-3, to=3-2]
\arrow[Rightarrow, no head, from=3-2, to=4-2]
\arrow["{\beta=\alpha^{-1}}"{description}, draw=none, from=3-2, to=4-1]
\arrow["\Delta"{description}, draw=none, from=2, to=0]
\arrow["\Delta"{description}, draw=none, from=3, to=1]
\end{tikzcd}
\end{equation*}
We just proved that $(\beta_!)^\op\o(\alpha^\op)_!=\id_{F\T}$. Similarly, one can prove also that converse and conclude that $(\alpha^\op)_!$ is the inverse of $(\beta_!)^\op$, as expected.
\end{proof}
\end{lemma}

\begin{remark}
\label{remark:strength-opposite-right-adjoints}
Given a pair of conjoints $(F,\beta_!)\dashv(G,\alpha)$ in the double category of tangent categories where $(G,\alpha)$ is a strong tangent morphism, Lemma~\ref{lemma:double-mates-preserve-strength} establishes that the pseudofunctor $(-)^\op$ maps $(F,\beta_!)\dashv(G,\alpha)$ to another pair of conjoints $(G^\op,\alpha^\op)\dashv(F^\op,(\beta_!)^\op)$ and that $(F^\op,(\beta_!)^\op)$ is also a strong tangent morphism. However, if $(G,\alpha)$ is strict this does not imply that $(F^\op,(\beta_!)^\op)$ is strict as well.
\end{remark}

In the following diagram, we represent the proof of Lemma~\ref{lemma:double-mates-preserve-strength}.
\begin{equation*}
\begin{tikzcd}
& \beta & {\beta_!} \\
\alpha && {(\beta_!)^\op} \\
{\alpha^\op} & {(\alpha^\op)_!}
\arrow["{F\dashv G}", from=1-2, to=1-3]
\arrow["{(-)^\op}", from=1-3, to=2-3]
\arrow["{(-)^\op}"', from=2-1, to=3-1]
\arrow["{F\dashv G}"', from=3-1, to=3-2]
\arrow["{\text{inverses}}", leftrightarrow, from=2-1, to=1-2]
\arrow["{\text{inverses}}"', dashed, leftrightarrow, from=3-2, to=2-3]
\end{tikzcd}
\end{equation*}
Starting from $\beta$, which is the inverse of the strong distributive law $\alpha$, by moving to the right, i.e. by mating along the adjunction $F\dashv G$, we obtain a lax distributive $\beta_!$, which, as noticed in Remark~\ref{remark:alpha-shriek-not-iso}, in general, is not invertible.  By moving down from $\beta_!$, by applying the pseudofunctor $(-)^\op$, we obtain a lax distributive law $(\beta_!)^\op$. Similarly, by starting from $\alpha$ and moving down, i.e. applying $(-)^\op$, we obtain a lax distributive law $\alpha^\op$, which, as mentioned in Remark~\ref{remark:pseudofunctoriality-adjoints}, in general, is not invertible. Finally, by moving from $\alpha^\op$ to the right, i.e. by mating along the adjunction $F\dashv G$, we obtain a colax distributive law $(\alpha^\op)_!$ which turns out to be the inverse of $(\beta_!)^\op$.
\par We can now prove the functoriality of the operation which takes an operad to its associated geometric tangent category. Similarly, as for the algebraic counterpart of this construction, this operation extends to two functors, one mapping operad morphisms $\varphi$ to a lax tangent morphism whose underlying functor is $(\varphi^\*)^\op$ and the second to a strong tangent morphism whose underlying functor is $\varphi_!^\op$.

\begin{proposition}
\label{proposition:functoriality-Geom}
The operation which takes an operad $\P$ to its associated geometric tangent category $\Geom(\P)$ extends to a contraviarant pseudofunctor $\Geom^\*\colon\Operad^\op\to\TngCat$ which sends a morphism of operads $\varphi\colon\P\to\OprQ$ to the lax tangent morphism $(\varphi^\*,\alpha^\*)\colon\Geom(\OprQ)\to\Geom(\P)$, where $\alpha^\*$ is defined as follows:
\begin{align*}
&\alpha^\*\colon\T\^\P\o\varphi^\*\xrightarrow{\T\^\P\varphi^\*\theta\^\OprQ}\T\^\P\o\varphi^\*\o\FootT\^\OprQ\o\T\^\OprQ=\T\^\P\o\FootT\^\P\o\varphi^\*\o\T\^\OprQ\xrightarrow{\tau\^\P_{\varphi^\*\T\^\OprQ}}\varphi^\*\o\T\^\OprQ
\end{align*}
where $(\theta\^\P,\tau\^\P)\colon\T\^\P\dashv\FootT\^\P$ and $(\theta\^\OprQ,\tau\^\OprQ)\colon\T\^\OprQ\dashv\FootT\^\OprQ$.
Moreover, the same operation extends also to a covariant pseudofunctor $\Geom_!\colon\Operad\to\TngCat_{\cong}$ which sends a morphism of operads $\varphi\colon\P\to\OprQ$ to the strong tangent morphism $(\varphi_!,\alpha_!)\colon\Geom(\P)\to\Geom(\OprQ)$, where $\alpha_!$ is defined as follows:
\begin{align*}
\alpha_!\colon\T\^\OprQ\o\varphi_!\xrightarrow{\T\^\OprQ\varphi_!\theta\^\P}\T\^\OprQ\o\varphi_!\o\FootT\^\P\o\T\^\P\xrightarrow{\T\^\OprQ(\beta_!)_{\T\^\P}}\T\^\OprQ\o\FootT\^\OprQ\o\varphi_!\o\T\^\P\xrightarrow{\tau\^\OprQ_{\varphi_!\T\^\P}}\varphi_!\o\T\^\P
\end{align*}
where $\beta_!$ is defined as in Proposition~\ref{proposition:functoriality-Alg-shriek}.
\end{proposition}

Concretely, given a morphism $\varphi\colon\P\to\OprQ$ and a $\OprQ$-algebra $B$, $\varphi^\*(\T\^\OprQ B)$ is a $\P$-algebra generated by all $b\in B$ and by symbols $\d\^\OprQ b$, for $b\in B$, satisfying suitable relations. On the other hand, $\T\^\P(\varphi^\* B)$ is generated by all $b\in B$ and by symbols $\d\^\P b$, for $b\in B$, satisfying suitable relations. Thus, the distributive law $\alpha^\*\colon\T\^\P(\varphi^\* B)\to\varphi^\*(\T\^\OprQ B)$ associated with $\varphi^\*$ sends each $b$ to $b$ and each $\d\^\P b$ to $\d\^\OprQ b$.
\par Similarly, given a $\P$-algebra $A$, $\varphi_!(\T\^\P A)$ is generated by all $a\in A$ and by $\d\^\P a$ for $a\in A$, satisfying suitable relations. On the other hand, $\T\^\OprQ(\varphi_!A)$ is generated by all $a\in A$ and by $\d\^\OprQ a$, for $a\in A$, satisfying suitable relations. Thus, the distributive law $\alpha_!\colon\varphi_!(\T\^\P A)\to\T\^\OprQ(\varphi_!A)$ sends each $a$ to $a$ and each $\d\^\P a$ to $\d\^\OprQ a$.

\begin{example}
\label{example:com-ass-Geom}
In Example~\ref{example:lie-ass-Alg} we showed how the canonical morphism of operads $\varphi\colon\Ass\to\Com$ is mapped by the functors $\Alg^\*$ and $\Alg_!$. The functor $\Geom^\*$ maps $\varphi$ to the lax tangent morphism defined over the pullback functor $\varphi^\*$. Interestingly, this lax tangent morphism is not strong, i.e. the distributive law $\T\^\Ass\o\varphi^\*\to\varphi^\*\o\T\^\Com$ (as a $\Ass$-algebra morphism) is not an isomorphism.
\par To prove that, notice that the module of K\"ahler differentials $\Kahler\^\Com A$ of a commutative algebra $A$ is given by quotienting the ideal $I\=\ker(\nu\colon A\otimes_RA\to A)$, where $\nu$ represents the multiplication of $A$, by $I^2$, i.e. $\Kahler A=I/I^2$. If $B$ is an associative algebra, the corresponding module of K\"ahler differentials $\Kahler\^\Ass B$ is simply given by the ideal $I\=\ker(\mu\colon B\otimes_RB\to B)$, where $\mu$ represents the multiplication of $B$ (see~\cite{ginzburg:notes-on-nc-geometry} for a detailed description of both the modules of K\"ahler differentials in the commutative and in the associative case).
\par Thus, for a commutative algebra $A$, there is a natural quotient map $\Kahler\^\Ass\pi^\*(A)=I\to I/I^2=\Kahler\^\Com A$. The distributive law is induced precisely by this quotient map since it maps the symbols $\d\^\Ass a$ to $\d\^\Com a$. If the distributive law was an isomorphism such a comparison map between the modules of K\"ahler differentials would be invertible, which is clearly not. We note that a similar argument was used by Ginzburg in~\cite{ginzburg:notes-on-nc-geometry} to distinguish between ``noncommutative geometry \textit{in the small}, and noncommutative geometry \textit{in the large}'', meaning that the former ``is a generalization of the conventional ‘commutative’ algebraic geometry to
the noncommutative world''. The latter instead  ``is not a generalization of commutative theory. The world
of noncommutative geometry ‘in the large’ does not contain commutative world
as a special case, but is only similar, parallel, to it.'' (\cite[Introduction]{ginzburg:notes-on-nc-geometry}).
\par Finally, the functor $\TT_!$ maps the morphism of operads $\varphi$ to the strong tangent morphism whose underlying functor is the (opposite of the) abelianization functor $\varphi_!$. The corresponding distributive law $\T\^\Com\o\varphi_!\Rightarrow\varphi_!\o\T\^\Ass$ (as a commutative algebra morphism) is the commutative algebra morphism:
\begin{align*}
&\T\^\Com\left({A}/{[A,A]}\right)\to{\T\^\Ass A}/{\left[\T\^\Ass A,\T\^\Ass A\right]}
\end{align*}
which sends the generators $[a]$ and $\d\^\Com[a]$ to $[a]$ and $[\d\^\Ass a]$, respectively, where we used the square brackets to indicate the left coset given by the commutator and an element of the associative algebra $A$. It is not hard to see that the algebra morphism $\T\^\Ass A\to\T\^\Com(A/[A,A])$ which sends each $a$ to $[a]$ and $\d\^\Ass a$ to $\d\^\Com[a]$ is well-defined and provides an inverse for the distributive law.
\end{example}

\begin{example}
\label{example:lie-ass-Geom}
In Example~\ref{example:lie-ass-Alg} we showed how the canonical morphism of operads $\varphi\colon\Lie\to\Ass$ is mapped by the functors $\Alg^\*$ and $\Alg_!$. The functor $\Geom^\*$ maps $\varphi$ to the lax tangent morphism whose underlying functor is (the opposite of) $\varphi^\*$. In order to understand the distributive law $\T\^\Lie\o\varphi^\*\Rightarrow\varphi^\*\o\T\^\Ass$ (as an associative algebra morphism), let's first take a closer look at $\T\^\Lie(\varphi^\*(A))$ and $\varphi^\*(\T\^\Ass(A))$ for an associative algebra $A$. The former one is the Lie algebra generated by $a\in A$ and by symbols $\d\^\Lie a$ for each $a\in A$, satisfying the following relations:
\begin{align*}
&[a,b]=ab-ba\\
&\d\^\Lie([a,b])=[\d\^\Lie a,b]+[a,\d\^\Lie b]
\end{align*}
The second algebra is generated by $a\in A$ and by symbols $\d\^\Ass a$ for each $a\in A$, satisfying the following relations:
\begin{align*}
&[a,b]=ab-ba\\
&\d\^\Ass(ab)=\d\^\Ass a\.b+a\d\^\Ass b\\
&[a,\d\^\Ass b]=a\d\^\Ass b-\d\^\Ass b\.a\\
&[\d\^\Ass a,\d\^\Ass b]=\d\^\Ass a\d\^\Ass b-\d\^\Ass b\d\^\Ass a
\end{align*}
Note that the relations of the former one are implied by the relations of the latter. The canonical quotient map $\T\^\Lie(\varphi^\*(A))\to\varphi^\*(\T\^\Ass(A))$ corresponds the distributive law. Note that such a map is not an isomorphism.
\par Finally, the functor $\Geom_!$ maps $\varphi$ to the lax tangent morphism whose underlying functor is the (opposite of the) universal enveloping algebra functor $\varphi_!$. To understand the distributive law $\T\^\Ass\o\varphi_!\Rightarrow\varphi_!\o\T\^\Lie$ (as an associative algebra morphism), we first take a closer look at $\T\^\Ass(\varphi_!(\Lieg))$ and $\varphi_!(\T\^\Lie(\Lieg))$ for a Lie algebra $\Lieg$. The former is the associative algebra generated by all $g\in\Lieg$ and by symbols $\d\^\Ass g$ for each $g\in\Lieg$ and satisfying the relations:
\begin{align*}
&gh-hg=[g,h]\\
&\d\^\Ass(gh)=\d\^\Ass g\.h+g\d\^\Ass h
\end{align*}
The latter is the associative algebra generated by $g\in\Lieg$ and by symbols $\d\^\Lie g$ for each $g\in\Lieg$, satisfying the relations:
\begin{align*}
&gh-hg=[g,h]\\
&\d\^\Lie g\.h-h\d\^\Lie g=[\d\^\Lie g,h]\\
&g\d\^\Lie h-\d\^\Lie h\.g=[g,\d\^\Lie h]\\
&\d\^\Lie g\.\d\^\Lie h-\d\^\Lie h\.\d\^\Lie g=[\d\^\Lie g,\d\^\Lie h]\\
&\d\^\Lie[g,h]=[\d\^\Lie g,h]+[h,\d\^\Lie g]
\end{align*}
Because the first set of relations implies the latter, this allows us to define a morphism of associative algebras $\varphi_!(\T\^\Lie(\Lieg))\to\T\^\Ass(\varphi_!(\Lieg))$, which corresponds to the (inverse of the) distributive law. Thanks to Lemma~\ref{lemma:double-mates-preserve-strength}, this morphism is an isomorphism.
\end{example}

\section{The slice tangent category as a right adjoint functor}
\label{section:slice-tangent-category}
Rosick\'y proved that, under mild assumptions, the slice of a tangent category $(\X,\TT)$ over an object $A\in\X$ is still a tangent category (cf.~\cite{rosicky:tangent-cats}). Cockett and Cruttwell further investigated this construction and related this to the notion of tangent fibrations (cf.~\cite{cockett:differential-bundles}).
\par In this section, we prove an important result that shows the deep relationship between operads and tangent categories. In a nutshell, we show that the slice tangent category of the geometric tangent category $\Geom(\P)$ of an operad $\P$ over a $\P$-affine scheme $A\in\Geom(\P)$ is still the geometric tangent category $\Geom(\P\^A)$ of an operad $\P\^A$. In particular, $\P\^A$ is the enveloping operad of the $\P$-algebra $A$.
\par To prove this result we are going to show that the functor which associates to each pair $((\X,\TT),A)$ formed by a tangent category (sliceable over $A$) and an object $A\in\X$ the corresponding slice tangent category $(\X,\TT)/A$ fulfills the same universality condition of the functor that associates to each pair $(\P,A)$ formed by an operad $\P$ and a $\P$-algebra $A$ the corresponding enveloping operad $\P\^A$. This is not just an important connection between the world of operads and the one of tangent categories, but it also provides a new characterization for the construction of the slice tangent category in terms of a right adjoint functor and therefore it also constitutes a new result in tangent category theory.
\par The section is organized as follows: first, we recall the original definition of the slice tangent category of a tangent category over an object. Then, we introduce the new characterization of this construction in terms of a right adjoint functor. Let's start with the main definitions.

\begin{definition}
\label{definition:sliceability}
A tangent category $(\X,\TT)$ is \textbf{sliceable over} an object $A\in\X$ if for any $E\in\X$ and $f\colon E\to A$ in $\X$, the $\T$-pullback of $\T f$ along the zero morphism is well-defined, that is the following diagram:
\begin{equation}
\label{equation:sliceable-pullback}
\begin{tikzcd}
{\T\^AE} & {\T E} \\
A & {\T A}
\arrow["{\T f}", from=1-2, to=2-2]
\arrow["z"', from=2-1, to=2-2]
\arrow["{f^\*}"', dashed, from=1-1, to=2-1]
\arrow["{v_f}", dashed, from=1-1, to=1-2]
\arrow["\lrcorner"{anchor=center, pos=0.125}, draw=none, from=1-1, to=2-2]
\end{tikzcd}
\end{equation}
is a well-defined pullback diagram and for every positive integer $m$ the functor $\T^m\=\T\o\T\o\dots\o\T$ preserves its universality. We also say that a tangent category is \textbf{sliceable} if it is sliceable over all of its objects.
\end{definition}

Given a sliceable tangent category $(\X,\TT)$ over an object $A$ we can define a tangent bundle functor:
\begin{align*}
&\T\^A\colon\X/A\to\X/A
\end{align*}
which maps each morphism $f\colon E\to A$ in $\X$ to the unique morphism $f^\*\colon\T\^AE\to A$. We adopt the following notation: we will write $\T\^Af$ for the tangent bundle over $f\in\X/A$, regarded as an object in the slice category over $A$. Abusing notation, we also denote by $\T\^AE$ the domain of $\T\^Af\colon\T\^AE\to A$, regarded as a morphism of $\X$. Notice that $\T\^A$ is functorial in the slice category but not in $\X$.
\par This characterization of the slice tangent bundle functor, also known as the vertical tangent bundle functor, is due to Cockett and Cruttwell in their article on differential bundles and tangent fibrations. The equivalent original characterization of $\T\^A$ is due to Rosick\'y. For our purposes the Rosick\'y version is more useful, therefore we recall briefly here this construction. First, notice that a tangent category $(\X,\TT)$ is sliceable over $A\in\X$ if and only if for any morphism $f\colon E\to A$, the equalizer:
\begin{equation*}
\begin{tikzcd}
{\T\^AE} & {\T E} & {\T A}
\arrow["{v_f}", dashed, from=1-1, to=1-2]
\arrow["{\T f}", shift left=2, from=1-2, to=1-3]
\arrow["{\T fpz}"', shift right=2, from=1-2, to=1-3]
\end{tikzcd}
\end{equation*}
is well-defined and is a $\T$-equalizer, which means that for every positive integer $m$ the functor $\T^m$ preserves its universality. In the following, we denote by $v_f$ the equalizer map $v_f\colon\T\^AE\to\T E$. We can then give the following characterization:
\begin{description}
\item[tangent bundle functor] The tangent bundle functor $\T\^A\colon\X/A\to\X/A$ is defined as follows:
\begin{align*}
&\T\^A(f\colon E\to A)\colon\T\^AE\xrightarrow{v_f}\T E\xrightarrow{p}E\xrightarrow{f}A
\end{align*}
for any $f\in\X/A$. Moreover, given a morphism $g\colon(f\colon E\to A)\to(f'\colon E'\to A)$, i.e. $g\colon E\to E'$ such that $gf'=f$, we can define:
\begin{equation*}
\begin{tikzcd}
{\T\^AE} & {\T E} & {\T A} \\
{\T\^AE'} & {\T E'} & {\T A}
\arrow["{v_f}", from=1-1, to=1-2]
\arrow["{\T f}", shift left=2, from=1-2, to=1-3]
\arrow["{\T fpz}"', shift right=2, from=1-2, to=1-3]
\arrow["{v_{f'}}"', from=2-1, to=2-2]
\arrow["{\T f'}", shift left=2, from=2-2, to=2-3]
\arrow["{\T f'pz}"', shift right=2, from=2-2, to=2-3]
\arrow[Rightarrow, no head, from=1-3, to=2-3]
\arrow["{\T g}"{description}, from=1-2, to=2-2]
\arrow["{\T\^A g}"', dashed, from=1-1, to=2-1]
\end{tikzcd}
\end{equation*}
In particular, $\T\^A$ is functorial.
\item[projection] The projection $p\^A\colon\T\^Af\to f$, is defined as:
\begin{equation*}
\begin{tikzcd}
{\T\^AE} & \T E & E
\arrow["{v_f}", from=1-1, to=1-2]
\arrow["p", from=1-2, to=1-3]
\end{tikzcd}
\end{equation*}
\item[zero morphism] The zero morphism $z\^A\colon f\to\T\^Af$ is defined as the unique morphism that makes commuting the following diagram:
\begin{equation*}
\begin{tikzcd}
{\T\^AE} & \T E & \T A \\
B
\arrow["{v_f}", from=1-1, to=1-2]
\arrow["\T f", shift left=2, from=1-2, to=1-3]
\arrow["\T fpz"', shift right=2, from=1-2, to=1-3]
\arrow["z"', from=2-1, to=1-2]
\arrow["{z\^A}", dashed, from=2-1, to=1-1]
\end{tikzcd}
\end{equation*}
where we employed the universality of $v_f$;
\item[sum morphism] The sum morphism $s\^A\colon\T\^A_2f\to\T\^Af$ is defined as:
\begin{equation*}
\begin{tikzcd}
& {\T\^AE} & \T E & \T A \\
{\T\^A_2E} & {\T_2E}
\arrow["{v_f}", from=1-2, to=1-3]
\arrow["\T f", shift left=2, from=1-3, to=1-4]
\arrow["\T fpz"', shift right=2, from=1-3, to=1-4]
\arrow["s"', from=2-2, to=1-3]
\arrow[dashed, from=2-2, to=1-2]
\arrow["{v_f\times v_f}"', from=2-1, to=2-2]
\arrow["{s\^A}", dashed, from=2-1, to=1-2]
\end{tikzcd}
\end{equation*}
where we employed the universality of $v_f$;
\item[vertical lift] The lift $l\^A\colon\T\^Af\to(\T\^A)^2f$ is defined as:
\begin{equation*}
\begin{tikzcd}
{(\T\^A)^2E} & {\T\T\^AE} & {\T^2E} & \T A \\
{\T\^AE} & \T E
\arrow["{\T v_f}", from=1-2, to=1-3]
\arrow["{\T^2f}", shift left=2, from=1-3, to=1-4]
\arrow["{\T^2f\T p\T z}"', shift right=2, from=1-3, to=1-4]
\arrow["l"', from=2-2, to=1-3]
\arrow[dashed, from=2-2, to=1-2]
\arrow["{v_{\T\^Af}}", from=1-1, to=1-2]
\arrow["{v_f}"', from=2-1, to=2-2]
\arrow["{l\^A}", dashed, from=2-1, to=1-1]
\end{tikzcd}
\end{equation*}
where we employed the universality of $\T v_f$ and $v_{\T\^Af}$;
\item[canonical flip] The canonical flip $c\^A\colon(\T\^A)^2f\to(\T\^A)^2f$ is defined as:
\begin{equation*}
\begin{tikzcd}
& {(\T\^A)^2E} & {\T\T\^AE} & {\T^2E} & \T A \\
{(\T\^A)^2E} & {\T\T\^AE} & {\T^2E}
\arrow["{\T v_f}", from=1-3, to=1-4]
\arrow["{\T^2f}", shift left=2, from=1-4, to=1-5]
\arrow["{\T^2f\T p\T z}"', shift right=2, from=1-4, to=1-5]
\arrow["c"', from=2-3, to=1-4]
\arrow[dashed, from=2-3, to=1-3]
\arrow["{v_{\T\^Af}}", from=1-2, to=1-3]
\arrow["{\T v_f}"', from=2-2, to=2-3]
\arrow[dashed, from=2-2, to=1-2]
\arrow["{v_{\T\^Af}}"', from=2-1, to=2-2]
\arrow["{c\^A}", dashed, from=2-1, to=1-2]
\end{tikzcd}
\end{equation*}
where we employed the universality of $\T v_f$ and $v_{\T\^Af}$.
\end{description}
If $(\X,\TT)$ has negatives with negation $n$, then we can also lift the negation morphism to the slice tangent category as follows:
\begin{description}
\item[negation] The negation morphism $n\^A\colon\T\^Af\to\T\^Af$ is defined by:
\begin{equation*}
\begin{tikzcd}
{\T\^AE} & {\T E} & {\T A} \\
{\T\^AE} & {\T E} & {\T A}
\arrow["{v_f}", from=1-1, to=1-2]
\arrow["{\T f}", shift left=2, from=1-2, to=1-3]
\arrow["{\T fpz}"', shift right=2, from=1-2, to=1-3]
\arrow["{v_f}", from=2-1, to=2-2]
\arrow["{\T f}", shift left=2, from=2-2, to=2-3]
\arrow["{\T fpz}"', shift right=2, from=2-2, to=2-3]
\arrow[Rightarrow, no head, from=1-3, to=2-3]
\arrow["n"{description}, from=2-2, to=1-2]
\arrow["{n\^A}", dashed, from=2-1, to=1-1]
\end{tikzcd}
\end{equation*}
where we employed the universality of $v_f$.
\end{description}
We refer to this tangent category as the \textbf{slice tangent category} of $(\X,\TT)$ over $A$ and we denote it by $(\X,\TT)/A$.

\begin{remark}
\label{remark:slice-tangent-category-indexed-tangent-cat}
Given a sliceable tangent category $(\X,\TT)$ it is not hard to see that the operation $A\mapsto(\X,\TT)/A$ extends to a pseudofunctor $\Slice\^\TT\colon\X^\op\to\TngCat_\cong$. Cockett and Cruttwell in~\cite[Theorem~5.3]{cockett:differential-bundles} proved that the fibres of a tangent fibration (cf.~\cite[Definition~5.2]{cockett:differential-bundles}) are tangent categories and that the substitution functors are strong tangent morphisms. This result extends to a correspondence between tangent fibrations and pseudofunctors like $\X^\op\to\TngCat_\cong$. Interestingly, \cite[ Proposition~5.7]{cockett:differential-bundles} shows that $\Slice\^{\TT}$ is the pseudofunctor associated to a suitable tangent fibration.
\end{remark}

\subsection{The universal property of slicing}
\label{subsection:universality-slicing}
The goal of this subsection is to prove that the operation which takes a pair $(\X,\TT; A)$ formed by a tangent category (sliceable over $A$) and an object $A\in\X$ to its associated slice tangent category extends to a right adjoint of the functor that sends each tangent category $(\X,\TT)$ with terminal object $\*$ to the pair $((\X,\TT),\*)$. Let's start by introducing some useful jargon.

\begin{definition}
\label{definition:tangent-pair}
A \textbf{tangent pair} is a pair formed by a tangent category $(\X,\TT)$ sliceable over an object $A$ and the object $A$ itself. We denote a tangent pair by $(\X,\TT;A)$. Moreover, given two tangent pairs $(\X,\TT;A)$ and $({\X'},\TT';B)$, a \textbf{morphism of tangent pairs} $(F,\alpha;\varphi)\colon(\X,\TT;A)\to({\X'},\TT';B)$ is a lax tangent morphism $(F,\alpha)\colon(\X,\TT)\to({\X'},\TT')$ together with a morphism $\varphi\colon FA\to B$ of ${\X'}$.
\end{definition}

\par Tangent pairs together with their morphisms form a category denoted by $\TngPair$. In particular, notice that the composition of two morphisms of tangent pairs $(F,\alpha;\varphi)\colon(\X,\TT;A)\to({\X'},\TT';B)$ and $(G,\beta;\psi)\colon({\X'},\TT';B)\to({\X''},\TT'';C)$ is the triple formed by $G\o F\colon\X\to{\X''}$, the associated lax distributive law $G\o F\o\T\xrightarrow{G\alpha}G\o\T'\o F\xrightarrow{\beta_F}\T''\o G\o F$, and the morphism $G(F(A))\xrightarrow{G\varphi}GB\xrightarrow{\psi}C$ of ${\X''}$.

\begin{remark}
\label{remark:tangent-pair-category-of-elements}
Consider the pseudofunctor $\ITC\colon\TngCat\to\Cat$ which sends each tangent category $(\X,\TT)$ to the category of objects $A$ of $\X$ such that $(\X,\TT)$ is sliceable over $A$. Via the Grothendieck construction, this produces a cofibration $\int^{\TngCat}\ITC\to\TngCat$. The category of elements of this cofibration coincides with the category $\TngPair$ of tangent pairs.
\end{remark}

\begin{proposition}
\label{proposition:lifting-tangent-pair-morphisms-to-slice}
Consider two tangent pairs $(\X,\TT;A)$ and $({\X'},\TT';B)$ and a morphism of tangent pairs \linebreak$(F,\alpha;\varphi)\colon(\X,\TT;A)\to(\X',\TT';B)$. Let $f\colon E\to A$ a morphism in $\X$. Finally, consider the morphism $\theta_f\colon F\T\^AE\to{\T'}\^BFE$, as the unique morphism which makes commuting the following diagram:
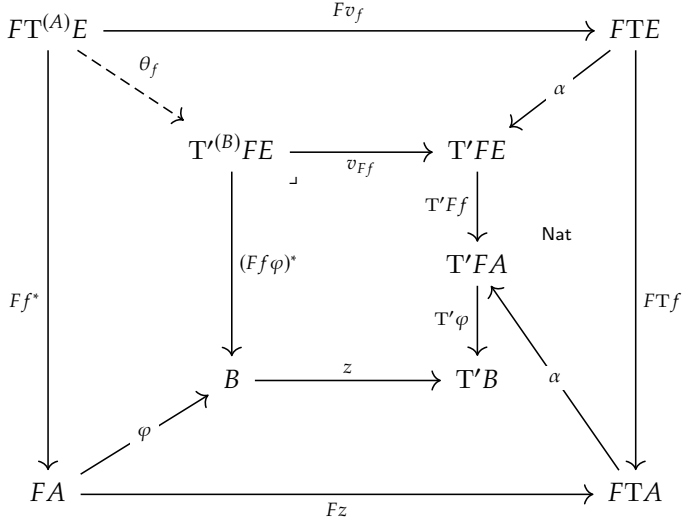
\begin{equation*}
\begin{tikzcd}
{F\T\^AE} &&&& {F\T E} \\
& {{\T'}\^BFE} && {\T' FE} \\
&&& {\T' FA} \\
& B && {\T' B} \\
FA &&&& {F\T A}
\arrow["{Fv_f}", from=1-1, to=1-5]
\arrow["{Ff^\*}"'{pos=0.6}, from=1-1, to=5-1]
\arrow["{(Ff\varphi)^\*}", from=2-2, to=4-2]
\arrow[""{name=0, anchor=center, inner sep=0}, "{\T' Ff}"', from=2-4, to=3-4]
\arrow["{v_{Ff}}"', from=2-2, to=2-4]
\arrow["\T'\varphi"', from=3-4, to=4-4]
\arrow["{\theta_f}", dashed, from=1-1, to=2-2]
\arrow["\alpha"{description}, from=1-5, to=2-4]
\arrow["z", from=4-2, to=4-4]
\arrow["\varphi"{description}, from=5-1, to=4-2]
\arrow["\lrcorner"{anchor=center, pos=0.125}, draw=none, from=2-2, to=3-4]
\arrow["Fz"', from=5-1, to=5-5]
\arrow[""{name=1, anchor=center, inner sep=0}, "{F\T f}"{pos=0.6}, from=1-5, to=5-5]
\arrow["\alpha"{description}, from=5-5, to=3-4]
\arrow["\Nat"{description}, draw=none, from=0, to=1]
\end{tikzcd}
\end{equation*}
Therefore, the functor:
\begin{align*}
&F\colon\X/A\to{\X'}/B\\
&F(f\colon E\to A)\mapsto(FE\xrightarrow{Ff}FA\xrightarrow{\varphi}B)\\
&F(g\colon(f\colon E\to A)\to(f'\colon E'\to A))\mapsto(Fg\colon(Ff\varphi)\to(Ff'\varphi))
\end{align*}
extends to a lax tangent morphism:
\begin{align*}
&(F,\alpha)/\varphi\colon(\X,\TT)/A\to({\X'},\TT')/B
\end{align*}
whose distributive law is defined by the natural transformation $\theta_f\colon F\T\^Af\to{\T'}\^BFf$.
\begin{proof}
For starters, let's prove the compatibility between $\theta$ and the projections:
\begin{equation*}
\begin{tikzcd}
{F\T\^Af} & {{\T'}\^BFf} \\
Ff
\arrow["{Fp\^A}"', from=1-1, to=2-1]
\arrow["{p\^B_F}", from=1-2, to=2-1]
\arrow["\theta", from=1-1, to=1-2]
\end{tikzcd}
\end{equation*}
which corresponds to the diagram:
\begin{equation*}
\begin{tikzcd}
{F\T\^AE} & {{\T'}\^BFE} \\
{F\T E} & {\T' FE} \\
FE & FE
\arrow[""{name=0, anchor=center, inner sep=0}, "\theta", from=1-1, to=1-2]
\arrow["{v_F}", from=1-2, to=2-2]
\arrow["{p_F}", from=2-2, to=3-2]
\arrow["Fv"', from=1-1, to=2-1]
\arrow["Fp"', from=2-1, to=3-1]
\arrow[""{name=1, anchor=center, inner sep=0}, Rightarrow, no head, from=3-1, to=3-2]
\arrow[""{name=2, anchor=center, inner sep=0}, "\alpha"{description}, from=2-1, to=2-2]
\arrow["{(\theta,\alpha;v)}"{description}, draw=none, from=0, to=2]
\arrow["{(\alpha;p)}"{description}, draw=none, from=1, to=2]
\end{tikzcd}
\end{equation*}
Let's take into consideration the compatibility diagram between $\theta$ and the zero morphisms:
\begin{equation*}
\begin{tikzcd}
{F\T\^Af} & {{\T'}\^BFf} \\
Ff
\arrow["\theta", from=1-1, to=1-2]
\arrow["{Fz\^A}", from=2-1, to=1-1]
\arrow["{z\^B_F}"', from=2-1, to=1-2]
\end{tikzcd}
\end{equation*}
To show that, first, consider the diagram:
\begin{equation*}
\begin{tikzcd}
{F\T\^AE} &&& {{\T'}\^BFE} \\
& {F\T E} & {\T' FE} \\
FE &&& FE
\arrow[""{name=0, anchor=center, inner sep=0}, "\theta", from=1-1, to=1-4]
\arrow["{v_F}"', from=1-4, to=2-3]
\arrow[""{name=1, anchor=center, inner sep=0}, "{Fz\^A}", from=3-1, to=1-1]
\arrow["Fv", from=1-1, to=2-2]
\arrow[""{name=2, anchor=center, inner sep=0}, "\alpha"', from=2-2, to=2-3]
\arrow["Fz"', from=3-1, to=2-2]
\arrow[""{name=3, anchor=center, inner sep=0}, Rightarrow, no head, from=3-1, to=3-4]
\arrow["{z_F}", from=3-4, to=2-3]
\arrow[""{name=4, anchor=center, inner sep=0}, "{z\^B_F}"', from=3-4, to=1-4]
\arrow["{(\alpha,\theta;v)}"{description}, draw=none, from=2, to=0]
\arrow["{(z;v)}"{description}, draw=none, from=2-3, to=4]
\arrow["{(z;v)}"{description}, draw=none, from=2-2, to=1]
\arrow["{(\alpha;z)}"{description}, draw=none, from=2, to=3]
\end{tikzcd}
\end{equation*}
Thus $Fz\^A\theta v_F=z\^B_Fv_F$ and from the universality of $v_F$ we conclude that $Fz\^A\theta=z\^B_F$, as expected. The next step is to prove the compatibility with the sum morphism:
\begin{equation*}
\begin{tikzcd}
{F\T\^A_2f} & {{\T'}\^B_2Ff} \\
{F\T\^Af} & {{\T'}\^BFf}
\arrow["\theta\times\theta", from=1-1, to=1-2]
\arrow["{Fs\^A}"', from=1-1, to=2-1]
\arrow["{s\^A_F}", from=1-2, to=2-2]
\arrow["\theta"', from=2-1, to=2-2]
\end{tikzcd}
\end{equation*}
Thus, consider the following diagram:
\begin{equation*}
\begin{tikzcd}
{F\T\^A_2E} &&& {{\T'}\^B_2FE} \\
& {F\T_2E} & {\T'_2FE} \\
& {F\T E} & {\T' FE} \\
{F\T\^AE} &&& {{\T'}\^BFE}
\arrow[""{name=0, anchor=center, inner sep=0}, "\theta"', from=4-1, to=4-4]
\arrow["{v_F}", from=4-4, to=3-3]
\arrow[""{name=1, anchor=center, inner sep=0}, "{Fs\^A}"', from=1-1, to=4-1]
\arrow["Fv"', from=4-1, to=3-2]
\arrow[""{name=2, anchor=center, inner sep=0}, "\alpha", from=3-2, to=3-3]
\arrow[""{name=3, anchor=center, inner sep=0}, "{s\^B_F}", from=1-4, to=4-4]
\arrow["{v_F\times v_F}"', from=1-4, to=2-3]
\arrow[""{name=4, anchor=center, inner sep=0}, "{s_F}", from=2-3, to=3-3]
\arrow["{Fv\times Fv}", from=1-1, to=2-2]
\arrow[""{name=5, anchor=center, inner sep=0}, "Fs"', from=2-2, to=3-2]
\arrow[""{name=6, anchor=center, inner sep=0}, "\alpha\times\alpha", from=2-2, to=2-3]
\arrow[""{name=7, anchor=center, inner sep=0}, "\theta\times\theta", from=1-1, to=1-4]
\arrow["{(s;v)}"{description}, draw=none, from=4, to=3]
\arrow["{(\alpha,\theta;v)}"{description}, draw=none, from=7, to=6]
\arrow["{(\alpha,\theta;v)}"{description}, draw=none, from=2, to=0]
\arrow["{(\alpha;s)}"{description}, draw=none, from=6, to=2]
\arrow["{(s;v)}"{description}, draw=none, from=1, to=5]
\end{tikzcd}
\end{equation*}
Thus, $Fs\^A\theta v_F=(\theta\times\theta)s\^B_Fv_F$ and from the universality of $v_F$ we conclude that $Fs\^A\theta=(\theta\times\theta)s\^B$, as expected. Let's prove the compatibility with the lift:
\begin{equation*}
\begin{tikzcd}
{F\T\^Af} && {{\T'}\^BFf} \\
{F(\T\^A)^2f} & {{\T'}\^BF\T\^Af} & {({\T'}\^B)^2Ff}
\arrow["\theta", from=1-1, to=1-3]
\arrow["{Fl\^A}"', from=1-1, to=2-1]
\arrow["{l\^B_F}", from=1-3, to=2-3]
\arrow["{\theta_{\T\^A}}"', from=2-1, to=2-2]
\arrow["{{\T'}\^B\theta}"', from=2-2, to=2-3]
\end{tikzcd}
\end{equation*}
As before, consider the following diagram:
\begin{equation*}
\begin{tikzcd}
{F\T\^AE} &&&& {{\T'}\^BFE} \\
& {F\T E} && {\T' FE} \\
& {F\T^2E} & {\T' F\T E} & {\T'^2FE} \\
& {F\T\^A\T E} & {{\T'}\^BF\T E} & {{\T'}\^B\T' FE} \\
{F(\T\^A)^2E} && {{\T'}\^BF\T\^AE} && {({\T'}\^B)^2FE}
\arrow[""{name=0, anchor=center, inner sep=0}, "\theta", from=1-1, to=1-5]
\arrow[""{name=1, anchor=center, inner sep=0}, "{l\^B_F}", from=1-5, to=5-5]
\arrow["{v_{\T' F}}"', from=4-4, to=3-4]
\arrow[""{name=2, anchor=center, inner sep=0}, "{Fl\^A}"', from=1-1, to=5-1]
\arrow["{v_F}"', from=1-5, to=2-4]
\arrow["{l_F}", from=2-4, to=3-4]
\arrow["Fv"', from=1-1, to=2-2]
\arrow[""{name=3, anchor=center, inner sep=0}, "\alpha", from=2-2, to=2-4]
\arrow["Fl"', from=2-2, to=3-2]
\arrow["{F\T\^Av}"{pos=0.3}, from=5-1, to=4-2]
\arrow["{Fv_\T}", from=4-2, to=3-2]
\arrow["{\alpha_\T}", from=3-2, to=3-3]
\arrow["\T'\alpha", from=3-3, to=3-4]
\arrow[""{name=4, anchor=center, inner sep=0}, "{\theta_{\T\^A}}"', from=5-1, to=5-3]
\arrow[""{name=5, anchor=center, inner sep=0}, "{{\T'}\^B\theta}"', from=5-3, to=5-5]
\arrow["{{\T'}\^Bv}"'{pos=0.3}, from=5-5, to=4-4]
\arrow["{{\T'}\^BFv}"{description}, from=5-3, to=4-3]
\arrow["{{\T'}\^B\alpha}", from=4-3, to=4-4]
\arrow["{\theta_\T}", from=4-2, to=4-3]
\arrow["{v_{F\T}}"{description}, from=4-3, to=3-3]
\arrow["\Nat"{description}, draw=none, from=3-3, to=4-4]
\arrow["{(\alpha,\theta;v)}"{description}, draw=none, from=3-3, to=4-2]
\arrow["{(\alpha,\theta;v)}"{description}, draw=none, from=0, to=3]
\arrow["{(\alpha;l)}"{description}, draw=none, from=3, to=3-3]
\arrow["{(\alpha,\theta;v)}"{description}, draw=none, from=5, to=4-4]
\arrow["\Nat"{description}, draw=none, from=4-2, to=4]
\arrow["{(l;v)}"{description}, draw=none, from=3-4, to=1]
\arrow["{(l;v)}"{description}, draw=none, from=3-2, to=2]
\end{tikzcd}
\end{equation*}
Therefore, $\theta l\^B_F{\T'}\^Bvv_{\T' F}=Fl\^A\theta_{\T\^A}{\T'}\^B\theta{\T'}\^Bvv_{\T' F}$. By the universality of ${\T'}\^Bvv_{\T' F}$ we conclude that $\theta l\^B_F=Fl\^A\theta_{\T\^A}{\T'}\^B\theta$, as expected. Finally, let's prove the compatibility with the canonical flip:
\begin{equation*}
\begin{tikzcd}
{F(\T\^A)^2f} & {{\T'}\^BF\T\^Af} & {({\T'}\^B)^2Ff} \\
{F(\T\^A)^2f} & {{\T'}\^BF\T\^Af} & {({\T'}\^B)^2Ff}
\arrow["{Fc\^A}"', from=1-1, to=2-1]
\arrow["{c\^B_F}", from=1-3, to=2-3]
\arrow["{\theta_{\T\^A}}"', from=2-1, to=2-2]
\arrow["{{\T'}\^B\theta}"', from=2-2, to=2-3]
\arrow["{\theta_{\T\^A}}", from=1-1, to=1-2]
\arrow["{{\T'}\^B\theta}", from=1-2, to=1-3]
\end{tikzcd}
\end{equation*}
Thus:
\begin{equation*}
\begin{tikzcd}
{F(\T\^A)^2E} && {{\T'}\^BF\T\^AE} && {({\T'}\^B)^2FE} \\
& {F\T\^A\T E} & {{\T'}\^BF\T E} & {{\T'}\^B\T' FE} \\
& {F\T E} & {\T' F\T E} & {\T'^2FE} \\
& {F\T^2E} & {\T' F\T E} & {\T'^2FE} \\
& {F\T\^A\T E} & {{\T'}\^BF\T E} & {{\T'}\^B\T' FE} \\
{F(\T\^A)^2E} && {{\T'}\^BF\T\^AE} && {({\T'}\^B)^2FE}
\arrow["{v_{\T' F}}"', from=5-4, to=4-4]
\arrow[""{name=0, anchor=center, inner sep=0}, "{Fc\^A}"', from=1-1, to=6-1]
\arrow["{c_F}", from=3-4, to=4-4]
\arrow["Fc"', from=3-2, to=4-2]
\arrow["{F\T\^Av}"{pos=0.3}, from=6-1, to=5-2]
\arrow["{Fv_\T}", from=5-2, to=4-2]
\arrow["{\alpha_\T}", from=4-2, to=4-3]
\arrow["\T'\alpha", from=4-3, to=4-4]
\arrow[""{name=1, anchor=center, inner sep=0}, "{\theta_{\T\^A}}"', from=6-1, to=6-3]
\arrow[""{name=2, anchor=center, inner sep=0}, "{{\T'}\^B\theta}"', from=6-3, to=6-5]
\arrow["{{\T'}\^Bv}"'{pos=0.3}, from=6-5, to=5-4]
\arrow["{{\T'}\^BFv}"{description}, from=6-3, to=5-3]
\arrow["{{\T'}\^B\alpha}", from=5-3, to=5-4]
\arrow["{\theta_\T}", from=5-2, to=5-3]
\arrow["{v_{F\T}}"{description}, from=5-3, to=4-3]
\arrow["\Nat"{description}, draw=none, from=4-3, to=5-4]
\arrow["{(\alpha,\theta;v)}"{description}, draw=none, from=4-3, to=5-2]
\arrow[""{name=3, anchor=center, inner sep=0}, "{c\^B_F}", from=1-5, to=6-5]
\arrow["{{\T'}\^Bv}"', from=1-5, to=2-4]
\arrow["{v_{\T' F}}", from=2-4, to=3-4]
\arrow["{F\T\^Av}"', from=1-1, to=2-2]
\arrow["{Fv_\T}"', from=2-2, to=3-2]
\arrow["{v_{F\T}}"{description}, from=2-3, to=3-3]
\arrow["{\alpha_\T}"', from=3-2, to=3-3]
\arrow["\T'\alpha"', from=3-3, to=3-4]
\arrow["{\theta_\T}"', from=2-2, to=2-3]
\arrow["{{\T'}\^B\alpha}"', from=2-3, to=2-4]
\arrow["{{\T'}\^BFv}"{description}, from=1-3, to=2-3]
\arrow[""{name=4, anchor=center, inner sep=0}, "{\theta_{\T\^A}}", from=1-1, to=1-3]
\arrow[""{name=5, anchor=center, inner sep=0}, "{{\T'}\^B\theta}", from=1-3, to=1-5]
\arrow["{(\alpha,\theta;v)}"{description}, draw=none, from=3-3, to=2-2]
\arrow["\Nat"{description}, draw=none, from=3-3, to=2-4]
\arrow["{(\alpha;c)}"{description}, draw=none, from=3-3, to=4-3]
\arrow["{(\alpha,\theta;v)}"{description}, draw=none, from=2, to=5-4]
\arrow["\Nat"{description}, draw=none, from=5-2, to=1]
\arrow["{(c;v)}"{description}, draw=none, from=4-2, to=0]
\arrow["{(c;v)}"{description}, draw=none, from=4-4, to=3]
\arrow["\Nat"{description}, draw=none, from=2-2, to=4]
\arrow["{(\alpha,\theta;v)}"{description}, draw=none, from=2-4, to=5]
\end{tikzcd}
\end{equation*}
This proves that $\theta_{\T\^A}{\T'}\^B\theta c\^B_F{\T'}\^Bvv_{\T' F}=Fc\^A\theta_{\T\^A}{\T'}\^B\theta{\T'}\^Bvv_{\T' F}$. Finally, using the universality of ${\T'}\^Bvv_{\T' F}$ we conclude that $\theta_{\T\^A}{\T'}\^B\theta c\^B_F=Fc\^A\theta_{\T\^A}{\T'}\^B\theta$, as expected.
\end{proof}
\end{proposition}

Proposition~\ref{proposition:lifting-tangent-pair-morphisms-to-slice} allows us to lift morphisms of tangent pairs to the corresponding slice tangent categories. The next step is to find sufficient conditions so that the corresponding tangent morphism over the slice categories is strong. This will play a key role in the next section. Let's introduce a definition.

\begin{definition}
\label{definition:cartesian-morphisms-of-tangent-pairs}
Given two tangent pairs $(\X,\TT;A)$ and $({\X'},\TT';B)$, a morphism of tangent pairs $(F,\alpha;\varphi)\colon(\X,\TT;A)\to({\X'},\TT';B)$ is \textbf{Cartesian} if the following diagrams:
\begin{equation*}
\begin{tikzcd}
{F\T E} & {\T' FE} \\
{F\T A} & {\T' FA}
\arrow["{F\T f}"', from=1-1, to=2-1]
\arrow["{\T' Ff}", from=1-2, to=2-2]
\arrow["\alpha"', from=2-1, to=2-2]
\arrow["\alpha", from=1-1, to=1-2]
\end{tikzcd}\hfill
\begin{tikzcd}
FA & {\T' FA} \\
B & {\T' B}
\arrow["\varphi"', from=1-1, to=2-1]
\arrow["\T'\varphi", from=1-2, to=2-2]
\arrow["z"', from=2-1, to=2-2]
\arrow["{z_F}", from=1-1, to=1-2]
\end{tikzcd}
\end{equation*}
are pullback diagrams, and moreover the functor $F$ preserves the pullbacks of Equation~\eqref{equation:sliceable-pullback}. Concretely, this last condition means that for every morphism $f\colon E\to A$ of $\X$, the diagram:
\begin{equation*}
\begin{tikzcd}
{F\T\^AE} & {F\T E} \\
FA & {F\T A}
\arrow["{F\T f}", from=1-2, to=2-2]
\arrow["Fz"', from=2-1, to=2-2]
\arrow["{Fv_f}", from=1-1, to=1-2]
\arrow["{Ff^\*}"', from=1-1, to=2-1]
\end{tikzcd}
\end{equation*}
must be a pullback diagram.
\end{definition}

\begin{lemma}
\label{lemma:cartesian-morphisms-become-strong}
A Cartesian morphism of tangent pairs $(F,\alpha;\varphi)\colon(\X,\TT;A)\to({\X'},\TT';B)$ lifts to the slice tangent categories as a strong tangent morphism. Concretely, this means that the natural transformation $\theta_f\colon F\T\^Af\to{\T'}\^BFf$ is invertible.
\begin{proof}
Consider the following diagram:
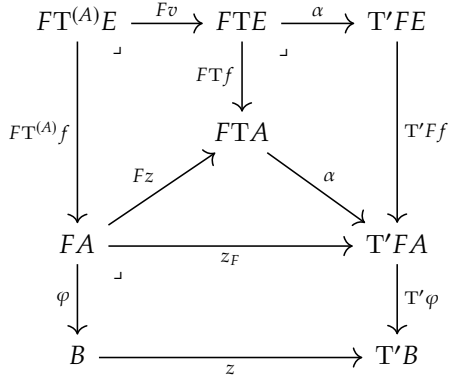
\begin{equation*}
\begin{tikzcd}
{F\T\^AE} & {F\T E} & {\T' FE} \\
& {F\T A} \\
FA && {\T' FA} \\
B && {\T' B}
\arrow["Fv", from=1-1, to=1-2]
\arrow["\alpha", from=1-2, to=1-3]
\arrow["{F\T\^Af}"', from=1-1, to=3-1]
\arrow["Fz", from=3-1, to=2-2]
\arrow["{F\T f}"', from=1-2, to=2-2]
\arrow["\alpha", from=2-2, to=3-3]
\arrow[""{name=0, anchor=center, inner sep=0}, "{\T' Ff}", from=1-3, to=3-3]
\arrow["{z_F}"', from=3-1, to=3-3]
\arrow["\varphi"', from=3-1, to=4-1]
\arrow[""{name=1, anchor=center, inner sep=0}, "z"', from=4-1, to=4-3]
\arrow["\T'\varphi", from=3-3, to=4-3]
\arrow["\lrcorner"{anchor=center, pos=0.125}, draw=none, from=1-1, to=2-2]
\arrow["\lrcorner"{anchor=center, pos=0.125}, draw=none, from=1-2, to=0]
\arrow["\lrcorner"{anchor=center, pos=0.125}, draw=none, from=3-1, to=1]
\end{tikzcd}
\end{equation*}
where we used that $Fz\alpha=z_F$. Thanks to the Cartesianity of $(F,\alpha;\varphi)$ this is a pullback diagram, since it is formed by pullback diagrams. However, by definition, $\theta$ is defined by the diagram:
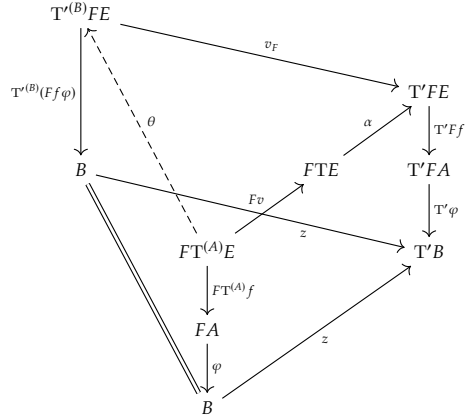
\begin{equation*}
\adjustbox{scale=.7,center}{
\begin{tikzcd}
{{\T'}\^BFE} \\
&&& {\T' FE} \\
B && {F\T E} & {\T' FA} \\
& {F\T\^AE} && {\T' B} \\
& FA \\
& B
\arrow["Fv", from=4-2, to=3-3]
\arrow["\alpha", from=3-3, to=2-4]
\arrow["{F\T\^Af}", from=4-2, to=5-2]
\arrow["{\T' Ff}", from=2-4, to=3-4]
\arrow["\varphi", from=5-2, to=6-2]
\arrow["z"', from=6-2, to=4-4]
\arrow["\T'\varphi", from=3-4, to=4-4]
\arrow["{v_F}", from=1-1, to=2-4]
\arrow["{{\T'}\^B(Ff\varphi)}"', from=1-1, to=3-1]
\arrow["z"'{pos=0.7}, from=3-1, to=4-4]
\arrow["\theta"', dashed, from=4-2, to=1-1]
\arrow[Rightarrow, no head, from=6-2, to=3-1]
\end{tikzcd}
}
\end{equation*}
However, the top and the right rectangular sides of this triangular diagram are pullbacks, so $\theta$ must be an isomorphism.
\end{proof}
\end{lemma}

The next is a key concept for our discussion.

\begin{definition}
\label{definition:tangent-category-terminal-object}
A tangent category \textbf{with terminal object} is a tangent category $(\X,\TT)$ equipped with a terminal object $\*$ so that the unique morphism $\T\*\to\*$ is an isomorphism. We also denote by $\trmTngCat$ the $2$-category of tangent categories with terminal objects, lax tangent morphisms and natural transformations compatible with the lax distributive laws.
\end{definition}

In the following, we denote by $\*$ the (unique up to unique isomorphism) terminal object of a category. Moreover, for any object $A$, the unique morphism from $A$ to $\*$ is denoted by $!\colon A\to\*$. It is straightforward to see that a tangent category $(\X,\TT)$ with terminal object is sliceable over $\*$ and that the slice tangent category $(\X,\TT)/\ast$ is isomorphic to $(\X,\TT)$ via $(!\colon A\to\*)\mapsto A$. This observation allows us to define the following pseudofunctor:
\begin{align*}
&\Term\colon\trmTngCat\to\TngPair\\
&\Term(\X,\TT)\=(\X,\TT;\ast)\\
&\Term((F,\alpha)\colon(\X,\TT)\to({\X'},\TT'))\=(F,\alpha;F\*\xrightarrow{!}\*)\colon(\X,\TT;\ast)\to({\X'},\TT';\ast)
\end{align*}
Thanks to Proposition~\ref{proposition:lifting-tangent-pair-morphisms-to-slice}, the operation which takes a tangent pair $(\X,\TT;A)$ to its slice tangent category extends to a functor. Observe that the slice tangent category of a tangent pair $(\X,\TT;A)$ is equipped with a terminal object, the terminal object being the identity over $A$. With this in mind, we are able to define the following pseudofunctor:
\begin{align*}
&\Slice\colon\TngPair\to\trmTngCat\\
&\Slice(\X,\TT;A)\mapsto(\X,\TT)/A\\
&\Slice((F,\alpha;\varphi)\colon(\X,\TT;A)\to({\X'},\TT';B))\=(F,\alpha)/\varphi\colon(\X,\TT)/A\to({\X'},\TT')/B
\end{align*}

\begin{remark}
\label{remark:pseudofunctoriality-slice-term}
$\Term$ and $\Slice$ are not strict functors but rather pseudofunctors. This comes from the fact that terminal objects and slice tangent structures are defined only up to unique isomorphisms. Thus, the associators and unitors are defined by these unique isomorphisms.
\end{remark}

We can finally characterize the operation which takes a tangent pair to its slice tangent category as an adjunction between pseudofunctors.

\begin{theorem}
\label{theorem:adjunction-Term-Slice}
The pseudofunctors $\Slice\colon\TngPair\leftrightarrows\trmTngCat\colon\Term$ form an adjunction whose left adjoint is $\Term$, the right adjoint is $\Slice$, the unit $(U,\eta)\colon(\X,\TT)\to\Slice(\Term(\X,\TT))=(\X,\TT)/\ast$, as a lax tangent morphism between tangent categories with terminal objects, is the isomorphism:
\begin{align*}
&U\colon\X\to\X/\ast\\
&U(A)\mapsto(!\colon A\to\*)\\
&U(f\colon A\to B)\mapsto(f(!\colon A\to\*)\to(!\colon B\to\*))\\
&\eta\colon(U(\T A))=(!\colon\T A\to\*)\xrightarrow{\id_{\T A}}(!\colon\T A\to\*)=\T(U(A))
\end{align*}
and the counit $(C,\epsilon;\varphi)\colon\Term(\Slice(\X,\TT;A))=((\X,\TT)/A,\id_A)\to(\X,\TT;A)$ is the morphism of tangent pairs:
\begin{align*}
&C\colon(\X,\TT)/A\mapsto(\X,\TT)\\
&C(f\colon E\to A)\mapsto E\\
&C(g\colon(f\colon E\to A)\to(f'\colon E'\to A))\mapsto(g\colon E\to E')\\
&\epsilon\colon C(\T\^A(f\colon E\to A))=\T\^AE\xrightarrow{v_f}\T E=\T(C(f\colon E\to A))\\
&\varphi\colon C(\id_A\colon A\to A)=A\xrightarrow{\id_A}A
\end{align*}
\begin{proof}
To prove the result we need to show that the unit $(U,\eta)$ and the counit $(C,\epsilon;\varphi)$ fulfill the triangle identities. Let's start by considering the following diagram:
\begin{equation*}
\begin{tikzcd}
{\Term(\X,\TT)} & {\Term(\Slice(\Term(\X,\TT)))} \\
& {\Term(\X,\TT)}
\arrow["{\Term(U,\eta)}", from=1-1, to=1-2]
\arrow["{(C,\epsilon;\varphi)_\Term}", from=1-2, to=2-2]
\arrow[Rightarrow, no head, from=1-1, to=2-2]
\end{tikzcd}
\end{equation*}
for a tangent category $(\X,\TT)$ with terminal object. However, it is straightforward to realize that the underlying tangent morphisms $(C,\epsilon)$ and $(U,\eta)$ of $(C,\epsilon;\varphi)_{\Term}$ and $\Term(U,\eta)$ define the equivalence between $(\X,\TT)$ and $(\X,\TT)/\*$ and that, by the universality of the terminal object, that the composition of the comparison morphisms $\varphi=\id_\*$ and $!\colon U\*\to\*$ is the identity over the terminal object. Similarly, by considering the diagram:
\begin{equation*}
\begin{tikzcd}
{\Slice(\X,\TT;A)} & {\Slice(\Term(\Slice(\X,\TT;A)))} \\
& {\Slice(\X,\TT;A)}
\arrow["{(U,\eta)_\Slice}", from=1-1, to=1-2]
\arrow["{\Slice(C,\epsilon;\varphi)}", from=1-2, to=2-2]
\arrow[Rightarrow, no head, from=1-1, to=2-2]
\end{tikzcd}
\end{equation*}
for a tangent pair $(\X,\TT;A)$, it is straightforward to show the underlying tangent morphisms of $\Slice(C,\epsilon;\varphi)$ and $(U,\eta)_{\Slice}$ define the equivalence between $(\X,\TT)/A$ and $((\X,\TT)/A)/\id_A$ and that the composition of the comparison morphisms gives the identity. Finally, notice that the unit is always an isomorphism.
\end{proof}
\end{theorem}

\section{The slice tangent categories of the affine schemes over an operad}
\label{section:slicing-operadic-affine-schemes}
The previous section was dedicated to characterizing the slicing of tangent categories via the adjunction between two pseudofunctors. A similar phenomenon happens in the realm of operads: given an operad $\P$ and a $\P$-algebra $A$ the enveloping operad $\P\^A$ of $\P$ over $A$ is the operad whose category of algebras is equivalent to the coslice category of $\Alg_\P$ under $A$.
\par The goal of this section is to prove that these two phenomena are two faces of the same coin: the geometric tangent category of the enveloping operad of $\P$ over $A$ is equivalent to the slice tangent category of the geometric tangent category of $\P$ over $A$.
\par Let's start by recalling the definition of the enveloping operad of a pair $(\P;A)$. We advise the interested reader to consult \cite{moerdijk:enveloping-operads}, \cite{loday:operads}, or \cite{fresse:operads}. For this purpose, recall that since the category of algebras of an operad $\P$ is cocomplete, each operad has an initial algebra, which corresponds to the $R$-module $\P(0)$ together with structure map $\P(m)\x\P(0)^{\x m}\to\P(0)$ defined by the operadic composition. This allows us to introduce an operation $\P\mapsto(\P;\P(0))$ between operads and operadic pairs. Notice that for a \textbf{operadic pair} we mean a pair $(\P;A)$ formed by an operad $\P$ and a $\P$-algebra $A$. Moreover, given two operadic pairs $(\P;A)$ and $(\OprQ;B)$ a \textbf{morphism of operadic pairs} $(f;\varphi)\colon(\P;A)\to(\OprQ;B)$ is a morphism of operads $f\colon\P\to\OprQ$ together with a morphism of $\P$-algebras $\varphi\colon A\to f^\*B$, $f^\*\colon\Alg_\OprQ\to\Alg_\P$ being the pullback functor induced by $f$. Operadic pairs together with their morphisms form a category that we denote by $\OprPair$. So, we have:
\begin{align*}
&\Init\colon\Operad\to\OprPair\\
&\Init(\P)\=(\P;\P(0))\\
&\Init(f\colon\P\mapsto\OprQ)\=(f;!\colon\P(0)\to f^\*\OprQ(0))
\end{align*}
$!$ being the unique morphism of $\P$-algebras induced by the universality of the initial algebra $\P(0)$. Concretely, $!$ sends an element $u\in\P(0)$ to $f_0(u)$.

\begin{remark}
\label{remark:operadic-pairs-category-of-elements}
Similarly as for tangent pairs (see Remark~\ref{remark:tangent-pair-category-of-elements}), also operadic pairs can be regarded as objects in the category of elements of a fibration. Consider the pseudofunctor $\ITC\colon\Operad^\op\to\Cat$ which sends each operad to the corresponding category of algebras. Via the Grothendieck construction, this is equivalent to a fibration $\int_{\Operad}\ITC\to\Operad$ and the category of elements $\int_{\Operad}\ITC$ is equivalent to the category $\OprPair$ of operadic pairs.
\end{remark}

\par $\Init$ admits a left adjoint $\Env$ (cf.~\cite{moerdijk:enveloping-operads}), which sends an operadic pair $(\P;A)$ to the corresponding enveloping operad $\Env(\P;A)\=\P\^A$. Following the description provided by~\cite[Section~4.1.3]{fresse:operads}, $\P\^A$ is generated by the symbols $(\mu;a_1\,a_k|$, for every $\mu\in\P(m+k)$, $a_1\,a_k\in A$ and every non-negative integer $k$ (when $k=0$, $(\mu|$ are the only terms) which fulfill the following relations:
\begin{align}
\label{equation:enveloping-operad-relations}
&(\mu;a_1\,\nu(a_i\,a_{i+n})\,a_{k+n}|=(\mu\o_i\nu;a_1\,a_{k+n}|
\end{align}
for $\mu\in\P(m+k)$, $\nu\in\P(n)$ and $a_1\,a_{k+n}\in A$, where we used the notation $\mu\o_i\nu$ for $\mu(1_\P\,\nu\,1_\P)$. In particular, it is not hard to see that $\P\^A(0)\cong A$. So, the functor $\Env$ sends a morphism of operadic pairs $(f,\varphi)\colon(\P;A)\to(\OprQ;B)$ to the morphism of operads $\Env(f;\varphi)\colon\P\^A\to\OprQ\^B$ defined on generators as follows:
\begin{align*}
&(\mu;a_1\,a_k|\mapsto(f(\mu);\varphi(a_1)\,\varphi(a_k)|
\end{align*}
From this description of the enveloping operad, it is not hard to see that an algebra $B$ of the enveloping operad $\P\^A$ is precisely given by a $\P$-algebra $C^\*B$, $C\colon\P\to\P\^A$ being the canonical inclusion $\mu\mapsto(\mu|$, together with a  morphism of $\P$-algebras $A\to C^\*B$ induced by the structure map $A=\P\^A(0)\to C^\*B$ of $B$.
\par Conversely, every morphism of $\P$-algebras $f\colon A\to B$ induces over $B$ a $\P\^A$-algebra structure defined as follows:
\begin{align*}
&(\mu;a_1\,a_k|(b_1\,b_m)\=\mu_B(f(a_1)\,f(a_k),b_1\,b_m)
\end{align*}
for $\mu\in\P(m+k)$, $a_1\,a_k\in A$ and $b_1\,b_m\in B$. This proves that the category of $\P\^A$-algebras is equivalent to the coslice category of $\P$-algebras over $A$ (cf.~\cite[Lemma~1.7]{moerdijk:enveloping-operads}).

\subsection{The geometric tangent category of the enveloping operad}
\label{subsection:enveloping-tangent-category}
Theorem~\ref{theorem:adjunction-Term-Slice} establishes that $\Term$ and $\Slice$ form an adjunction and from our discussion on the enveloping operad we also know that also $\Env$ and $\Init$ form an adjunction. We would like to compare $\Term$ with $\Init$ and $\Slice$ with $\Env$. However, $\Term$ is a left adjoint, while $\Init$ is a right adjoint and similarly, $\Slice$ is a right adjoint and $\Env$ is a left adjoint. To solve this issue, we transpose the adjunction $\Env\dashv\Init$ to the opposite categories. To compare these functors, notice that $\Geom^\*$ extends to operadic pairs as follows:
\begin{align*}
&\Geom^\*\colon\OprPair^\op\to\TngPair\\
&\Geom^\*(\P;A)\=(\Geom(\P);A)\\
&\Geom^\*((f,\varphi)\colon(\P;A)\to(\OprQ;B))\=(\Geom^\*(f)=(f^\*,\alpha^\*);\varphi^\op\colon A\leftarrow\varphi^\*B)\colon(\Geom(\P);A)\to(\Geom(\OprQ);B)
\end{align*}
Note that, since $\Alg_\P$ is cocomplete, $\Geom(\P)$ is sliceable.

\begin{lemma}
\label{lemma:init-is-term}
The following diagram:
\begin{equation*}
\begin{tikzcd}
{\Operad^\op} & {\OprPair^\op} \\
\trmTngCat & \TngPair
\arrow["{\Geom^\*}"', from=1-1, to=2-1]
\arrow["{\Geom^\*}", from=1-2, to=2-2]
\arrow["\Init", from=1-1, to=1-2]
\arrow["\Term"', from=2-1, to=2-2]
\end{tikzcd}
\end{equation*}
commutes.
\begin{proof}
It is straightforward to see that, for an operad $\P$:
\begin{align*}
&\Geom^\*(\Init((\P))=(\Geom(\P);\P(0))=\Term(\Geom(\P))=\Term(\Geom^\*(\P))
\end{align*}
and for a morphism of operads $f\colon\P\to\OprQ$:
\begin{align*}
&\Geom^\*(\Init(f))=\Geom^\*(f;!\colon f^\*\OprQ(0)\leftarrow\P(0))=\\
&=(f^\*,\alpha^\*;!\colon\P(0)\to f^\*\OprQ(0))=\Term(f^\*,\alpha^\*)=\Term(\Geom^\*(f))
\end{align*}
So, the diagram commutes.
\end{proof}
\end{lemma}

Thanks to Lemma~\ref{lemma:init-is-term} we can now also compare the functors $\Env$ and $\Slice$. Crucially, to do that we are going to use that $\Init\dashv\Env$ (on the opposite categories) and that $\Term\dashv\Slice$ form adjunctions. In general, given a square diagram as follows:
\begin{equation*}
\begin{tikzcd}
\bullet & \bullet \\
\bullet & \bullet
\arrow["{F'}", shift left=2, from=1-1, to=1-2]
\arrow["{U'}", shift left=2, from=1-2, to=1-1]
\arrow["F", shift left=2, from=2-1, to=2-2]
\arrow["U", shift left=2, from=2-2, to=2-1]
\arrow["G"', from=1-1, to=2-1]
\arrow["H", from=1-2, to=2-2]
\end{tikzcd}
\end{equation*}
with $(\eta,\epsilon)\colon F\dashv U$ and $(\eta',\epsilon')\colon F'\dashv U'$ forming adjunctions, then if the diagram:
\begin{equation*}
\begin{tikzcd}
\bullet & \bullet \\
\bullet & \bullet
\arrow["{F'}", from=1-1, to=1-2]
\arrow["F", from=2-1, to=2-2]
\arrow["G"', from=1-1, to=2-1]
\arrow["H", from=1-2, to=2-2]
\end{tikzcd}
\end{equation*}
commutes, then, by using mates, we can define the following natural transformation:
\begin{align*}
&G\o U'\xrightarrow{\eta_{GU'}}U\o F\o G\o U'=U\o H\o F'\o U'\xrightarrow{UH\epsilon'}U\o H
\end{align*}
A priori, there is no reason to conclude that such a natural transformation is a natural isomorphism. In order to prove that the natural transformation induced by the adjunctions $\Init\dashv\Env$, $\Term\dashv\Slice$, and by Lemma~\ref{lemma:init-is-term} is an isomorphism, we need to show that the counit of $\Init\dashv\Env$ induces a Cartesian morphism of tangent pairs over the geometric tangent pairs.

\begin{lemma}
\label{lemma:cartesianity-counit}
The counit, regarded as a morphism of $\OprPair$, $(C,\epsilon)\colon(\P;A)\to\Init(\Env(\P;A))=(\P^A;\P\^A(0))$ of the adjunction $\Init\dashv\Env$ induces a Cartesian morphism of tangent pairs:
\begin{align*}
&\Geom^\*(C,\epsilon)\colon\Geom^\*(\P\^A;\P\^A(0))\to\Geom^\*(\P;A)
\end{align*}
\begin{proof}
Let's start by recalling the definition of the counit. $C\colon\P\to\P\^A$ is the morphism of operads which includes $\P$ into $\P\^A$ by mapping $\mu\in\P(m)$ into $(\mu|\in\P\^A(m)$. Moreover, $\epsilon\colon A\to C^\*\P\^A(0)$ is the isomorphism $A\ni a\mapsto(1_\P;a|\in C^\*\P\^A(0)$, where $1_\P\in\P(1)$ is the unit of $\P$. To see that this is an isomorphism, notice that the generators of $\P\^A(0)$ are all the symbols $(\mu;a_1\,a_m|$ for every $\mu\in\P(m)$ and $a_1\,a_m\in A$, but thanks to the relations~\eqref{equation:enveloping-operad-relations} we also have:
\begin{align*}
&(\mu;a_1\,a_m|=(1_\P(\mu);a_1\,a_m|=(1_\P;\mu(a_1\,a_m)|
\end{align*}
So, with the identification $a=(1_\P;a|$ we have that $\P\^A(0)$ is equal to $A$. Notice also that, given a $\P\^A$-algebra $B$, $C^\*B$ is the $\P$-algebra over $B$ with structure map defined by:
\begin{align*}
&\mu(b_1\,b_m)\=(\mu|_B(b_1\,b_m)
\end{align*}
\par To distinguish between the different tangent structures, for this proof we adopt the following convention: we denote by $\TT$ the geometric tangent structure of $\P$, by $\TT\^B$ the slice tangent structure over $B$, and by $\TT_A$ the geometric tangent structure of $\P\^A$.
\par The Cartesianity of $\Geom^\*(C,\epsilon)$ means that for a morphism $f\colon B\to E$ of $\P\^A$-algebras the diagrams in the category of $\P$-algebras:
\begin{equation*}
\begin{tikzcd}
{C^\*T_AB} & {C^\*T_AE} \\
{C^\*B} & {C^\*(\T_A)\^BE}
\arrow["{C^\*z_A}"', from=1-1, to=2-1]
\arrow["{C^\*v_f}", from=1-2, to=2-2]
\arrow["{C^\*\T_Af}", from=1-1, to=1-2]
\arrow["{C^\*f_\ast}"', from=2-1, to=2-2]
\end{tikzcd}\hfill
\begin{tikzcd}
{\T C^\*B} & {C^\*\T_AB} \\
{\T C^\*E} & {C^\*\T_AE}
\arrow["{\alpha^\*}", from=1-1, to=1-2]
\arrow["{\alpha^\*}"', from=2-1, to=2-2]
\arrow["{\T C\*f}"', from=1-1, to=2-1]
\arrow["{C^\*\T_Af}", from=1-2, to=2-2]
\end{tikzcd}\hfill
\begin{tikzcd}
{\T A} & A \\
{\T C^\*\P\^A(0)} & {C^\*\P\^A(0)}
\arrow["z", from=1-1, to=1-2]
\arrow["\epsilon", from=1-2, to=2-2]
\arrow["\T\epsilon"', from=1-1, to=2-1]
\arrow["{z_{C^\*}}"', from=2-1, to=2-2]
\end{tikzcd}
\end{equation*}
are all pushout diagrams, where $f_\*$ is the morphism defined by the pushout diagram in $\Alg_{\P\^A}$:
\begin{equation*}
\begin{tikzcd}
{T_AB} & {T_AE} \\
B & {(\T_A)\^BE}
\arrow["{z_A}"', from=1-1, to=2-1]
\arrow["{v_f}", from=1-2, to=2-2]
\arrow["{\T_Af}", from=1-1, to=1-2]
\arrow["{f_\ast}"', from=2-1, to=2-2]
\arrow["\lrcorner"{anchor=center, pos=0.125, rotate=180}, draw=none, from=2-2, to=1-1]
\end{tikzcd}
\end{equation*}
Notice that the third diagram is trivially a pushout since $\epsilon$ is an isomorphism. Let's focus on the first diagram and let's consider two morphisms $g\colon C^\*\T E\to K$ and $h\colon C^\*A\to K$ of $\P$-algebras satisfying the commutativity of the diagram:
\begin{equation*}
\begin{tikzcd}
{C^\*T_AB} & {C^\*T_AE} \\
{C^\*B} & {C^\*(\T_A)\^BE} \\
&& {K}
\arrow["{C^\*z_A}"', from=1-1, to=2-1]
\arrow["{C^\*v_f}", from=1-2, to=2-2]
\arrow["{C^\*\T_Af}", from=1-1, to=1-2]
\arrow["{C^\*f_\ast}"', from=2-1, to=2-2]
\arrow["g"', bend right, from=2-1, to=3-3]
\arrow["h", bend left, from=1-2, to=3-3]
\end{tikzcd}
\end{equation*}
So, we have that:
\begin{align*}
&g(b)=h(f(b))\\
&h(\d f(b))=0
\end{align*}
for every $b\in B$. Recall that $\P\^A$-algebras are equivalent to $\P$-algebra morphisms with $A$ for domain. Since $B$ is a $\P\^A$-algebra, we obtain a $\P$-algebra morphism $\beta\colon A\to C^\*B$. So, by post-composing by $g$ we get a new $\P$-algebra morphism $A\xrightarrow{\beta}C^\*B\xrightarrow{g}K$, thus, we get a $\P\^A$-algebra structure over $K$. Concretely, the structure map of this $\P\^A$-algebra $\bar K$ is defined by:
\begin{align*}
&(\mu;a_1\,a_k|_{\bar K}(x_1\,x_m)\=\mu_K(g(\beta(a_1)\,g(\beta(a_k)),x_1\,x_m)
\end{align*}
Moreover, we can also lift $g$ to a morphism of $\P\^A$-algebras $\bar g\colon B\to\bar K$, defined simply by $b\mapsto g(b)$. Let's now do the same for $h$: define a morphism of $\P\^A$-algebras $\bar h\colon\T_AE\to\bar K$. To do so, note also that $\T_AE$ is a $\P\^A$-algebra, which corresponds to a morphism of $\P$-algebras $\gamma\colon A\to C^\*\T_AE$ but because $f$ is a morphism of $\P\^A$-algebras, we have that for any $a\in A$:
\begin{align*}
&\gamma(a)=p_A(f(\beta(a)))=f(\beta(a))
\end{align*}
where we used that the projection $E\to\T_AE$ sends each element $x\in E$ to itself, $\T_AE$ being generated by all $x$ and all $\d x$. This implies that we can define $\bar h$ as the morphism which sends each $y\in\T_AE$ to $h(y)$. To see that this is a morphism of $\P\^A$-algebras note that:
\begin{eqnarray*}
& &(\mu;a_1\,a_k|_{\bar K}(\bar h(y_1)\,\bar h(y_m))\\
&=&\mu_K(g(\beta(a_1))\,g(\beta(a_k)),h(y_1)\,h(y_m))\\
&=&\mu_K(h(f(\beta(a_1)))\,h(f(\beta(a_k))),h(y_1)\,h(y_m))\\
&=&h(\mu_{C^\*\T_AE}(f(\beta(a_1))\,f(\beta(a_k)),y_1\,y_m))\\
&=&h(\mu_{C^\*\T_AE}(\gamma(a_1)\,\gamma(a_k),y_1\,y_m))\\
&=&h((\mu;a_1\,a_k|_{\T_AE}(y_1\,y_m))
\end{eqnarray*}
where we used that $g(b)=h(f(b))$, for any $b\in B$. Moreover, note that $C^\*\bar K=K$, $C^\*\bar g=g$ and that $C^\*\bar h=h$. Thus, we now have the following commutative diagram in $\Alg_{\P\^A}$:
\begin{equation*}
\begin{tikzcd}
{T_AB} & {T_AE} \\
B & {(\T_A)\^BE} \\
&& {\bar K}
\arrow["{z_A}"', from=1-1, to=2-1]
\arrow["{v_f}", from=1-2, to=2-2]
\arrow["{\T_Af}", from=1-1, to=1-2]
\arrow["{f_\ast}"', from=2-1, to=2-2]
\arrow["{\bar g}"', bend right, from=2-1, to=3-3]
\arrow["{\bar h}", bend left, from=1-2, to=3-3]
\arrow["\lrcorner"{anchor=center, pos=0.125, rotate=180}, draw=none, from=2-2, to=1-1]
\arrow["{\overline{[g,h]}}"{description}, dashed, from=2-2, to=3-3]
\end{tikzcd}
\end{equation*}
To see that recall that $g(b)=h(f(b))$ and that $h(\d f(b))=0$, which precisely implies the commutativity of this diagram. Therefore, we have a unique morphism $\overline{[g,h]}\colon(\T_A)\^BE\to\bar K$ of $\P\^A$-algebras. So, let's introduce:
\begin{align*}
&[g,h]\=C^\*(\overline{[g,h]})\colon C^\*(\T_A)\^BE\to C^\*\bar K=K
\end{align*}
However:
\begin{align*}
&C^\*f[g,h]=C^\*fC^\*(\overline{[g,h]})=C^\*(f\overline{[g,h]})=C^\*\bar g=g\\
&C^\*v_f[g,h]=C^\*v_fC^\*(\overline{[g,h]})=C^\*(v_f\overline{[g,h]})=C^\*\bar h=h
\end{align*}
Finally, suppose that $r\colon C^\*(\T_A)\^BE\to K$ is a second morphism of $\P$-algebra such that $(C^\*f)r=g$ and $(C^\*v_f)r=h$. However, in a similar fashion we can also lift $r$ to a morphism of $\P\^A$-algebras $\bar r\colon(\T_A)\^BE\to\bar K$ such that $C^\*\bar r=r$. But this implies that $\bar r=\overline{[g,h]}$ and thus, $r=C^\*\bar r=C^\*(\overline{[g,h]})=[g,h]$. This proves that the first diagram is a pushout.
\par Finally, let's prove that the diagram which expresses the naturality of $\alpha^\*$ is also a pushout. The first step is to lift $\alpha^\*$ to a morphism of $\P\^A$-algebras $\overline{\alpha^\*}$ so that $C^\*(\overline{\alpha^\*})=\alpha^\*$. Secondly, we are going to show that $\overline{\alpha^\*}$ is a coequalizer morphism from direct inspection, and finally, we use that $C^\*$ preserves the universality property of $\overline{\alpha^\*}$ to conclude our result.
\par Let's start by noticing that, since $B$ is a $\P\^A$-algebra it corresponds to a morphism of $\P$-algebras $\beta\colon A\to C^\*B$. Moreover, using the projection we obtain a morphism $A\xrightarrow{\beta}C^\*B\xrightarrow{p}\T C^\*B$ of $\P$-algebras which defines a new $\P\^A$-algebra $\overline{\T C^\*B}$. Concretely, this is the $\P\^A$-algebra defined over $\T C^\*B$ whose structure map is defined by:
\begin{align*}
&(\mu;a_1\,a_k|(x_1\,x_m)\=\mu_{\T C^\*B}(\beta(a_1)\,\beta(a_k),x_1\,x_m)
\end{align*}
Then, it is not hard to see that $\alpha^\*$ can be lifted to a morphism of $\P\^A$-algebras $\overline{\alpha^\*}\colon\overline{\T C^\*B}\to\T_AB$, which sends an element $y\in\T C^\*B$ to $\alpha^\*(y)\in\T_AB$. Recall also that, by construction, $\alpha^\*$ sends the generators $b$ and $\d b$ of $\T C^\*B$ to the corresponding generators $b$ and $\d_Ab$ of $C^\*\T_AB$.
\par By direct inspection we see that the $\P\^A$-algebra $\T_AB$ is generated by all $b\in B$ and by symbols $\d_Ab$ for $b\in B$, satisfying the following properties:
\begin{align*}
&(\mu;a_1\,a_k|_{\T_AB}(b_1\,b_m)=(\mu;a_1\,a_k|_B(b_1\,b_m)=\mu_{C^\*B}(\beta(a_1)\,\beta(a_k),b_1\,b_m)\\
&\d_A((\mu;a_1\,a_k|(b_1\,b_m))=\sum_{j=1}^m(\mu;a_1\,a_k|(b_1\,\d_Ab_j\,b_m)\\
&=\sum_{j=1}^m\mu_{C^\*\T_AB}(\beta(a_1)\,\beta(a_k),b_1\,\d_Ab_j\,b_m)
\end{align*}
Similarly, it is not hard to see that $\overline{\T C^\*B}$ is also generated by $b\in B$ and by symbols $\d b$, for $b\in B$, satisfying the following properties:
\begin{align*}
&(\mu;a_1\,a_k|_{\overline{\T C^\*B}}(b_1\,b_m)=\mu_{\T C^\*B}(\beta(a_1)\,\beta(a_k),b_1\,b_m)=\mu_{C^\*B}(\beta(a_1)\,\beta(a_k),b_1\,b_m)\\
&\d(\mu(b_1\,b_m))=\sum_{j=1}^m\mu(b_1\,\d b_j\,b_m)
\end{align*}
It is clear from this that the relations of $\T_AB$ imply the ones of $\overline{\T C^\*B}$. Since $\overline{\alpha^\*}$ sends generators to corresponding generators, this implies that $\T_AB$ can be represented as a quotient algebra of $\overline{\T C^\*B}$ over a specific ideal $I$, that is $\T_AB\cong\overline{\T C^\*B}/I$, and that $\overline{\alpha^\*}$ is the quotient map $\overline{\T C^\*B}\to\overline{\T C^\*B}/I$. Direct inspection shows that the ideal $I$ is generated by all the $\d_A(\beta(a))$ for every $a\in A$, that is in $\T_AB$, $\d_A(\beta(a))=0$.
\par Using a similar argument as the one we used to prove that the first diagram was a pushout, we conclude also that $\alpha^\*$ is a quotient map $\T C^\*B\to C^\*\T_AB$, so that $C^\*\T_AB$ is a quotient algebra of $\T C^\*B$ over an ideal $I$ generated by $\d_A(\beta(a))=0$.
\par Let's now come back to the naturality diagram and consider $g\colon\T C^\*E\to K$ and $h\colon C^\*\T_AB\to K$ as follows:
\begin{equation*}
\begin{tikzcd}
{\T C^\*B} & {C^\*\T_AB} \\
{\T C^\*E} & {C^\*\T_AE} \\
&& K
\arrow["{\alpha^*}", from=1-1, to=1-2]
\arrow["{\alpha^*}"', from=2-1, to=2-2]
\arrow["{\T C^\*f}"', from=1-1, to=2-1]
\arrow["{C^\*\T_Af}", from=1-2, to=2-2]
\arrow["g"', bend right, from=2-1, to=3-3]
\arrow["h", bend left, from=1-2, to=3-3]
\end{tikzcd}
\end{equation*}
This implies that:
\begin{align*}
&h(b)=g(f(b))\\
&h(\d_Ab)=g(f(\d b))=g(\d f(b))
\end{align*}
for every $b\in B$. Notice that, since $E$ is a $\P\^A$-algebra, we can also define a morphism of $\P$-algebras $\gamma\colon A\to E$ and that since $f$ is a morphism of $\P\^A$-algebras we have that $f(\beta(a))=\gamma(a)$. So, to lift $g$ to $C^\*\T_AE$ we need to show that $g(\d\gamma(a))=0$, however, we have the following:
\begin{eqnarray*}
& &g(\d\gamma(a))\\
&=&g(\d f(\beta(a))\\
&=&g(\d h(\d_A\beta(a)))\\
&=&0
\end{eqnarray*}
where we used that $\d_A\beta(a)=0$. This finally proves that we can lift $g$ to $\T C^\*E/I=C^\*\T_AE$, that is we find a morphism $\overline{[g,h]}\colon C^\*\T_AB\to K$. We leave to the reader to prove that such a morphism is the unique morphism which makes commuting the following diagram:
\begin{equation*}
\begin{tikzcd}
{\T C^\*B} & {C^\*\T_AB} \\
{\T C^\*E} & {C^\*\T_AE} \\
&& K
\arrow["{\alpha^*}", from=1-1, to=1-2]
\arrow["{\alpha^*}"', from=2-1, to=2-2]
\arrow["{\T C^\*f}"', from=1-1, to=2-1]
\arrow["{C^\*\T_Af}", from=1-2, to=2-2]
\arrow["g"', bend right, from=2-1, to=3-3]
\arrow["h", bend left, from=1-2, to=3-3]
\arrow["{[g,h]}"{description}, dashed, from=2-2, to=3-3]
\end{tikzcd}
\end{equation*}
This concludes the proof.
\end{proof}
\end{lemma}

We can finally prove the main result of this paper.

\begin{proposition}
\label{proposition:equivalence-slice-with-enveloping}
Consider the tangent morphism obtained as follows:
\begin{align*}
&\Geom^\*\o\Env\xrightarrow{(U;\eta)_{\Geom^\*\Env}}\Slice\o\Term\o\Geom^\*\o\Env\cong\\
&\cong\Slice\o\Geom^\*\o\Init\o\Env\xrightarrow{\Slice(\Geom^\*(C,\epsilon))}\Slice\o\Geom^\*
\end{align*}
This defines an equivalence of pseudofunctors which makes commutative the following diagram:
\begin{equation*}
\begin{tikzcd}
{\Operad^\op} & {\OprPair^\op} \\
\trmTngCat & \TngPair
\arrow["{\Geom^\*}"', from=1-1, to=2-1]
\arrow["{\Geom^\*}", from=1-2, to=2-2]
\arrow["\Env"', from=1-2, to=1-1]
\arrow["\Slice", from=2-2, to=2-1]
\end{tikzcd}
\end{equation*}
\begin{proof}
By Theorem~\ref{theorem:adjunction-Term-Slice}, $(U,\eta)$ is an equivalence of tangent categories. Moreover, thanks to Lemma~\ref{lemma:cartesianity-counit}, $\Geom^\*(C,\epsilon)$ is a Cartesian morphism of tangent pairs. By Lemma~\ref{lemma:cartesian-morphisms-become-strong}, $\Slice$ maps Cartesian morphisms into strong tangent morphisms. Thus, $\Slice(\Geom^\*(C,\epsilon))$ is strong. Finally, thanks to~\cite[Lemma 1.7]{moerdijk:enveloping-operads} the functorial component of $\Slice(\Geom^\*(C,\epsilon))$ is an isomorphism between the categories of $\P\^A$-algebras and the coslice category of $\P$-algebras under $A$, i.e. the slice category $\Alg_\P^\op/A$. Therefore, $\Slice(\Geom^\*(C,\epsilon))$ is an equivalence of tangent categories.
\end{proof}
\end{proposition}

\begin{theorem}
\label{theorem:geometric-tangent-category-enveloping-operad}
Given an operad $\P$ and a $\P$-algebra $A$, the geometric tangent category of the enveloping operad $\P\^A$ of $\P$ over $A$ is equivalent, as a tangent category, to the slice tangent category over $A$ of the geometric tangent category of $\P$. In formulas:
\begin{align*}
&\Geom(\P\^A)=\Geom(\P)/A
\end{align*}
\end{theorem}

Thanks to this characterization, we can now understand the vector fields over a $\P\^A$-algebra. For this purpose, recall that for a  morphism of $\P$-algebras $\beta\colon A\to B$ and a $B$-module $M$ (see Section~\ref{subsection:differential-bundles} for details) an $\beta$-relative derivation is a derivation $\delta\colon B\to M$, i.e. an $R$-linear morphism which satisfies the Leibniz rule:
\begin{align*}
&\delta(\mu(b_1\,b_m))=\sum_{k=1}^m\mu(b_1\,\delta(b_k)\,b_m)
\end{align*}
and moreover $\delta\o\beta=0$.

\begin{corollary}
\label{corollary:geometric-tangent-category-enveloping-operad}
For an operad $\P$, a $\P$-algebra $A$, and a $\P\^A$-algebra $B$, the vector fields over $B$ in the geometric tangent category of $\P\^A$ are in bijective correspondence with $\beta$-relative derivations, $\beta\colon A\to C^\*B$ being the morphism of $\P$-algebras corresponding to the $\P\^A$-algebra $B$.
\begin{proof}
Recall that in~\cite[Corollary~4.5.3]{ikonicoff:operadic-algebras-tagent-cats} it was proved that vector fields in a geometric tangent category of an operad correspond to derivations over the operadic algebras. Concretely, a vector field $v\colon\T A\to A$, regarded a morphism of $\P$-algebras, corresponds to a derivation $\delta_v\colon A\to A$ defined by:
\begin{align*}
&\delta_v(a)\=v(\d a)
\end{align*}
Viceversa, a derivation $\delta$ defines a vector field $v_\delta\colon\T A\to A$ by:
\begin{align*}
&v(a)\=a\\
&v(\d a)\=\delta(a)
\end{align*}
Thanks to Theorem~\ref{theorem:geometric-tangent-category-enveloping-operad}, we have that $\Geom(\P\^A)\cong\Geom(\P)/A$, thus, given a morphism $\beta\colon A\to C^\*B$, by definition of the slice tangent category, the tangent bundle functor $\T\^A$ of $\Geom(\P)/A$ is given by the coequalizer (in the category of $\P$-algebras):
\begin{equation*}
\begin{tikzcd}
{\T C^\*A} & {\T C^\*B} & {\T\^AB}
\arrow["\T\beta", shift left=2, from=1-1, to=1-2]
\arrow["{\T\beta zp}"', shift right=2, from=1-1, to=1-2]
\arrow["{v_\beta}", dashed, from=1-2, to=1-3]
\end{tikzcd}
\end{equation*}
or equivalently, by the pushout diagram:
\begin{equation*}
\begin{tikzcd}
{\T A} & {\T C^\*B} \\
A & {\T\^AB}
\arrow["{v_\beta}", from=1-2, to=2-2]
\arrow["\T\beta", from=1-1, to=1-2]
\arrow["z"', from=1-1, to=2-1]
\arrow["{\beta_\ast}"', from=2-1, to=2-2]
\arrow["\lrcorner"{anchor=center, pos=0.125, rotate=180}, draw=none, from=2-2, to=1-1]
\end{tikzcd}
\end{equation*}
This implies that $\T\^AB$ is the quotient of $\T C^\*B$ by the ideal generated by $\d\beta(a)$, for every $a\in A$. Therefore, a vector field $v\colon\T\^AB\to B$ corresponds to a derivation $\delta_v\colon B\to B$ defined by $\delta_v(b)\=v(\d b)$, and satisfying the following:
\begin{eqnarray*}
& &\delta_v(\beta(a))=v(\d\beta(a))=0
\end{eqnarray*}
that is a $\beta$-relative derivation of $B$. Conversely, a $\beta$-relative derivation $\delta\colon B\to B$ being a derivation over $B$, defines a vector field $v_\delta\colon\T C^\*B\to C^\*B$ over $C^\*B$ by $v_\delta(b)\=b$ and $v_\delta(\d b)=\delta(b)$, but since $\delta$ is $\beta$-relative, $v_\delta(\d\beta(a))=\delta(\beta(a))=0$, thus $v_\delta$ lifts to $\T\^AB\to B$.
\end{proof}
\end{corollary}

\subsection{The differential bundles of affine schemes over an operad}
\label{subsection:differential-bundles}
In~\cite[Section~4.6]{ikonicoff:operadic-algebras-tagent-cats}, the classification of differential objects for the geometric tangent category of an operad $\P$ was given. Roughly speaking, differential objects in a tangent category, first introduced by Cockett and Cruttwell in~\cite[Definition~4.8]{cockett:tangent-cats}, are the objects whose tangent bundle is trivial. In the tangent category of (connected) finite-dimensional smooth manifolds, differential objects correspond to the manifolds $\mathbb{R}^m$, for all integers $m$. For the geometric tangent category $\Geom(\P)$ of an operad $\P$, differential objects are in bijective correspondence with $\P(1)$-left modules, where we recall that $\P(1)$ becomes a unital and associative ring, once equipped with the unit and the composition of $\P$.
\par A related concept is the notion of differential bundles, introduced by Cockett and Cruttwell in~\cite[Definition~2.3]{cockett:differential-bundles}. Roughly speaking, differential bundles are bundles whose fibres are differential objects (cf.~\cite[Corollary~3.5]{cockett:differential-bundles}). More precisely, a differential bundle over $A\in\X$ in a tangent category $(\X,\TT)$ consists of a morphism $q\colon E\to A$ which admits pullbacks along any other morphism $B\to A$, together with a zero morphism $z_q\colon A\to E$, a sum morphism $s_q\colon E_2\to E$, $E_2$ being the pullback of $q$ along itself, and a vertical lift $l_q\colon E\to\T E$ satisfying a similar universality property of the vertical lift of the tangent structure $\TT$.
In~\cite{macadam:vector-bundles}, MacAdam proved that in the tangent category of finite-dimensional smooth manifolds, differential bundles are precisely vector bundles.\footnote{There is a slight difference between vector bundles and differential bundles over smooth manifolds. Vector bundles are defined as fibre bundles whose \textit{typical fibre} is a vector space. In general, differential bundles don't have a typical fibre and, when the manifold is not connected, they allow different connected components to have fibres with different dimensions. The two notions coincide for connected smooth manifolds.}
\par We also recall that a linear morphism $f\colon(q\colon E\to A)\to(q'\colon E'\to A)$ of differential bundles over $A\in\X$ is a morphism $f\colon E\to E'$ compatible with the lifts.
\par In this section, we are going to prove an important result: differential bundles over an operadic affine scheme $A$ in the geometric tangent category $\Geom(\P)$ of an operad $\P$ are equivalent to $A$-modules in the operadic sense. We recall that a module over a $\P$-algebra $A$ consists of an $R$-module $M$ equipped with a collection of $R$-linear morphisms $\P(m+1)\x A^{\x m}\x M\to M$ satisfying an equivariance condition with respect to the symmetric action, and associativity and unitality with respect to the structure map of $A$. We invite the interested reader to consult~\cite{loday:operads}, \cite{moerdijk:enveloping-operads} and~\cite{fresse:operads} for a detailed definition of modules over operadic algebras.
\par First, we prove that the correspondence between differential objects and left $\P(1)$-modules shown in~\cite[Theorem~4.6.8]{ikonicoff:operadic-algebras-tagent-cats} extends to a functorial equivalence between the category $\DObj_\lnr(\P)$ of differential objects and linear morphisms of the geometric tangent category $\Geom(\P)$ of $\P$ and the opposite of the category of left $\P(1)$-modules. We also prove that this equivalence is indeed an equivalence of tangent categories.
\par To understand what is the tangent structure over $\Mod_{\P(1)}^\op$, notice that, for any associative and unital $R$-algebra $A$, $\Mod_A$ is a semi-additive category, that is it has finite biproducts, denoted by $\oplus$. Thus, it comes with the canonical tangent structure $\FootTT_A$, whose tangent bundle functor is the diagonal functor $\FootT_AM\=M\oplus M$. It is straightforward to see that $\FootT_A$ is left-adjoint to itself, thus it also defines an adjoint tangent structure $\TT_A$ over the opposite category $\Mod_A^\op$.
\par It is interesting to note that (see~\cite[Example~4.6.10]{ikonicoff:operadic-algebras-tagent-cats}) given an associative and unital $R$-algebra $A$, the geometric tangent category of the associated operad $A^\bullet$ whose only non-trivial entry is $A^\bullet(1)=A$, is precisely $(\Mod_A^\op,\TT_A)$. So, in particular, $(\Mod_{\P(1)}^\op,\TT_{\P(1)})$ is the geometric tangent category of $\P(1)^\bullet$.
\par Note also that this construction extends to operadic algebras. To see that take into account a $\P$-algebra $A$ and let $\Env_\P(A)$ be its enveloping algebras. Concretely, the enveloping algebra $\Env_\P(A)$ of $A$ corresponds to the associative and unital algebra $\P\^A(1)$ which satisfies the following property: the category of modules over $A$ is equivalent to the category of left modules over $\Env_\P(A)$. Thus, let $A^\bullet$ be the operad whose only non-trivial entry is $A^\bullet(1)\=\Env_\P(A)$. Thus, $\Geom(A^\bullet)\cong(\Mod_{\Env_\P(A)}^\op,\TT_{\Env_\P(A)})\cong(\Mod_A^\op,\TT_A)$.

\begin{proposition}
\label{proposition:differential-objects-are-P(1)-modules}
For an operad $\P$, the tangent category $\DObj_\lnr(\P)$ of differential objects and linear morphisms of the geometric tangent category $\Geom(\P)$ of $\P$ is equivalent to the geometric tangent category $\Geom(\P(1)^\bullet)=(\Mod_{\P(1)}^\op,\TT_{\P(1)})$ associated with the associative and unital $R$-algebra $\P(1)$.
\begin{proof}
First, recall that $\P(1)$ is the enveloping algebra of the initial $\P$-algebra $\P(0)$ (cf.~\cite[Lemma~1.4]{moerdijk:enveloping-operads}), thus $\Mod_{\P(1)}\cong\Mod_{\P(0)}\^\P$. Recall also that in~\cite{ikonicoff:operadic-algebras-tagent-cats} it was proved the existence of a functor, for every $\P$-algebra $A$, $\Free_A\colon\Mod_A\to A/\Alg_\P$, which sends an $A$-module $M$, in the operadic sense, to a morphism of $\P$-algebras $A\to\Free_AM$. In particular, this defines a functor $\Free_A\colon\Mod_A\to\Alg_\P$, which maps each $M$ to $\Free_AM$. In~\cite[Theorem~4.6.8]{ikonicoff:operadic-algebras-tagent-cats} it was proved that, for every $\P(0)$-module $M$, $\Free_{\P(0)}M$ comes equipped with a canonical differential structure, so that $\Free_{\P(0)}M$ is a differential object of $\Geom(\P)$. Conversely, for a differential object $A\in\Geom(\P)$, there is a canonical vertical lift $l\colon\T A\to A$ (regarded as a morphism of $\P$-algebras). In particular, $l$ defines a derivation over $A$, as follows:
\begin{align*}
&\delta_l(a)\mapsto l(\d a)
\end{align*}
It was shown that the image of $\delta_l$ gives a $\P(0)$-module $UA$ and that the correspondence $M\mapsto\Free_{\P(0)}M$ and $A\to UA$ are inverses to each other, up to a canonical isomorphism. It is not hard to prove that this correspondence extends to a correspondence between linear morphisms. In particular, this means that given a $\P(0)$-linear morphism $f\colon M\to N$ of $\P(0)$-modules, $\Free_{\P(0)}f$ is again linear, in the sense that is compatible with the lifts of the corresponding differential objects. Similarly, given a linear morphism of differential objects $g\colon A\to B$ (regarded as a morphism of $\P$-algebras), define $Ug$ as the morphism whose domain is the image of $\delta_l$, $l$ being the vertical lift of $A$. So, for each $a\in A$, $Ug(\delta_l(a))\=g(l(\d a))$. However, since $g$ is compatible with the lifts, we have that $g(l(\d a))=l'(\d g(a))=\delta_{l'}(g(a))$, $l'$ being the vertical lift of $B$. Thus, $Ug$ is well-defined and also $\P(0)$-linear. Finally, this correspondence is functorial and it extends to an equivalence of categories. Finally, notice that the tangent structure over differential objects reduces to a Cartesian differential structure (cf.~\cite[Theorem~4.11]{cockett:tangent-cats}), thus the tangent bundle functor $\T$ sends a differential object $A$ to $A\times A$, being $\times$ the Cartesian product. Moreover, since all morphisms are linear, the same is true for morphisms as well, i.e. $\T f\cong f\times f$. However, Cartesian products in $\Geom(\P)$ are coproducts in $\Alg_\P$ and $\Free_{\P(0)}$ preserves coproducts, thus $\Free_{\P(0)}\T$ is isomorphic to the tangent bundle functor over $\Mod_{\P(0)}^\op$. Similarly, $\Free_{\P(0)}$ maps all the natural transformations of the tangent structure $\TT\^\P$ to the ones of $\TT_{\P(0)}$. This concludes the proof.
\end{proof}
\end{proposition}

Cockett and Cruttwell in~\cite[Proposition~5.12]{cockett:differential-bundles}) proved that differential bundles over an object $A\in\X$ of a tangent category $(\X,\TT)$ are equivalent to differential objects of the slice tangent category $(\X,\TT)/A$ of $(\X,\TT)$ over $A$. It is not hard to see that this correspondence extends to an equivalence of tangent categories:
\begin{align*}
&\DBnd(\X,\TT;A)\cong\DObj((\X,\TT)/A)
\end{align*}
between the tangent category $\DBnd(\X,\TT;A)$ of differential bundles over $A$ and the tangent category $\DObj((\X,\TT)/A)$ of differential objects of the slice tangent category $(\X,\TT)/A$. Moreover, this equivalence restricts to linear morphisms, that is:
\begin{align*}
&\DBnd_\lnr(\X,\TT;A)\cong\DObj_\lnr((\X,\TT)/A)
\end{align*}
where $\lnr$ indicates that morphisms are only linear morphisms (cf.~\cite[Definition~2.3]{cockett:differential-bundles}).
\par Let's denote by $\DBnd_\lnr(\P;A)$ the tangent category of differential bundles and linear morphisms over a $\P$-affine scheme $A$ in the geometric tangent category $\Geom(\P)$ of an operad $\P$.

\begin{theorem}
\label{theorem:classification-differential-bundles}
Let $\P$ be an operad and $A$ a $\P$-affine scheme. Then the tangent category $\DBnd_\lnr(\P;A)$ of differential bundles over $A$ and linear morphisms in the geometric tangent category of $\P$ is equivalent to the geometric tangent category of the operad $A^\bullet$:
\begin{align*}
&\DBnd_\lnr(\P;A)\cong\Geom(A^\bullet)\cong(\Mod_A^\op,\TT_A)
\end{align*}
In particular, differential bundles over $A$ are equivalent to $A$-modules in the operadic sense and linear morphisms of differential bundles over $A$ are equivalent to $A$-linear morphisms of $A$-modules (in the opposite of the category of $A$-modules).
\begin{proof}
Take into account an operad $\P$ and a $\P$-algebra $A$. Then, the tangent category $\DBnd_\lnr(\P;A)$ of differential bundles over $A$ and linear morphisms in the geometric tangent category $\Geom(\P)$ of $\P$ is equivalent to the tangent category $\DObj_\lnr(\Geom(\P)/A)$ of differential objects and linear morphisms of the slice tangent category $\Geom(\P)/A$. Thanks to Theorem~\ref{theorem:geometric-tangent-category-enveloping-operad}, $\Geom(\P)/A\cong\Geom(\P\^A)$, $\P\^A$ being the enveloping operad of $\P$ over $A$. By Proposition~\ref{proposition:differential-objects-are-P(1)-modules}, differential objects over $\Geom(\P\^A)$ are $\P\^A(1)$-left modules; in particular, $\DObj_\lnr(\Geom(\P\^A))\cong\Geom(\P\^A(1)^\bullet)$, but $\P\^A(1)$ is the enveloping algebra of $A$ (cf.~\cite[Definition~1.11]{moerdijk:enveloping-operads}), thus $\Geom(\P\^A(1)^\bullet)\cong\Geom(A^\bullet)$:
\begin{align*}
&\DBnd_\lnr(\P;A)                 &&&\\
=\quad&\DBnd_\lnr(\Geom(\P);A)         &&&\text{Diff. bundles are diff. objects in the slice tangent cat.}\\
\cong\quad&\DObj_\lnr(\Geom(\P)/A))    &&&\text{Theorem~\ref{theorem:geometric-tangent-category-enveloping-operad}}\\
\cong\quad&\DObj_\lnr(\Geom(\P\^A))    &&&\text{Proposition~\ref{proposition:differential-objects-are-P(1)-modules}}\\
\cong\quad&\Geom(\P\^A(1)^\bullet)     &&&\text{$\P\^A(1)=\Env_\P(A)$ (cf.~\cite[Definition~1.11]{moerdijk:enveloping-operads})}\\
\cong\quad&\Geom(A^\bullet)
\end{align*}
This concludes the proof.
\end{proof}
\end{theorem}

\section{Conclusion}
\label{section:conclusion}
The main results of this paper are the following:
\begin{description}
\item[Theorem~\ref{theorem:adjunction-Term-Slice}] In Section~\ref{subsection:universality-slicing} we gave a new characterization for the operation which takes a tangent pair $(\X,\TT;A)$ to its associated slice tangent category $(\X,\TT)/A$ in terms of the adjunction $\Term\dashv\Slice$.
\item[Theorem~\ref{theorem:geometric-tangent-category-enveloping-operad}] In Section~\ref{subsection:enveloping-tangent-category} we proved that the geometric tangent category of the enveloping operad of a $\P$-algebra $A$ is equivalent to the slice tangent category $\Geom(\P)/A$ of the geometric tangent category of $\P$ over $A$.
\item[Theorem~\ref{theorem:classification-differential-bundles}] In Section~\ref{subsection:differential-bundles} we classified differential bundles over operadic affine schemes in the geometric tangent category of an operad $\P$. We showed that differential bundles correspond to modules over the $\P$-algebras.
\end{description}
We also proved some minor but striking results:
\begin{description}
\item[Propositions~\ref{proposition:double-category-tangent-categories} and~\ref{proposition:conjuctions-tangent-morphisms}] We showed that tangent categories are organized in a double category whose horizontal and vertical morphisms are respectively lax and colax tangent morphisms. We also classified conjunctions in this double category in terms of a colax and a lax tangent morphism whose underlying functors form an adjunction and whose distributive laws are mates along this adjunction.
\item[Propositions~\ref{proposition:functoriality-Alg-star} and~\ref{proposition:functoriality-Alg-shriek}] We proved that the operation which takes an operad to its algebraic tangent category extends to a pair of functors $\Alg^\*$ and $\Alg_!$.
\item[Proposition~\ref{proposition:functoriality-Geom}] We proved that the operation which takes an operad to its geometric tangent category extends to a pair of functors $\Geom^\*$ and $\Geom_!$.
\item[Lemma~\ref{lemma:cartesian-morphisms-become-strong}] We proved that Cartesian morphisms of tangent pairs lift to the slice tangent categories as strong tangent morphisms.
\item[Corollary~\ref{corollary:geometric-tangent-category-enveloping-operad}] We classified vector fields over the geometric tangent category of the enveloping operad $\P\^A$ as relative derivations.
\end{description}

\subsection{Future work}
\label{subsection:future-work}
This paper is not just a natural continuation of \cite{ikonicoff:operadic-algebras-tagent-cats} but also the beginning of a fruitful program of research dedicated to understanding the intimate relationship between operads and the geometrical features of their corresponding operadic affine schemes. The classification of differential bundles, of vector fields (already covered in~\cite{ikonicoff:operadic-algebras-tagent-cats}), and the classification of the geometric tangent category of the enveloping operads represent the starting point of this program. Here is a list of some possible future directions of research:
\begin{enumerate}
\item Classifications of connections. Connections, introduced in~\cite{cockett:connections} are probably one of the most important geometrical tools available in a tangent category;
\item Classifications of principal bundles and principal connection. Principal bundles and principal connections were first introduced in the context of tangent categories by Cruttwell during a talk at Foundational Methods in Computer Science 2017 (General connections in tangent categories - FMCS 2017);
\item Study of sector forms and cohomology for operadic affine schemes (cf.~\cite{cruttwell:sector-cohomology});
\item Study of curve objects and of differential equations for operadic affine schemes (cf.~\cite{cockett:differential-equations});
\item An important application of this program is the study of associative affine schemes, which lead to a description of non-commutative algebraic geometry via tangent categories;
\item An important construction in the theory of operads is Kozsul duality (cf.~\cite[Chapter~7]{loday:operads}). A natural question is what is the geometric tangent category of the Kozsul dual of an operad $\P$;
\item There is a notion of distributive laws between operads which allows two operads to be composed together. An example of an operad obtained via a distributive law between the operad $\Ass$ and $\Lie$ is the operad $\Pois$, whose algebras are Poisson algebras. What kind of relationship exists between the geometric tangent categories of two operads $\P$ and $\OprQ$ and the one of the operad obtained by composing $\P$ and $\OprQ$ provided a distributive law between them?
\end{enumerate}
These are only a few of the possible new paths of research that this paper inspires.

\printbibliography[title=Bibliography]

\end{document}